\documentclass{dcds-sOF}
\usepackage{amsmath}
  \usepackage{paralist}
  \usepackage{graphics} 
  \usepackage{epsfig} 
\usepackage{graphicx}  \usepackage{epstopdf} 
 \usepackage[colorlinks=true]{hyperref}
\hypersetup{urlcolor=blue, citecolor=red}

\def\currentvolume{X}
 \def\currentissue{X}
  \def\currentyear{200X}
   \def\currentmonth{XX}
    \def\ppages{X--XX}
    \def\DOI{ }

\theoremstyle{definition}

\usepackage{comment}
\usepackage{bm}
\usepackage{amsmath}
\usepackage{amssymb}
\usepackage{caption}
\usepackage{subcaption}
\usepackage{graphicx}
\usepackage{float}
\usepackage{cancel}
\usepackage{empheq,etoolbox}
\usepackage{algorithm}
\usepackage{algorithmic}
\usepackage{hyperref}
\usepackage{todonotes}
\usepackage{multirow}
\usepackage{empheq}
\usepackage[pagewise]{lineno}
\usepackage{longtable}

\newcommand{\xhat}{{\widehat{\bm x}}}

\title[Geometric multiscale model for the blood flow in the left heart] 
      {A geometric multiscale model for the numerical simulation of blood flow in the human\\ left heart} 

\author[Alberto Zingaro, Ivan Fumagalli, Luca Dede' et al.]{}


 \email{alberto.zingaro@polimi.it}
 \email{ivan.fumagalli@polimi.it}
 \email{luca.dede@polimi.it}
 \email{marco.fedele@polimi.it}
 \email{pasqualeclaudio.africa@polimi.it}
 \email{tonycorno2@gmail.com}
 \email{alfio.quarteroni@polimi.it}
 
\thanks{$^*$ Corresponding author: Alberto Zingaro}

\begin{document}
\maketitle

\centerline{\scshape Alberto Zingaro$^{*, 1}$,  Ivan Fumagalli$^1$,  Luca Dede'$^1$}
\centerline{\scshape Marco Fedele$^1$, Pasquale C. Africa$^1$}
\centerline{\scshape Antonio F. Corno$^{2}$ and Alfio Quarteroni$^{1, 3}$}
\medskip
{\footnotesize
 \centerline{$^1$ MOX, Dipartimento di Matematica, Politecnico di Milano,	Piazza Leonardo da Vinci, 32}
 \centerline{20133 Milan, Italy}
 \centerline{$^2$ Children’s Heart Institute, Hermann Children’s Hospital, University of Texas Health}
 \centerline{McGovern Medical School, Houston, TX, USA}
 \centerline{$^3$ Institute of Mathematics, \'Ecole Polytechnique F\'ed\'erale de Lausanne}
 \centerline{Station 8, Av. Piccard}
 \centerline{CH-1015 Lausanne, Switzerland  (Professor Emeritus)}
}



\begin{abstract}
 We present a new computational model for the numerical simulation of blood flow in the human left heart. To this aim, we use the Navier-Stokes equations in an Arbitrary Lagrangian Eulerian formulation to account for the endocardium motion and we model the cardiac valves by means of the Resistive Immersed Implicit Surface method. To impose a physiological displacement of the domain boundary, we use a 3D cardiac electromechanical model of the left ventricle coupled to a lumped-parameter (0D) closed-loop model of the remaining circulation. We thus obtain a one-way coupled electromechanics-fluid dynamics model in the left ventricle. To extend the left ventricle motion to the endocardium of the left atrium and to that of the ascending aorta, we introduce a preprocessing procedure according to which an harmonic extension of the left ventricle displacement is combined with the motion of the left atrium based on the 0D model. To better match the 3D cardiac fluid flow with the external blood circulation, we couple the 3D Navier-Stokes equations to the 0D circulation model, obtaining a multiscale coupled 3D-0D fluid dynamics model that we solve via a segregated numerical scheme. We carry out numerical simulations for a healthy left heart and we validate our model by showing that meaningful hemodynamic indicators are correctly reproduced.
\end{abstract}
\vspace{1cm}
\textbf{List of acronyms}
\begin{longtable}{ll}
        AA & Ascending Aorta \\
        ALE & Arbitrary Lagrangian Eulerian \\
        AV & Aortic Valve \\
        CFD & Computational Fluid Dynamics \\
        EM & Electromechanics \\
        FE & Finite Element \\
        FSI & Fluid Structure Interaction \\
        GMRES & Generalized Minimal Residual \\
        LA & Left Atrium \\
        LES & Large Eddy Simulation \\
        LH & Left Heart \\
        LV & Left Ventricle \\
        MV & Mitral Valve \\
        NS & Navier-Stokes \\
        PV & Pulmonary Valve \\
        RA & Right Atrium \\
        RIIS & Resistive Immersed Implicit Surface \\
        RV & Right Ventricle \\
        TV & Tricuspid Valve \\
        VMS & Variational Multiscale \\
        WSS & Wall shear Stress
\end{longtable}

\section{Introduction}
\label{SEC:INTRO}
The study of cardiac blood flow aims at improving the knowledge of the heart physiology, assessing the pathological conditions and potentially helping therapeutical treatment and the design of surgical interventions. In the clinical routine, blood flow analysis is conventionally based on non-invasive imaging techniques. However, the space and time resolution of the available techniques is not accurate enough to capture small-scales features like recirculation regions, possible regions of transition to turbulence and small coherent structures \cite{ngo2019four}. Moreover, imaging techniques cannot accurately provide relevant fluid dynamics indicators such as wall shear stress (WSS), turbulent kinetic energy dissipation, or the oscillatory stress index, which are correlated with the function and remodeling of the heart and vessels \cite{quarteroni2017cardiovascular, Chnafa_2014, ngo2019four, zhao2003comparative, hollnagel2009comparative}. In this respect, numerical simulations -- also known as in silico models -- of the heart and circulation represent a valuable tool to quantitatively assess the cardiac function and to enhance the understanding of cardiac dysfunction.

The hemodynamics in the heart chambers is characterized by different complex features that computational modeling needs to take into account \cite{Chnafa_2014}. The mathematical problem is defined in complex geometries, and there is a strong interaction between the myocardial structure and the blood flow due to the electromechanical activity of the heart. This yields a complex coupled problem among electrophysiology, mechanics and fluid dynamics. Furthermore, the topology of the domain changes during the heartbeat due to the presence of unidirectional cardiac valves, affecting the dynamics of the intracardiac flow. In addition, the blood flow regime is known to be neither laminar, nor fully turbulent, but rather transitional \cite{bluestein1994transition, verkaik2012coupled, ZDMQ_2020, vignon2010outflow}. Eventually, the flow in the heart is strictly coupled with the flow in the pulmonary and systemic circulation.

The aim of this work is to introduce an accurate computational model accounting for all the features of the hemodynamics in the left heart (LH): the motion of the surrounding cardiac tissue, the dynamics of aortic valve (AV) and mitral valve (MV), transitional flow effects, and the flows and pressures in the  rest of the cardiocirculatory system. Moreover, the present study represents a significant step towards a high fidelity fluid dynamics model of the whole human heart.

A key point in 3D hemodynamic models is the treatment of the boundary displacement, especially in the cardiac chambers, were the blood motion is driven by the heart contraction and relaxation. Boundary displacement can be mainly modeled  according to the following two paradigms: Fluid-Structure-Interaction (FSI) models and Computational Fluid Dynamics (CFD) simulations under a prescribed wall motion. In cardiac FSI problems, the fluid model is coupled with an electromechanics (EM) model and the coupling is explicitly solved, entailing a high computational effort \cite{viola2021fsei, viola2020fluid, karabelas2018towards, viola2021effects, choi2015new, gerbi2018numerical, bucelli2021partitioned, santiago2018fully, santiago2018fluid}. In prescribed displacement CFD simulations, on the other hand, the influence of the cardiac walls motion is modeled by analytical laws \cite{TagliabueLV,DMQ_2019, tagliabue2017complex, ZDMQ_2020, domenichini2005three, baccani2002vortex, zheng2012computational, seo2013effect}, by patient-specific image-based reconstructions \cite{Fumagalli2020, This_regurgitation, Chnafa_2014, MasciLA, masci2017patient, masci2020proof}, or via the displacement field computed from a previous EM simulation \cite{augustin2016patient, karabelas2018towards, fernandez_2019}.

A possible approach to account for the coupling between the flow field of the region of interest and the one of the remaining circulation is the geometric multiscale modeling \cite{quarteroni2016geometric}. The specific region of interest (a vessel or one of the heart chambers) is represented by a 3D model, while the remaining part of the circulation is simulated by models featuring a lower geometric dimension, as 0D \cite{Blanco_2010, shi2006numerical, quarteroni2016geometric, milivsic2004analysis, kim2009coupling} or 1D \cite{quarteroni2016geometric, van2011pulse, formaggia2003one, formaggia2001coupling} models. With the term 0D models, we refer to lumped-parameter models where the dependence on the spatial coordinates is completely neglected, and a uniform spatial distribution of pressures and flowrates in any specific compartment is assumed \cite{quarteroni2016geometric, shi2011review}.

We model the blood flow by the incompressible Navier-Stokes (NS) equations in an Arbitrary Lagrangian Eulerian (ALE) formulation to account for the moving boundary, whereas the AV and MV are immersed in the domain by the Resistive Immersed Implicit Surface (RIIS) method \cite{Fedele_RIIS, Fumagalli2020}.
In order to account for the transitional flow regime, the Variational Multiscale - Large Eddy Simulation (VMS-LES) turbulence model is considered \cite{bazilevs_2007,DF_2015, ZDMQ_2020}. This model has the advantage of acting also as a stabilization method for the CFD numerical scheme. The motion of the wall is derived from an EM simulation on the left ventricle (LV)  \cite{REGAZZONI2022111083, REGAZZONI2022111083} and then extended to the whole boundary of the domain of interest by means of an original preprocessing procedure, suitably considering a volume-based definition of the displacement of the left atrium (LA).
By prescribing the EM-based velocity at the endocardial wall, we enforce a one-way (kinematic) coupling condition between EM and CFD in the LV. Furthermore, to address the interdependence between the fluid dynamics of the LH and the remaining cardiovascular system, we couple the 3D CFD model to the lumped-parameter (0D) model proposed in  \cite{REGAZZONI2022111083} representing the whole cardiovascular system.

We numerically simulate the hemodynamics of the LH in physiological conditions. We analyze the obtained complex blood flow pattern, and we validate our model by comparing meaningful hemodynamic indicators with data available in the literature.

This paper is organized as follows: in Section~\ref{SEC:0D_CIRCULATION}, we briefly recall the lumped-parameter circulation model of the cardiovascular system we employ; in Section~\ref{SEC:3D_CFD}, we describe the mathematical model for the moving-domain hemodynamics in the LH with the immersed valves. In particular, we introduce the preprocessing procedure to reconstruct the displacement on the whole domain boundary based on an EM simulation of the LV in Section~\ref{SEC:displacement_modeling}, while we present the reduced valve dynamics model in Section~\ref{SEC:valves_modeling}. The coupling between the 0D circulation model and the 3D CFD model is described in Section~ \ref{SEC:coupling_0D_3D}.
Section~\ref{SEC:numerical_methods} is devoted to the description of the numerical approximation of each component of the model and the segregated scheme we propose for their coupling.
Numerical results are reported in Section~\ref{SEC:numerical_results} and conclusions are drawn in Section~\ref{SEC:conclusions}.

\section{A 0D circulation model of the whole cardiovascular system}
\label{SEC:0D_CIRCULATION}
\begin{figure}[t]
	\centering
	\includegraphics[trim={4 4 4 4},
	clip,
	width=0.7\textwidth]
	{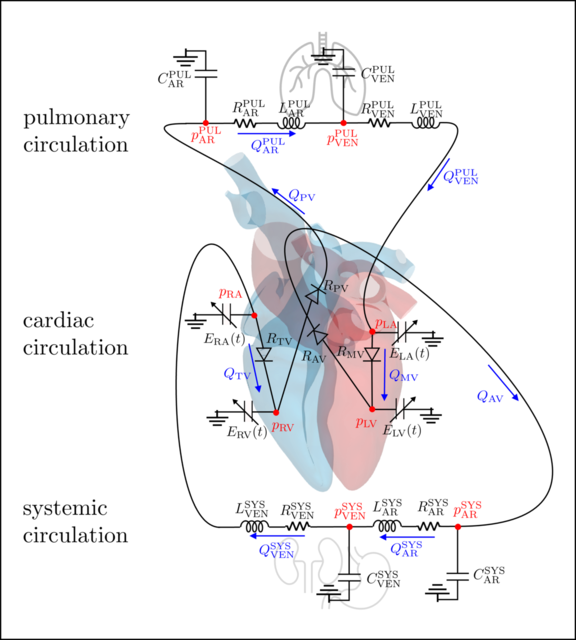}
	\caption{The 0D circulation model.}
	\label{circuit_0D}
\end{figure}

We employ the closed-loop lumped-parameter (0D) circulation model that was proposed in \cite{REGAZZONI2022111083} (see also \cite{Blanco_2010, Hirschvogel_2017}).
As represented in Figure \ref{circuit_0D}, in this model the systemic and pulmonary circulations are represented by resistance-inductance-capacitance (RLC) circuits, one for the arterial part and one for the venous part.
The resistance $R$ models viscosity effects, the inductance $L$ inertial effects, and the conductance $C$ the compliance of vessels.
The four heart chambers (cardiac circulation, Figure~\ref{circuit_0D}) are modeled by elements with time-varying elastance $E(t)$; the four valves by non-ideal diodes.
Let $V_\mathrm{i}(t)$, $p_\mathrm{i}(t)$ and $E_\mathrm{i}(t)$, with $\mathrm i =\mathrm{LA, \, LV, \, RA, \, RV}$ be respectively the volume, pressure and elastances of the four heart chambers: left atrium (LA), left ventricle (LV), right atrium (RA), right ventricle (RV); $Q_\text{j}(t)$, with $\mathrm j= \text{MV, AV, TV, PV}$ the flowrates through the mitral valve (MV), aortic valve (AV), tricuspid valve (TV) and pulmonary valve (PV). $Q_\text{AR}^\text{SYS}(t)$,  $Q_\text{VEN}^\text{SYS}(t)$ and $p_\text{AR}^\text{SYS}(t)$,  $p_\text{VEN}^\text{SYS}(t)$ are the flowrates and pressures in the systemic circulation in the arterial and venous parts, respectively; similarly for $Q_\text{AR}^\text{PUL}(t)$,  $Q_\text{VEN}^\text{PUL}(t)$ and  $p_\text{AR}^\text{PUL}(t)$,  $p_\text{VEN}^\text{PUL}(t)$ in the pulmonary circulation.
The 0D closed-loop circulation model of the whole cardiovascular system reads \cite{REGAZZONI2022111083}: for any $t \in (0, T_f)$:
\begin{subequations}
	\allowdisplaybreaks
	\begin{align}
		\allowdisplaybreaks
		\frac{\mathrm d V_\text{LA}(t)}{\mathrm d t} &  =  Q_\text{VEN}^\text{PUL}(t) - Q_\text{MV}(t),
		\label{0D_VLA}  \\
		\frac{\mathrm dV_\text{LV}(t)}{\mathrm  dt} &  =  Q_\text{MV}(t) - Q_\text{AV}(t),
		\label{0D_VLV} \\
		\frac{\mathrm  dV_\text{RA}(t)}{\mathrm  dt} &  =  Q_\text{VEN}^\text{SYS}(t) - Q_\text{TV}(t) ,
		\label{0D_VRA}  \\
		\frac{\mathrm  dV_\text{RV}(t)}{\mathrm dt} &  =  Q_\text{TV}(t) - Q_\text{PV}(t),
		\label{0D_VRV}  \\
		C_\text{AR}^\text{SYS}\frac{\mathrm dp_\text{AR}^\text{SYS}(t)}{\mathrm dt} & = Q_\text{AV}(t) - Q_\text{AR}^\text{SYS}(t),
		\label{0D_pARSYS} \\
		C_\text{VEN}^\text{SYS}\frac{\mathrm dp_\text{VEN}^\text{SYS}(t)}{\mathrm dt} & = Q_\text{AR}^\text{SYS}(t) - Q_\text{VEN}^\text{SYS}(t),
		\label{0D_pVENSYS} \\
		C_\text{AR}^\text{PUL}\frac{\mathrm dp_\text{AR}^\text{PUL}(t)}{\mathrm dt} & = Q_\text{PV}(t) - Q_\text{AR}^\text{PUL}(t),
		\label{0D_pARPUL} \\
		C_\text{VEN}^\text{PUL}\frac{\mathrm dp_\text{VEN}^\text{PUL}(t)}{\mathrm dt} & = Q_\text{AR}^\text{PUL}(t) - Q_\text{VEN}^\text{PUL}(t),
		\label{0D_pVENPUL} \\
		\frac{L_\text{AR}^\text{SYS}}{R_\text{AR}^\text{SYS}} \frac{\mathrm dQ_\text{AR}^\text{SYS}(t)}{\mathrm dt} & = - Q_\text{AR}^\text{SYS}(t) - \frac{p_\text{VEN}^\text{SYS}(t)-p_\text{AR}^\text{SYS}(t)}{R_\text{AR}^\text{SYS}},
		\label{OD_QARSYS} \\
		\frac{L_\text{VEN}^\text{SYS}}{R_\text{VEN}^\text{SYS}} \frac{\mathrm dQ_\text{VEN}^\text{SYS}(t)}{\mathrm dt} & = - Q_\text{VEN}^\text{SYS}(t) - \frac{p_\text{RA}(t)-p_\text{VEN}^\text{SYS}(t)}{R_\text{VEN}^\text{SYS}},
		\label{OD_QVENSYS} \\
		\frac{L_\text{AR}^\text{PUL}}{R_\text{AR}^\text{PUL}} \frac{\mathrm dQ_\text{AR}^\text{PUL}(t)}{\mathrm dt} & = - Q_\text{AR}^\text{PUL}(t) - \frac{p_\text{VEN}^\text{PUL}(t)-p_\text{AR}^\text{PUL}(t)}{R_\text{AR}^\text{PUL}},
		\label{OD_QARPUL} \\
		\frac{L_\text{VEN}^\text{PUL}}{R_\text{VEN}^\text{PUL}} \frac{\mathrm dQ_\text{VEN}^\text{PUL}(t)}{\mathrm dt} & = - Q_\text{VEN}^\text{PUL}(t) - \frac{p_\text{LA}(t)-p_\text{VEN}^\text{PUL}(t)}{R_\text{VEN}^\text{PUL}},
		\label{0D_QVENPUL}
	\end{align}
	\label{0D_ODE}
\end{subequations}
where, by denoting with $p_\mathrm{EX}$ the external pressure and $V_{0, \mathrm{i}}$ the resting volume of each cardiac chamber, the atrial and ventricular pressures are defined as
\begin{equation}
	p_{\mathrm i}(t) = p_\text{EX}(t) + E_\text{i}(t) \left ( V_\text{i}(t) - V_{0, \text{i}}\right ), \,  \mathrm{with} \; \mathrm{i} = \mathrm{LA},\, \mathrm{LV},\, \mathrm{RA}, \, \mathrm{RV},
	\label{0D_pressurechambers}
\end{equation}
and the flowrates across the valves as
\begin{equation}
	Q_{\mathrm{j}}(t) = \frac{p_{\mathrm{up, j}}(t) - p_{\mathrm{down, j}}(t)}{R_\mathrm{j}(p_{\mathrm{up, j}}(t), p_{\mathrm{down, j}}(t))}, \, \mathrm{with} \; \mathrm{j} = \mathrm{MV},\, \mathrm{AV},\, \mathrm{TV}, \, \mathrm{PV}.
	\label{0D_flowrates}
\end{equation}
In Eq.~\eqref{0D_flowrates}, $R_\mathrm{j}(p_{\mathrm{up, j}}(t), p_{\mathrm{down, j}}(t))$ is the valve resistance, being $p_{\mathrm{up, j}}, \, p_{\mathrm{down, j}}$ the pressure upstream and downstream the valve $\mathrm j$ respectively, with $\mathrm j = \mathrm{MV},\, \mathrm{AV},\, \mathrm{TV},\, \mathrm{PV}$ (see Figure~\ref{circuit_0D}).

\section{A 3D fluid dynamics model of the left heart}
\label{SEC:3D_CFD}
In this section, we introduce the 3D fluid dynamics model of the LH. Specifically, we present the NS equations in ALE frawework with RIIS modelling in Section~\ref{SEC:NS_ALE_RIIS}; the LH model and the boundary conditions (BCs) we employ are described in Section~\ref{SEC:LH_geometry_model_BCs}. We introduce the preprocessing procedure to compute the LH displacement in Section~\ref{SEC:displacement_modeling}, and the valves dynamics in Section~\ref{SEC:valves_modeling}.

\subsection{The NS-ALE-RIIS equations}
\label{SEC:NS_ALE_RIIS}
\begin{figure}[t!]
	\centering
	\includegraphics[trim={2 4cm 2 2},
	clip,
	width=\textwidth]
	{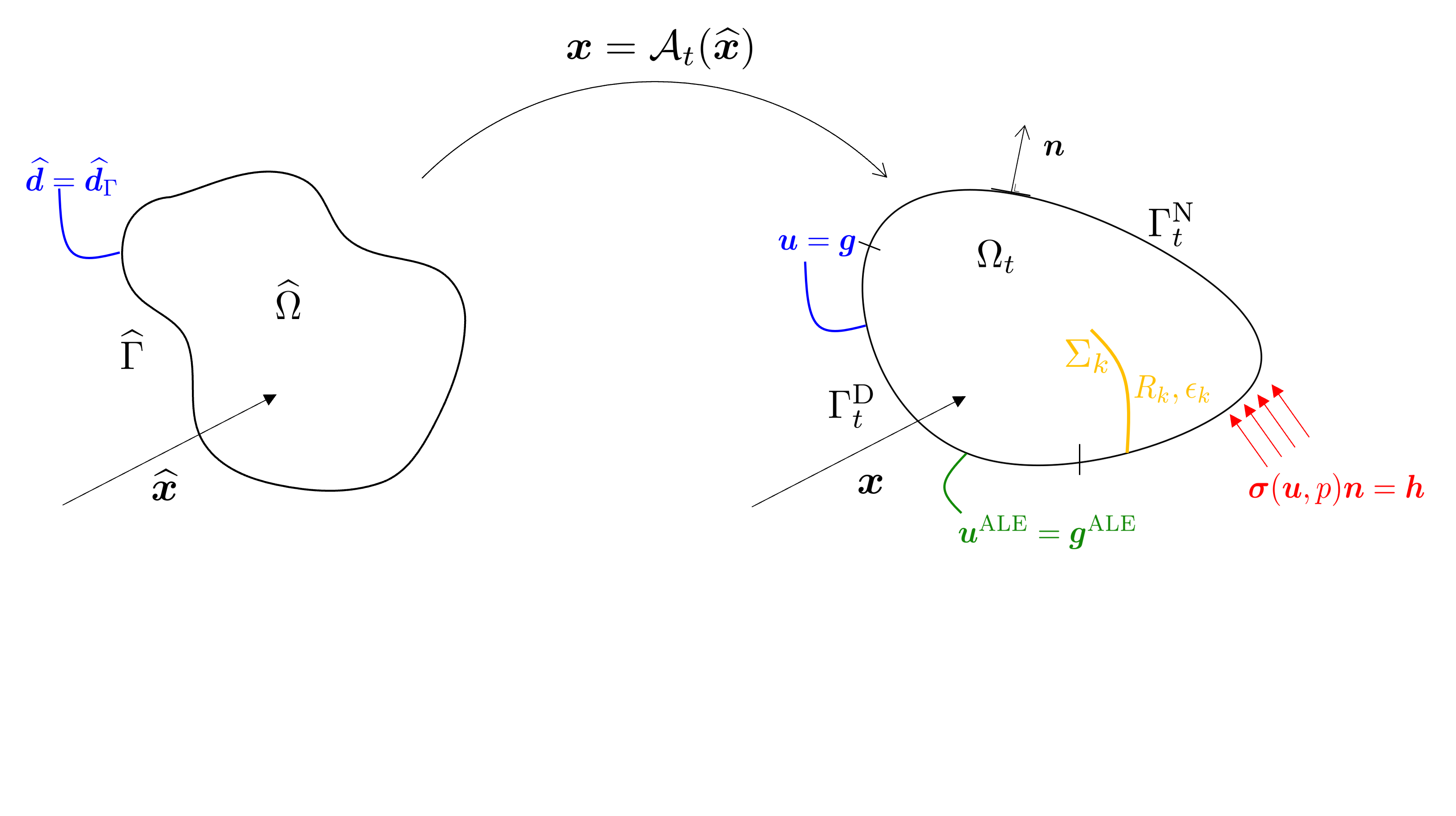}
	\caption{Fluid domain in reference configuration (left), ALE map $\bm x = \mathcal{A}_t (\xhat)$ and domain in current configuration (right). In the current configuration, the domain $\Omega_t$ is bounded by $\Gamma_t = \Gamma_t^\mathrm{D} \cup \Gamma_t^\mathrm{N}$; $\Sigma_k$ is an immersed surface modeled by means of the RIIS method. }
	\label{domain_ns_ale_riis}
\end{figure}

Let $\Omega_t \subset \mathbb{R}^d$ be the fluid domain at a specific time instant $t>0$ (current configuration), with a sufficiently regular boundary  $ \partial \Omega_t \equiv \Gamma_t $ oriented by the outward pointing normal unit vector ${\bm {n}}$. We denote as $\Gamma_{t}^\mathrm D$ and $\Gamma_{t}^\mathrm N$ the portions of the boundary where, respectively, Dirichlet and Neumann type BCs are prescribed, with $\Gamma_t = \overline{\Gamma_t^\mathrm D} \cup \overline{\Gamma_t^\mathrm N}$ and $\overset{\circ}{\Gamma_t^\mathrm D} \cap \overset{\circ}{\Gamma_t^\mathrm N} = \emptyset$.
Let $\widehat{\Omega} \subset \mathbb{R}^d $ be the reference configuration of $\Omega_t$, and $\partial \widehat \Omega \equiv \widehat \Gamma$ its sufficiently regular boundary, as we display in Figure \ref{domain_ns_ale_riis}.

We introduce the ALE map $\mathcal{A}_t$ which associates, at each $t \in (0, T_f)$, a point $\widehat{\bm x} \in \widehat{\Omega}$ to a point $\bm x \in \Omega_t$ \cite{Formaggia_Nobile_1999}:
\begin{equation}
	\mathcal A_t: \widehat{\Omega} \to \Omega_t: \quad \bm x = \mathcal{A}_t(\widehat{\bm x}) = \widehat {\bm x} + \widehat {\bm d}(\widehat {\bm x}, t),
	\label{ALEmap}\end{equation}
being $\widehat{\bm d}$ the domain displacement with respect to the reference configuration $\widehat{\Omega}$. For every function $w$ defined in the current configuration, we denote by $\widehat{w} = w \circ \mathcal{A}_t$ the corresponding function in the reference frame. Similarly, $w = \widehat{w} \circ \mathcal{A}_t^{-1}$. Assuming $\widehat{\bm{d}}_{\Gamma} (\widehat{\bm x}, t)$ as known on the whole boundary $ \widehat{\Gamma}$ at any time $t \in (0, T_f)$, we can prolongate it to the fluid domain by solving, at each time, the following harmonic extension problem:
\begin{subequations}
	\begin{empheq}[left=\empheqlbrace]{align}
		- \Delta  \widehat{\bm{d}} & =  \,  \bm 0 & \text{ in } \widehat{\Omega},  \label{ALE_laplacian}\\
		\widehat{\bm d} & =  \, \widehat{\bm{d}}_{\Gamma} & \text{ on } \widehat{\Gamma}. \label{ALE_laplacian_BC}
	\end{empheq}
	\label{geometric_strong}
\end{subequations}
Eq.~\eqref{geometric_strong}, together with Eq.~\eqref{ALEmap},  allows to compute the current domain $\Omega_t = \mathcal A_t(\widehat \Omega)$, for all $t \in (0, T_f)$. The ALE velocity is eventually obtained as
\begin{equation}
	\bm u^\text{ALE} = \left ( \frac{\partial \widehat{\bm d}}{\partial t}\right) \circ \mathcal{A}_t^{-1}. \label{uALE_dd_dt}
\end{equation}

In the heart chambers, blood can be regarded as a Newtonian, incompressible fluid and modeled by means of the incompressible Navier-Stokes (NS) equations \cite{quarteroni2017cardiovascular, perktold1994flow, taylor1996finite, taylor1998finite}. Let $\bm u$ be the fluid velocity, $p$ the pressure. The motion of the domain is accounted for by expressing the NS equations in an ALE framework \cite{Formaggia_Nobile_1999, duarte2004arbitrary}. The action of the cardiac valves on the fluid field is simulated thanks to the RIIS method  \cite{Fedele_RIIS, Fumagalli2020}.
It lays in the class of immersed boundary - fictitious domain methods, and it represents a moving immersed surface in an Eulerian framework, without requiring the use of a surface-conforming mesh. Specifically, we introduce an additional resistive term into the momentum balance of the NS equations, penalizing a kinematic condition, i.e. the adherence of the blood to $m$ immersed moving surfaces $\Sigma_k$:
\begin{equation}
	\Sigma_k = \left \{ \bm x \, : \, \varphi_k (\bm x) = 0 \right \}, \quad \text{with } k = 1, \dots, m.
	\label{riis_sigma_k}
\end{equation}
In Eq. \eqref{riis_sigma_k}, $\varphi_k(\bm x)$ is a signed-distance function that implicitly describes the $k$--th immersed surface. In particular, we denote by $R_k$ the resistance coefficient; the  resistive term has support in a narrow layer around $\Sigma_k$, represented by a smoothed Dirac delta type function:
\begin{equation*}
	\delta_{\Sigma_k, \varepsilon_k} (\varphi_k(\bm x)) =
	\begin{cases}
		\dfrac{1 + \cos(\pi \varphi_k(\bm x) / \varepsilon_k)}{2 \varepsilon_k} & \text{ if } |\varphi_k(\bm x)| \leq \varepsilon_k, \\
		0 & \text{ if } |\varphi_k(\bm x)| > \varepsilon_k, \\
	\end{cases}
	\label{delta_dirac_RIIS}
\end{equation*}
and $\varepsilon_k > 0$ is a suitable parameter representing half of the thickness of the leaflet, $\text{with } k = 1, \dots, m$. 

The incompressible NS equations in ALE framework endowed with the RIIS method read: find $\bm u, \, p$ such that:
\begin{subequations}
	\begin{empheq}[left=\empheqlbrace]{align}
		& \rho \frac{\widehat \partial \bm u}{\partial t} + \rho \left( \left( \bm{u} - \bm{u}^{\text{ALE}} \right) \cdot \nabla \right) \bm{u} -  \nabla \cdot \bm{\sigma} (\bm u, p) & \notag
		\\
		& + \sum_{k=1}^{m}\frac{R_k}{\varepsilon_k}\delta_{\Sigma_k, \varepsilon_k}(\varphi_k)\left(\bm u - \bm u ^\mathrm{ALE}\right )  =  \bm 0 & \quad \text{ in } \Omega_t \times (0, T_f),  \label{eq_ns}
		\\
		& \nabla \cdot \bm u  =   0  &  \quad \text{ in } \Omega_t \times (0, T_f), \label{eq_div}
		\\
		& \bm u  = \bm g & \quad \text{ on } \Gamma_t^\mathrm{D} \times (0, T_f), \label{eq_dirichletbc} \\
		& \bm \sigma (\bm u, p) {\bm n}  = \bm h & \quad \text{ on } \Gamma_t^\mathrm{N} \times (0, T_f), \label{eq_neumannbc} \\
		& \bm u  = \bm u_0 & \quad \text{ in } \Omega_0 \times \{0\}.
	\end{empheq}
	\label{NS_ALE_RIIS_strong}
\end{subequations}
We will refer to this formulation as NS-ALE-RIIS formulation; in particular, $\frac{\widehat \partial \bm u }{\partial t} =  \frac{\partial \bm u }{\partial t} + (\bm u^{\text{ALE}} \cdot \nabla ) \bm u  $ is the ALE derivative, $\rho$ the fluid density and $\bm{\sigma} (\bm u, p)$ the stress tensor defined for Newtonian, incompressible and viscous fluids as
$
\bm \sigma (\bm u, p) = -p \bm I + 2 \mu \bm \epsilon (\bm u)
$.
$\mu$ is the dynamic viscosity and $
\bm \epsilon (\bm u) = \frac{1}{2} \left ( \nabla \bm u + \left ( \nabla \bm u \right )^T\right)$ the strain-rate tensor. The functions $\bm g$ and $\bm h$ are the Dirichlet and Neumann data that will be discussed in the next section, while  $\bm u_0$  the initial velocity.

\subsection{The LH geometry model and boundary conditions}
\label{SEC:LH_geometry_model_BCs}

\begin{figure}[t]
	\centering
	\begin{subfigure}{.32\textwidth}
		\centering
		\includegraphics[trim={2 2 2 2},clip, width=0.95\textwidth]{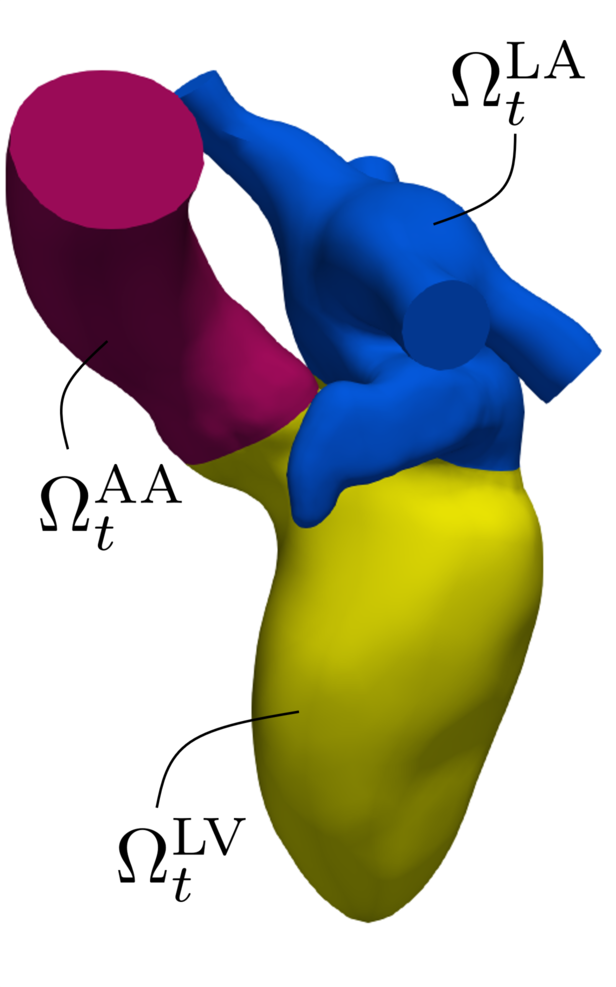}
		\caption{ }
		\label{LH_volumetags}
	\end{subfigure}
	\begin{subfigure}{.32\textwidth}
		\includegraphics[trim={2 2 2 2 },clip,width=0.95\textwidth]{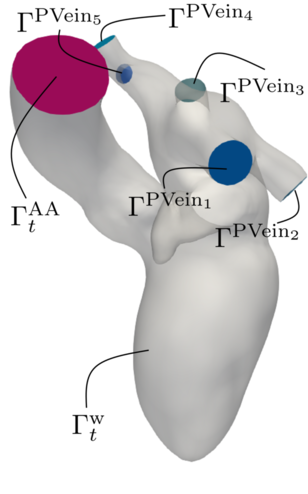}
		\caption{ }
		\label{LH_boundarytags}
	\end{subfigure}
	\begin{subfigure}{.32\textwidth}
		\includegraphics[trim={2 2 2 2 },clip,width=0.95\textwidth]{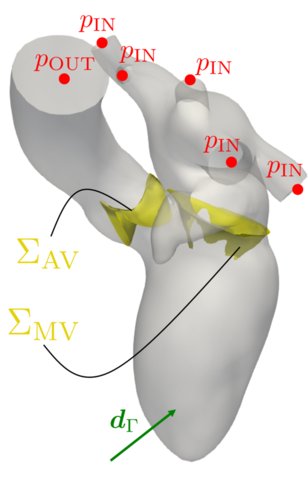}
		\caption{ }
		\label{LH_valves_bcs}
	\end{subfigure}
	\caption{The LH geometry: (A) the three subdomains $\Omega_t = \overline \Omega_t^{\text{LA}} \cup \overline \Omega_t^{\text{LV}} \cup \overline \Omega_t^{\text{AA}}$; (B) the boundary portions of the LH geometry $ \Gamma_t = \left (\bigcup_{i=1}^5 \overline \Gamma^{\text{PVein}_i} \right )\cup \overline \Gamma_t^\text{AA} \cup \overline \Gamma_t^{\text{w}}$; (C) in yellow, the immersed surfaces $\Sigma_{\text{MV}}$ and $\Sigma_{\text{AV}}$ (respectively in their open and closed configurations); in red, the Neumann data (for both inlet and outlet sections); in green the Dirichlet datum at wall.}
	\label{LH_volumetags_boundarytags_bcs}
\end{figure}

We consider a realistic LH geometry provided by Zygote \cite{zygote}, an accurate 3D model of the heart obtained with CT scan data. The LH is constituted by the LA, LV and a portion of the ascending aorta (AA). The two cardiac chambers (LA, LV) are separated by the MV, whereas the AV separates the LV from the AA. The oxygenated blood is collected from the pulmonary veins, the inlets of our domain. The four pulmonary veins are connected to the upper part of the LA. In the geometry considered there are four pulmonary veins, but one of them is splitted into two inlets, thus, our LH geometry is characterized by five inlet sections, as displayed in Figure \ref{LH_volumetags_boundarytags_bcs}. The blood is then pushed into the systemic circulation through the outlet section of the AA. As reported in Figure \ref{LH_volumetags}, we decompose the geometry into three subdomains: $\Omega_t = \overline\Omega_t^{\text{LA}} \cup \overline\Omega_t^{\text{LV}} \cup \overline\Omega_t^{\text{AA}}$, namely the LA $\Omega_t^{\text{LA}}$, the LV $\Omega_t^{\text{LV}}$ and the AA $\Omega_t^{\text{AA}}$. The boundary is split as $\Gamma_t = \left (\bigcup_{i=1}^5 \overline \Gamma^{\text{PVein}_i} \right) \cup \overline \Gamma_t^\text{AA} \cup \overline \Gamma_t^{\text{w}}$, where we denote with $\Gamma^{\text{PVein}_i}$, $i=1, \dots, 5$ the five inlet sections of the four pulmonary veins, with $\Gamma_t^\text{AA}$ the outlet section of the ascending aorta, and with $\Gamma_t^{\text{w}}$ the endocardium (see Figure \ref{LH_boundarytags}). The boundary portions $\Gamma^{\text{PVein}_i}$, $i = 1, \dots, 5$  are fixed. As displayed in Figure \ref{LH_valves_bcs}, we immerse two surfaces $\Sigma_{\text{MV}}$ and $\Sigma_{\text{AV}}$, namely the MV and the AV, in the LH domain. More details on valves dynamics are given in Section \ref{SEC:valves_modeling}.


We prescribe the pressure in the pulmonary veins $p_{\text{IN}}(t)$ on the five inlet sections of our computational domain, yielding the following Neumann BC:
\begin{equation}
	\bm \sigma (\bm u, p) \bm n = - p_{\text{IN}}\bm n, \quad \quad \text{ on } \Gamma^{\text{PVein}_i}\times (0, T_f), \quad i=1, \dots, 5.
	\label{neumann_bc_in}
\end{equation}
On the outlet section $\Gamma_t^\text{AA}$, we prescribe a Neumann BC by setting the outlet pressure $p_{\text{OUT}}(t)$:
\begin{equation}
	\bm \sigma  ( \bm u, p  ) \bm n = - p_{\text{OUT}}\bm n,  \quad \quad \text{ on } \Gamma_t^\text{AA}\times (0, T_f).
	\label{neumann_bc_out}
\end{equation}
Thus, the Neumann boundaries are $\Gamma_t^\mathrm{N}=\bigcup_{i=1}^5\Gamma^{\mathrm{PVein}_i} \cup \Gamma_t^\mathrm{AA}$.
More details on the way we compute inlet and outlet pressures $p_\text{IN}(t)$, $p_\text{OUT}(t)$ are given in Section \ref{SEC:coupling_0D_3D}.
Eventually,  we prescribe the boundary ALE velocity on the endocardium (wall) by time differentiating the wall displacement $\bm d_{\Gamma}$ (that will be introduced in Section \ref{SEC:displacement_modeling}) in the current configuration:
\begin{equation}
	\bm u = \bm u^\mathrm{ALE} = \frac{\partial \widehat{\bm d}_{\Gamma} }{\partial t} \circ \mathcal A_t^{-1}, \quad \quad \text{ on } \Gamma^{\text{w}}_t\times (0, T_f).
	\label{u_gamma_t}
\end{equation}
Hence, the Dirichlet boundary is $\Gamma_t^\mathrm{D} = \Gamma_t^\mathrm{w}$. A graphical sketch of the whole set of BCs is given in Figure \ref{LH_valves_bcs}.

\subsection{Displacement modeling}
\label{SEC:displacement_modeling}

\begin{figure}[!t]
	\centering
	\includegraphics[trim={4 4 4 4},
	clip,
	width=0.5\textwidth]
	{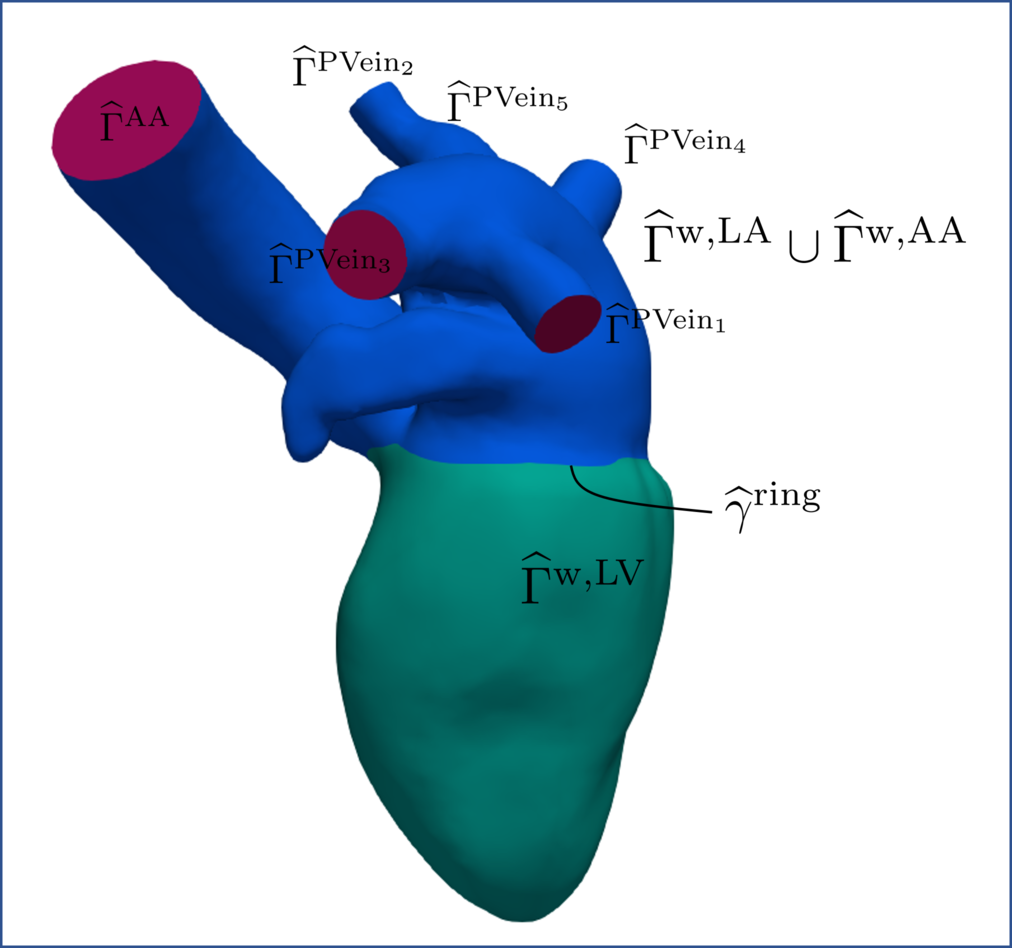}
	\caption{Boundary portions of the LH geometry in reference configuration. This splitting of the domain is used to define the Laplacian problem~\eqref{laplace_beltrami_d_star}. 
	}
	\label{split_wall_preproc}
\end{figure}

The mathematical problem describing the hemodynamics of the LH is solved by prescribing a given displacement field denoted as $\widehat{\bm d}_{\Gamma}$ on the endocardium. In particular, the latter is used as Dirichlet datum on $\Gamma_t^\text{w}$ for the NS-ALE-RIIS problem in Eq. \eqref{u_gamma_t}, and for the geometric problem in Eq. \eqref{ALE_laplacian_BC} on the whole boundary $\widehat{\Gamma}$. In the following, we present a procedure aimed at computing a physiological LH displacement $\widehat{\bm d}_\Gamma $ on $\widehat{\Gamma}$ starting from an EM simulation of the LV.

\begin{figure}[t]
	\centering
	\includegraphics[trim={1 1 1 1},
	clip,
	width=\textwidth]
	{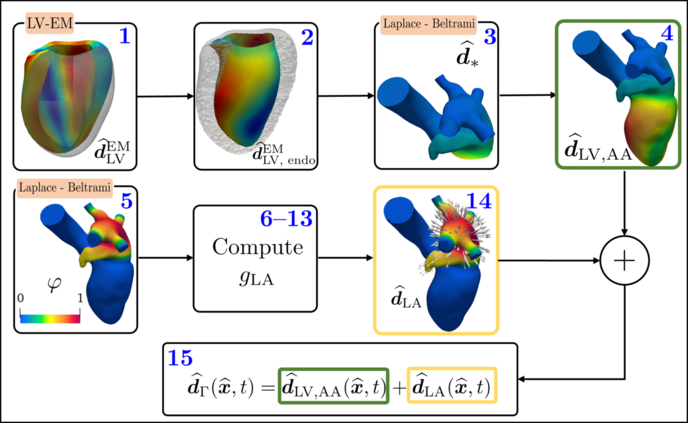}
	\caption{Displacement procedure. Boxes numbers are referred to lines in Algorithm \ref{algo_preprocessing}.}
	\label{displacement_pipeline}
\end{figure}

We split the wall as $\widehat{\Gamma}^\text{w} = \widehat{\Gamma}^\text{w,LV} \cup  \widehat{\Gamma}^\text{w,LA} \cup \widehat{\Gamma}^\text{w,AA} $, being respectively the walls of LV, LA and AA, as shown in Figure \ref{split_wall_preproc}.
Let $\widehat{\bm d}_\text{LV,AA}$ and $\widehat{\bm d}_\text{LA}$ two displacement fields acting on the whole LH which account respectively for the motion of the LV, AA ($\widehat{\bm d}_\text{LV,AA}$) and the LA ($\widehat{\bm d}_\text{LA}$), respectively. We define the displacement on the whole LH as
\begin{equation}
	\widehat{\bm d}_{\Gamma} (\widehat{\bm x}, t) =
	\widehat{\bm d}_\text{LV,AA}(\widehat{\bm x}, t) +
	\widehat{\bm d}_\text{LA}(\widehat{\bm x}, t), \,  \text{ on } \widehat{\Gamma} \times (0, T_f).
	\label{d_gamma_t}
\end{equation}
Thus, we model the LH displacement as the sum of two contributions: a displacement coming from an EM simulation of the LV that we extend to the whole LH domain and an ad-hoc designed displacement for the LA.
In the following, we describe how we compute the aforementioned displacement fields.  In Algorithm~\ref{algo_preprocessing}, we summarize  the main steps required in the preprocessing procedure; in Figure \ref{displacement_pipeline}, we represent these steps with boxes numbered as the lines in Algorithm \ref{algo_preprocessing}.
We point out that, since the displacement at the LV is computed through an EM simulation, and we prescribe this data at walls of the fluid domain, we are enforcing a  one-way (kinematic) coupling condition between EM and CFD in the LV.

\begin{algorithm}[t]
	\caption{Preprocessing procedure to compute LH displacement}
	\label{algo_preprocessing}
	\flushleft
	\textbf{Input}{ $\widehat{\Omega}$, $\widehat \Gamma$}
	\\
	\textbf{Output}{ $\widehat{\bm d}_{\Gamma}$}
	\\
	\begin{algorithmic}[1]
		\STATE{LV-EM simulation $\to \widehat{\bm d}_\text{LV}^\text{EM} \text{ on } \widehat \Omega^\text{LV, s} \times (0, T_f) $} \\
		\STATE{Extract solution on LV endocardium $\to \widehat{\bm d}_\text{LV, endo}^\text{EM} \text{ on } \widehat \Gamma^\text{w, LV} \times (0, T_f)$} \\
		\STATE {$\widehat{\bm d}_\text{LV, endo}^\text{EM} \to $ Laplace - Beltrami (Eq. \eqref{laplace_beltrami_d_star} $\to \widehat{\bm d}_*$ on $\widehat \Gamma \setminus \widehat \Gamma^\text{w, LV} \times (0, T_f)$} \\
		\STATE {Compute $\widehat{\bm d}_\text{LV,AA}$ on $\widehat \Gamma \times (0, T_f)$ (Eq. \eqref{d_LV})}\\
		\STATE {Compute $\widehat \varphi \text{ on } \widehat \Gamma $} \\
		\STATE{ Solve the geometric problem \eqref{geometric_strong} with Dirichlet datum: $\widehat{\bm d}_\text{LV,AA} \text{ on } \widehat \Gamma \times (0, T_f)$ } \\
		\STATE{ Compute $V_\text{LV}(t)$, $V_\text{AA}(t), \,  \forall t \in (0, T_f)$ }
		\\
		\STATE{ Compute $\mathcal{A} (t),\, \forall t \in (0, T_f)$ (Eq. \eqref{flux_A})} \\
		\STATE{Solve the 0D circulation model } $\to V_\text{LA}(t) =  V_\text{LA}^\text{0D}(t) $ \\
		\STATE{Compute $\Phi(t),\,  \forall t \in (0, T_f) $} \\
		\STATE{Compute $\widehat{\bm e}_\text{G}^\text{LA} \text{ on } \widehat \Gamma$} 
		\STATE{Compute $\mathcal B(t),\,  \forall t \in (0, T_f) $} \\
		\STATE{Solve the ODE in Eq. \eqref{ODE_preproc}} $ \to $ $g_\text{LA}(t), \,  \forall t \in (0, T_f) $
		\STATE{$\widehat{\bm d}_\text{LA}=\widehat \varphi \, \widehat{\bm e} _\text{G}^\text{LA}g_\text{LA} \text{ on } \widehat \Gamma \times (0, T_f) $}
		\STATE{$\widehat{\bm d}_{\Gamma}= \widehat{\bm d}_\text{LV,AA}  + \widehat{\bm d}_\text{LA}, \,  \text{ on } \widehat \Gamma \times (0, T_f)$}
	\end{algorithmic}
\end{algorithm}

\subsubsection{EM of a LV and harmonic extension on the whole LH geometry}

We describe how we compute the LV displacement through an EM simulation and then its extension onto the whole LH geometry. The steps here introduced correspond to boxes 1-4 in Figure \ref{displacement_pipeline} and, analogously, to lines 1-4 of Algorithm \ref{algo_preprocessing}.

We use the EM model developed in \cite{REGAZZONI2022111083} consisting of a cardiac EM model of the LV coupled to the 0D circulation model introduced in Section~\ref{SEC:0D_CIRCULATION}.  More specifics on this model, and on the setting we employ to carry out the cardiac EM simulations are given in Appendix~\ref{appendix_em}. In addition, since the one-way coupling approach ignores the dynamic balance among CFD and EM, the pressure in the LV during the isovolumetric phases (i.e. when both valves are closed) is not well defined \cite{fernandez_2019, quarteroni2009numerical}. To overcome this issue, that would however regard only a short portion of the hearbeat, we neglect the isovolumetric phases obtained in the EM simulation for our CFD simulation.

\begin{figure}[t]
	\centering
	\begin{subfigure}{0.32\textwidth}
		\includegraphics[width=\textwidth]{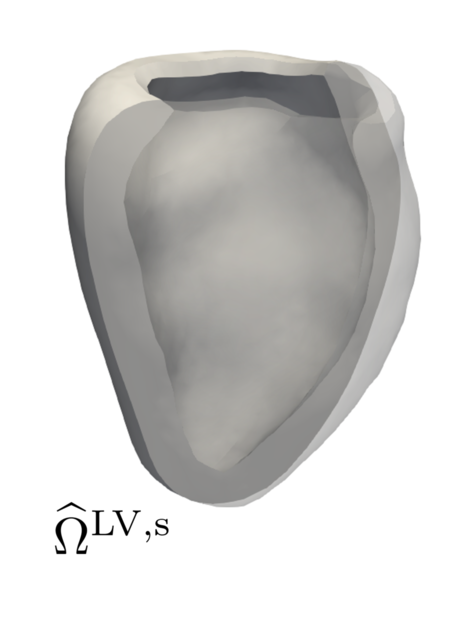}
		\caption{Reference configuration}
		\label{EM_refconf}
	\end{subfigure}
	\begin{subfigure}{0.32\textwidth}
		\includegraphics[width=\textwidth]{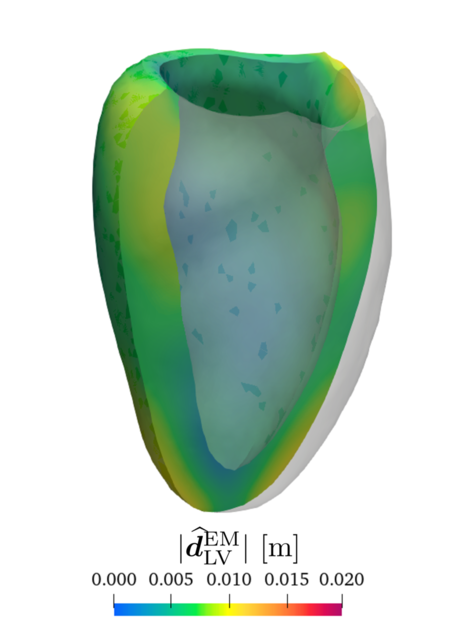}
		\caption{Systole}
		\label{EM_ejection}
	\end{subfigure}
	\begin{subfigure}{0.32\textwidth}
		\includegraphics[width=\textwidth]{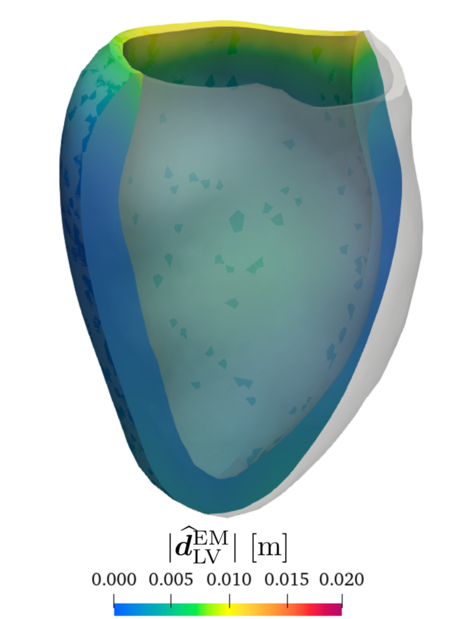}
		\caption{Diastole}
		\label{EM_relaxing}
	\end{subfigure}
	\caption{EM simulation of the LV: (A) LV in its reference configuration; (B), (C) LV during systole and diastole colored by displacement magnitude.}
	\label{EM_LV}
\end{figure}

Let $\widehat{\bm{d}}_\text{LV}^\text{EM}(\widehat {\bm x}, t)$  be the displacement of the LV, solution of the EM model and defined in the LV myocardium in $\widehat{\Omega}^\text{LV, s} \times (0, T_f)$.  We show the LV in its reference configuration ($\widehat{\Omega}^\text{LV,s}$) in Figure~\ref{EM_refconf};  we display snapshots of the numerical solution during systole and diastole in Figures \ref{EM_ejection} and \ref{EM_relaxing}, respectively. Moreover, we also report the EM solution in Figure \ref{displacement_pipeline}, box 1. Let then $\widehat{\bm{d}}_\text{LV,endo}^\text{EM}(\widehat{\bm x}, t)$ on $\widehat{\Gamma}^\text{w,LV} \times (0, T_f)$ be the displacement field restricted to the LV endocardium -- wall of the LV fluid domain -- as shown in Figure \ref{displacement_pipeline}, box 2.

We compute a displacement $\widehat{\bm{d}}_\text{*}$, acting on LA and AA only, as the solution of a vectorial Laplace-Beltrami problem. Specifically, we extend  the ventricular displacement $\widehat{\bm{d}}_\text{LV,endo}^\text{EM}(\widehat{\bm x}, t)$ on LA and AA, by keeping the pulmonary veins fixed (see Figure \ref{displacement_pipeline}, box 3):
\begin{subequations}
	\begin{empheq}[left=\empheqlbrace]{align}
		&-\Delta_{\widehat \Gamma} \widehat{\bm{d}}_\text{*} = \bm 0 & \text{ on } \widehat{\Gamma} \setminus \left( \widehat{\Gamma}^\text{w, LV} \cup \left(\bigcup_{i=1}^5 \widehat{\Gamma}^{\text{PVein}_i}\right)\right )\times (0, T_f), \\
		&\widehat{\bm{d}}_\text{*} = \widehat{\bm{d}}_\text{LV, endo}^\text{EM} & \text{ on } \widehat{\gamma}^\text{ring} \times (0, T_f), \\
		&\widehat{\bm{d}}_\text{*} = \bm 0& \text{ on } \partial\widehat{\Gamma}^{\text{PVein}_i}\times (0, T_f), \quad i=1,\dots, 5,
	\end{empheq}
	\label{laplace_beltrami_d_star}
\end{subequations}
being $\widehat{\gamma}^\text{ring} = \overline{\widehat\Gamma^\text{w, LV}}   \cap \left (\widehat\Gamma \setminus \widehat\Gamma^\text{w, LV} \right )$. We then define the displacement $\widehat{\bm{d}}_\text{LV,AA}(\widehat{\bm x}, t)$ $\text{ on } \widehat{\Gamma} \times (0, T_f)$ as:
\begin{equation}
	\widehat{\bm{d}}_\text{LV,AA}=
	\begin{cases}
		\widehat{\bm{d}}_\text{LV, endo}^\text{EM} & \text{ on } \widehat{\Gamma}^\text{w, LV}\times (0, T_f), \\
		\widehat{\bm{d}}_\text{*} & \text{ on } \widehat{\Gamma} \setminus  \widehat \Gamma^\text{w, LV} \times (0, T_f).
	\end{cases}
	\label{d_LV}
\end{equation}
By combining $\widehat{\bm{d}}_\text{*}$ and $ \widehat{\bm{d}}_\text{LV, endo}^\text{EM}$  as in Eq. \eqref{d_LV}, the resulting $\widehat {\bm d}_\mathrm{LV,AA}$ is hence defined on the whole LH boundary $\widehat{\Gamma}$, as shown in Figure \ref{displacement_pipeline}, box 4. Thus, $\widehat {\bm d}_\mathrm{LV,AA}$ is obtained by extending the EM-based LV displacement on the LH: the LV is moving according to the EM result and the AA according to the harmonic extension of the ventricular motion. As a consequence, at this stage, we are still not accounting for any contribution coming from the LA motion. As a matter of fact, a simple extension of the LV displacement would bring to a non physiological LA volume behavior, having then a direct impact on the fluid dynamics simulation, for instance in terms of non-physiological flowrates and pressures. For this reason,  in Section \ref{sec:modeling_LA_motion}, we propose a simplified model to account for a physiological LA motion.

\subsubsection{Modeling the LA motion}
\label{sec:modeling_LA_motion}
In absence of an EM model of the whole LH (or of the isolated LA), we propose a simplified model to compute the displacement of the whole LH, by also considering the LA physiological motion.

Let $\widehat{\bm d}_\text{LA}(\xhat, t)$ be the (unknown) displacement of the LA introduced in Eq. \eqref{d_gamma_t}. As in \cite{ZDMQ_2020}, we assume that the LA displacement can be modelled by separation of variables, and that it is directed towards its center of volume $\xhat_\text{G}^\text{LA}$. We introduce $\widehat{\bm e}_{\text{G}}^\text{LA}$ as the unit vector directed towards $\xhat_\text{G}^\text{LA}$   as
$	\widehat{\bm e}_{\text{G}}^\text{LA}(\xhat) = \frac{\xhat - \xhat_\text{G}^\text{LA}}{|\xhat - \xhat_\text{G}^\text{LA}|}, $
where $|\cdot|$ is the Euclidean norm. Thus, we define  $\widehat{\bm d}_\text{LA}$ as
\begin{equation}
	\widehat{\bm d}_\text{LA}(\xhat, t) =  \widehat \varphi (\xhat)	\widehat{\bm e}_{\text{G}}^\text{LA}(\xhat) g_\text{LA}(t) , \, \text{on } \widehat{\Gamma} \times (0, T_f),
	\label{d_LA_sep_var}
\end{equation}
being $g_\text{LA}(t)$ a time-dependent function; $\widehat \varphi (\xhat)$ is a smooth scalar function -- computed via Laplace-Beltrami problem -- that we introduce to limit the action of $\widehat{\bm d}_\text{LA}$ on the LA. Specifically, as displayed in Figure~\ref{displacement_pipeline}, Box 5,  $\widehat \varphi:\, \widehat \Gamma \to [0, \, 1]$ is null on the LV, AA, positive on the LA, smoothly vanishing on the pulmonary veins.

Eq. \eqref{d_LA_sep_var} introduces some simplifications since -- differently from the LV -- we are not considering the solution of a suitable EM model. We make these assumptions in order to model the LA displacement in all those cases in which this might not be available as input data: cases in which an EM model is available for the LV only (as in this case) or imaging data given for the LV solely. As a matter of fact, imaging data routinely acquired in diagnostic exams are suitable for reconstructing only the LV as a 3D domain, while in most cases the LA is only visible through single 2D images unsuitable for a 3D reconstruction. Thus, CFD simulations are usually carried out on the LV only, discarding the effects of the atrial flow. This extension procedure can hence be used to fill the missing data and to perform numerical simulations on the whole left part of the heart.

\begin{figure}[t]
	\centering
	\begin{subfigure}{0.19\textwidth}
		\includegraphics[trim={1 1 1 1 },clip,width=\textwidth]{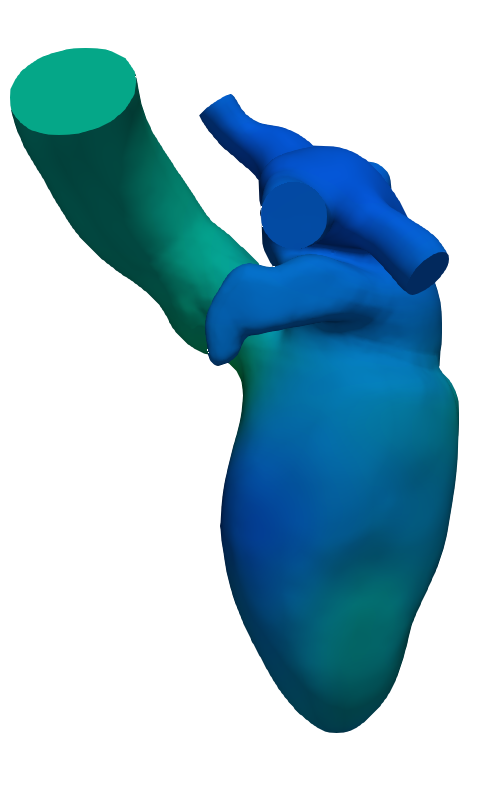}
		\caption{$t = 0.0$ s}
		\label{displacement_0p00}
	\end{subfigure}
	\begin{subfigure}{0.19\textwidth}
		\includegraphics[trim={1 1 1 1 },clip,width=\textwidth]{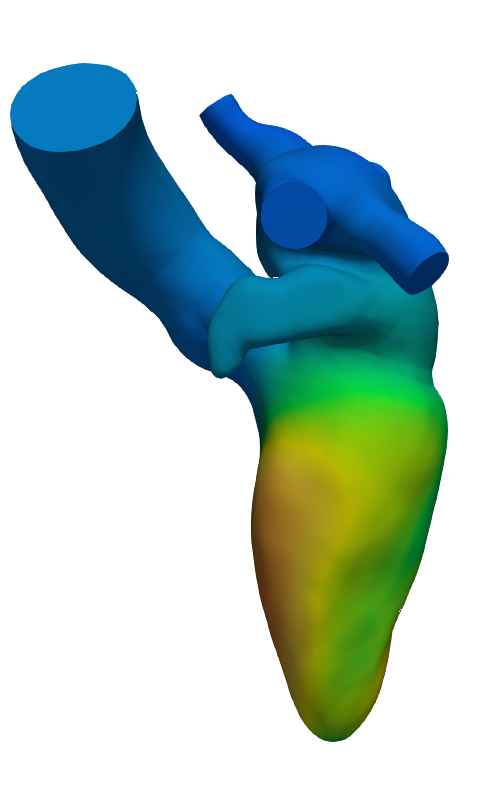}
		\caption{$t = 0.15$ s}
		\label{displacement_0p15}
	\end{subfigure}
	\begin{subfigure}{0.19\textwidth}
		\includegraphics[trim={1 1 1 1 },clip,width=\textwidth]{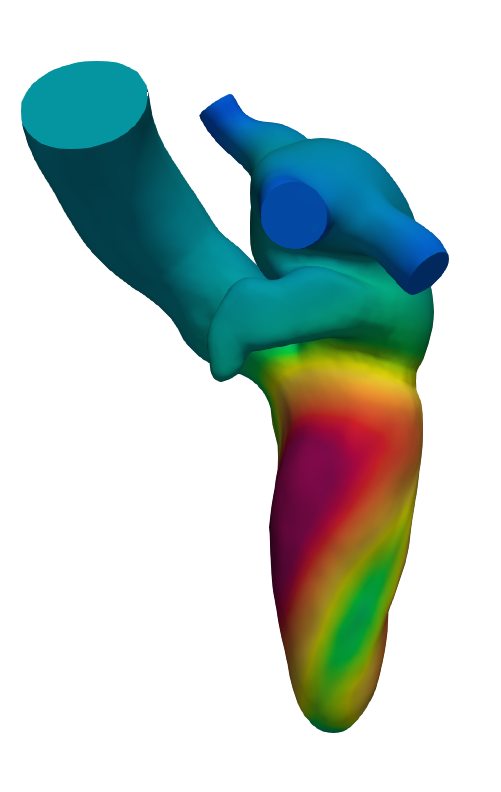}
		\caption{$t = 0.30$ s}
		\label{displacement_0p30}
	\end{subfigure}
	\begin{subfigure}{0.19\textwidth}
		\includegraphics[trim={1 1 1 1 },clip,width=\textwidth]{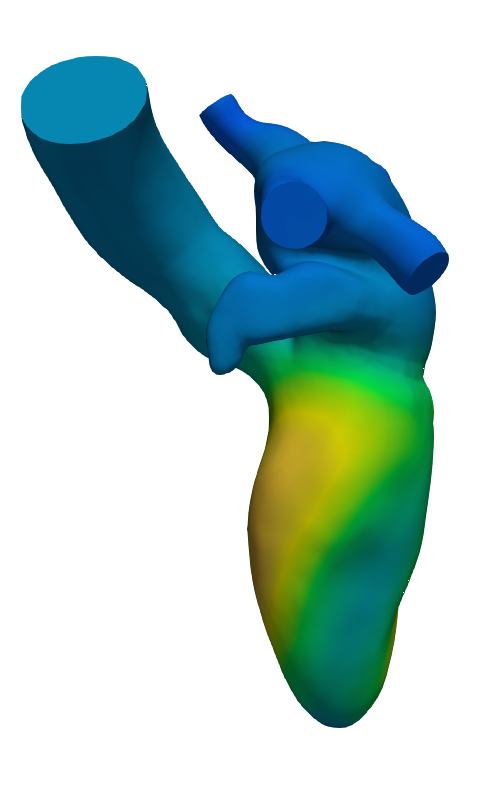}
		\caption{$t = 0.50$ s}
		\label{displacement_0p50}
	\end{subfigure}
	\begin{subfigure}{0.19\textwidth}
		\includegraphics[trim={1 1 1 1 },clip,width=\textwidth]{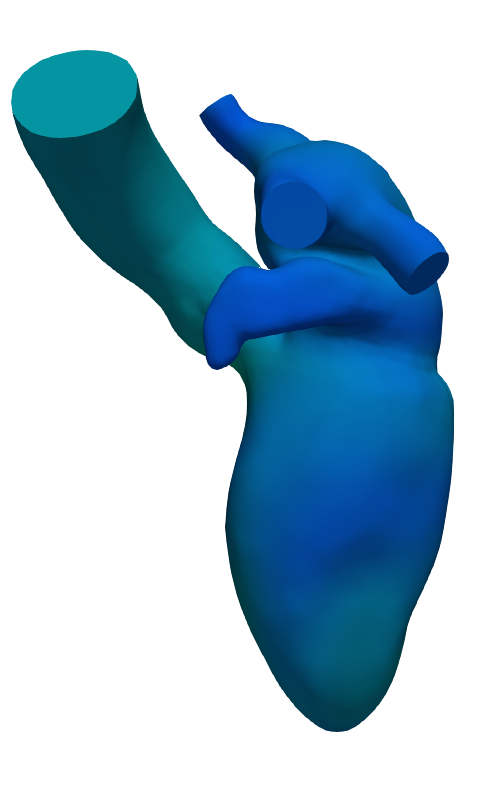}
		\caption{$t = 0.90$ s}
		\label{displacement_0p90}
	\end{subfigure}
	\\
	\begin{subfigure}{0.3\textwidth}
		\includegraphics[trim={0 0 0 0 },clip,width=\textwidth]{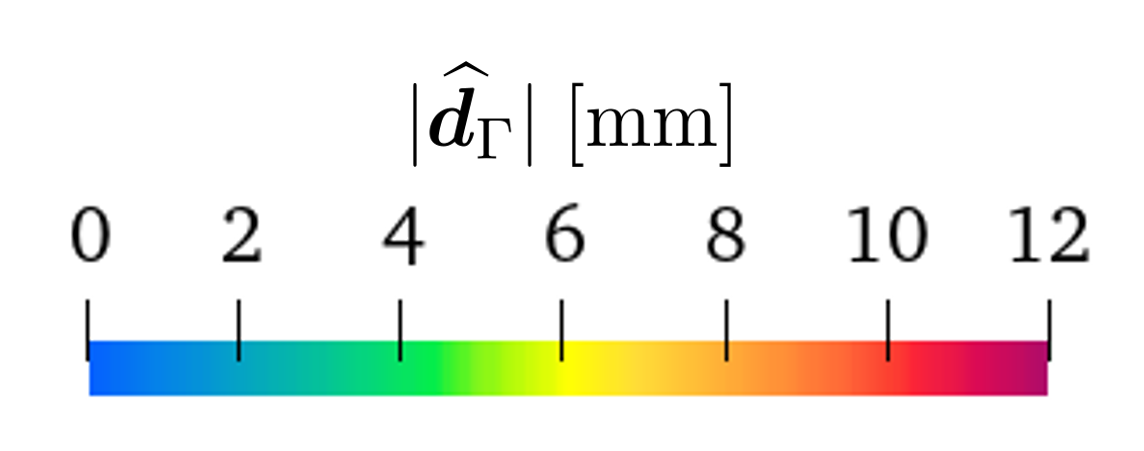}
	\end{subfigure}
	\caption{LH geometry warped by $\widehat{\bm d}_{\Gamma}$ at different times during a heart cycle.}
	\label{displacement_snapshot}
\end{figure}

To compute $g_\text{LA}(t)$, we consider the LH volume time-derivative and we express it through the Reynolds transport theorem (RTT) \cite{kundu} as
\begin{equation}
	\frac{\mathrm d V_\text{LH}}{\mathrm dt} =  \frac{\mathrm d}{\mathrm  d t} \int_{\Omega_t} \mathrm  d\bm x \overset{\text{RTT}}{=} \oint_{\Gamma_t} \bm u^\mathrm{ALE} \cdot \bm n \mathrm d\bm x = \oint_{\Gamma_t} \frac{\partial }{\partial t} {\bm d}_{\Gamma} \cdot \bm n \mathrm d\bm x,
	\label{reynolds_transport_theorem}
\end{equation}
being $V_\text{LH}$ the LH volume. By recalling that for a generic function $w$ holds that $w = \widehat{w} \circ \mathcal{A}_t^{-1}$,
we use Eq. \eqref{d_gamma_t} and \eqref{d_LA_sep_var}
mapped to the current configuration $\Gamma_t$ in Eq. \eqref{reynolds_transport_theorem} to get:
\begin{equation}
	\begin{aligned}
		\frac{\mathrm d V_\text{LH}}{\mathrm dt}  =
		& \oint_{\Gamma_t} \frac{\partial }{\partial t}   {\bm d}_\text{LV,AA}  \cdot \bm n  \mathrm d \bm x + \frac{\mathrm d g_\text{LA}}{\mathrm d t} \oint_{\Gamma_t} \varphi \, {\bm e}_{\text{G}}^\text{LA}  \cdot \bm n \mathrm d \bm x.
	\end{aligned}
	\label{RTT_2}
\end{equation}
Let $V_\text{LA}$, $V_\text{LV}$ and $V_\text{AA}$ be the volumes of LA, LV and AA, with $V_\text{LH}(t) = V_\text{LA}(t)  + V_\text{LV}(t)  + V_\text{AA}(t) $; we define the fluxes in \eqref{RTT_2} as
\begin{align}
	& \mathcal{A} = \oint_{\Gamma_t} \frac{\partial }{\partial t} {\bm d}_\text{LV,AA}   \cdot \bm n \mathrm  d \bm x, \quad
	\mathcal{B} =  \oint_{\Gamma_t} \varphi\, {\bm e}_{\text{G}}^\text{LA} \cdot \bm n \mathrm  d \bm x, \label{flux_A}
	\\
	& \Phi(t)  = \dfrac{\mathrm d V_\text{LH}(t)}{\mathrm d t} - \mathcal{A}(t)  = \dfrac{\mathrm d V_\text{LA}(t)}{\mathrm  dt} + \dfrac{\mathrm d V_\text{LV}(t)}{\mathrm dt} + \dfrac{\mathrm d V_\text{AA}(t)}{\mathrm dt} - \mathcal{A}(t). \label{PHI}
\end{align}
Solving Eq. \eqref{RTT_2} for $\frac{\mathrm  d g_\text{LA}(t)}{\mathrm  d t}$, yields the following Cauchy problem:
\begin{equation}
	\begin{cases}
		\dfrac{\mathrm  d g_\text{LA}(t)}{\mathrm d t} & =  \dfrac{\Phi(t)}{\mathcal{B}(t)}, \quad t \in (0, T_f),
		\\
		g_\text{LA}(0) & = g_{\text{LA}_0}.
	\end{cases}
	\label{ODE_preproc}
\end{equation}
Since the Zygote's geometry is generated at a time instant corresponding to the 70\% of diastole in a heartbeat \cite{zygote}, we consider such instant as initial time $t=0$ whence $\widehat{\Omega}={\Omega}_0$, and $g_{\text{LA}_0}=0$.
To compute $V_\text{LA}(t)$ in Eq. \eqref{PHI}, we solve the 0D closed-loop circulation model -- tuned on the basis of the EM simulation\footnote{In principle, one could directly adopt the LA volume obtained with the 3D-0D EM model \cite{REGAZZONI2022111083}. However, as previously explained, since we neglect the isovolumetric phases, we calibrate the 0D model to be consistent with the LV volume achieved.}-- and we denote the LA volume computed as $V_\mathrm{LA}^\mathrm{0D}(t)$. Thus, we set:
\begin{equation*}
	V_\mathrm{LA}(t) = V_\mathrm{LA}^\mathrm{0D}(t).
\end{equation*}

\begin{figure}[t]
	\centering
	\includegraphics[trim={1 1 1 1},clip,width=\textwidth]{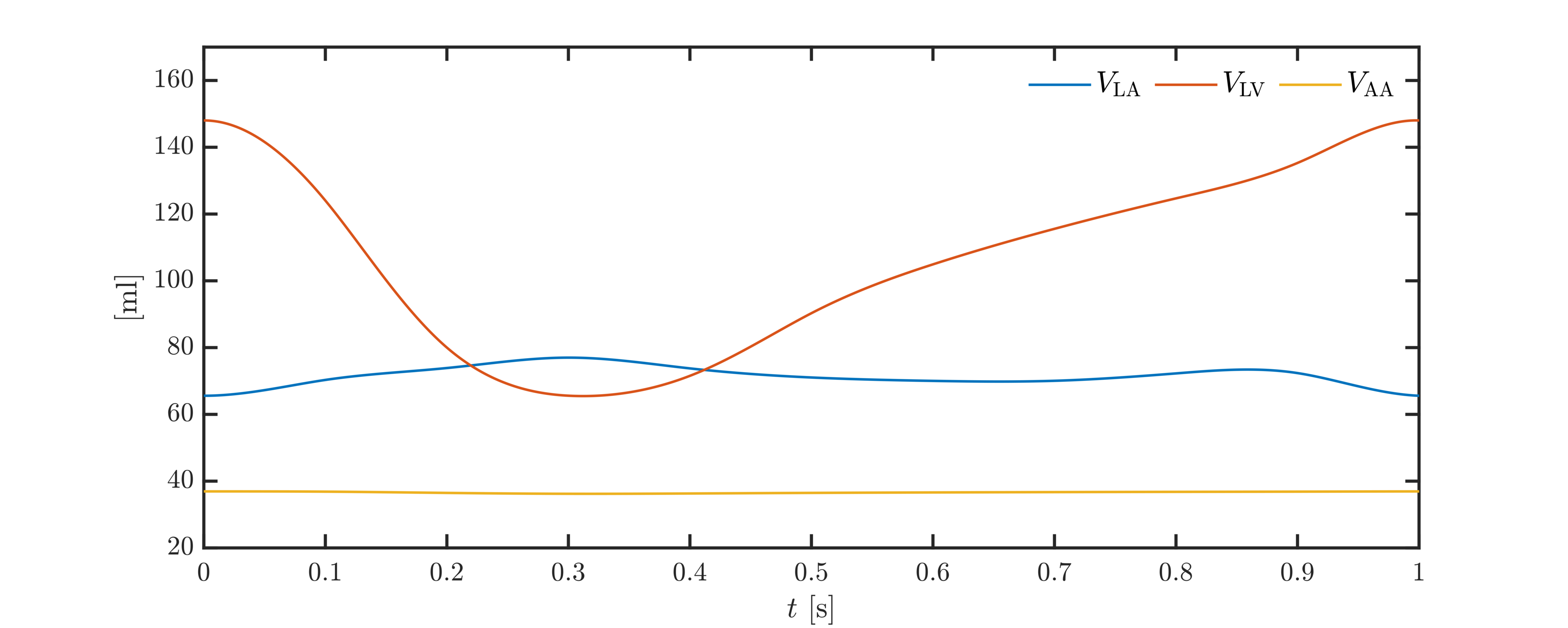}
	\caption{Volumes of LA, LV and AA achieved applying the displacement -- $\widehat{\bm d}_{\Gamma}$ defined in Eq. \eqref{d_gamma_t} -- to $\widehat \Gamma$ .}
	\label{volumes}
\end{figure}

As reported in Algorithm \ref{algo_preprocessing}, in order to compute $g_\text{LA}(t)$, we first solve a preliminary geometric problem, as the one in Eq. \eqref{geometric_strong} but using $\widehat{\bm d}_\text{LV,AA}$ as Dirichlet datum on $\widehat{\Gamma} \times (0, T_f)$. In this way, we compute once and for all the volumes $V_\text{LV}(t)$, $V_\text{AA}(t)$ and the flux $\mathcal{A}(t)$. We solve the 0D circulation model to get $V_\text{LA}(t)$, and we compute the flux $\Phi(t)$. Once we compute $\widehat{\bm e}_\mathrm{G}^\mathrm{LA}$, we calculate the flux $\mathcal B(t)$ to finally get $g_\text{LA}(t)$ as solution of the Cauchy problem \eqref{ODE_preproc}. The steps here mentioned are aimed at computing $g_\text{LA}(t)$: they are shown in lines 6-13 of Algorithm \ref{algo_preprocessing} and gathered in a single box in Figure \ref{displacement_pipeline}.
Once we get ${g}_\text{LA}(t)$, we can compute $\widehat{\bm d}_\mathrm{LA}$ as in Eq. \eqref{d_LA_sep_var} and displayed in box 14 in Figure \ref{displacement_pipeline}. The magnitude of the displacement $\widehat{\bm d}_\mathrm{LA}$ is shown in the yellow box in Figure \ref{displacement_pipeline}: it is non-null on the LA only, smoothly vanishing towards the pulmonary veins and on the LV and the AA. We also report glyphs of $\widehat{\bm d}_\mathrm{LA}$ showing that the LA displacement direction coincides with $\widehat{\bm e}_\mathrm{G}^\mathrm{LA}$.
The displacement $\bm d_{\Gamma}$ can be eventually computed as in Eq. \eqref{d_gamma_t} and as reported in Figure \ref{displacement_pipeline}, box 15.

We implemented the whole procedure in $\mathsf{PvPython}$ \cite{pvpython}. EM simulations of the LV and the geometric problem  \eqref{geometric_strong} with Dirichlet datum $\widehat{\bm d}_\text{LV,AA}$ are carried out using the in-house finite element library $\mathsf{life^x}$ \cite{lifex} (more details will be provided in Section \ref{SEC:numerical_results}). We solved the Laplace-Beltrami problems to get $\widehat{\bm d}_*$, and $\widehat \varphi$ using the harmonic extension algorithms proposed in \cite{FQ_2021} and implemented in a public form of
$\mathsf{vmtk}$ \cite{vmtk, vmtk_url, vmtk_fedele}.

In Figure \ref{displacement_snapshot}, we finally report the LH geometry warped by $\widehat{\bm d}_{\Gamma}$ at different time steps during the heart cycle starting from the time of end diastole. Moreover, in Figure \ref{volumes}, we show the volumes of LA, LV and AA obtained with our preprocessing procedure. We measure ventricular end diastolic volume EDV = 148.04 ml, end systolic volume ESV = 65.48 ml, stroke volume SV = EDV - ESV = 82.56 ml and the ejection fraction EF = 100 SV/EDV = 55.77$\%$. The results obtained are consistent with values routinely acquired in healthy subjects \cite{maceira2006normalized, viola2020fluid}.

\subsection{Valves dynamics}
\label{SEC:valves_modeling}

\begin{figure}[t]
	\centering
	\includegraphics[trim={4 4 4 4 },clip,width=0.8\textwidth]{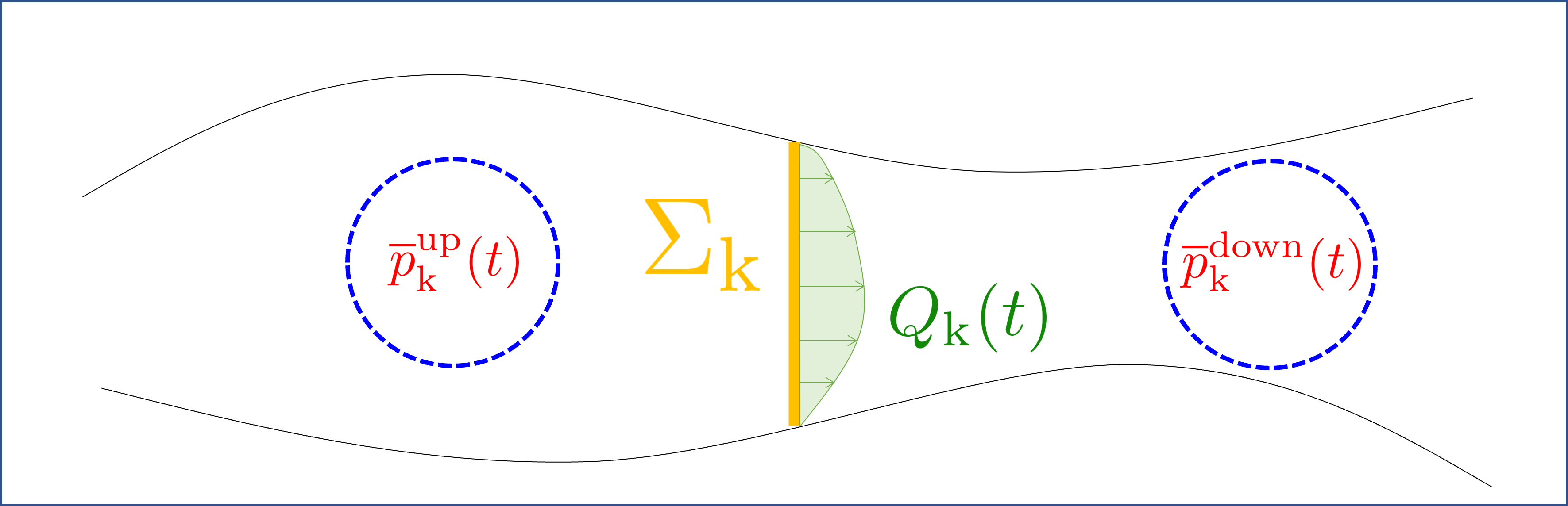}
	\caption{Immersed valve $\Sigma_\mathrm k$ with upwind and downwind control volumes where average pressures are computed. $Q_\mathrm k$ is the flowrate across $\Sigma_\mathrm k$. This picture corresponds to a simple a two-dimensional fluid domain for the sake of simplicity.}
	\label{control_volumes_valve}
\end{figure}

\begin{figure}[!t]
	\centering
	\begin{subfigure}{0.49\textwidth}
		\centering
		\includegraphics[trim={1 1 1 1 },clip,width=0.9\textwidth]{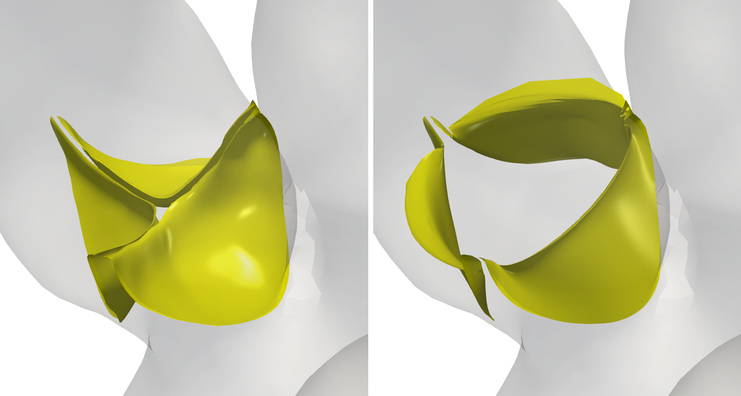}
		\caption{$\Sigma_\mathrm{AV}$}
		\label{av}
	\end{subfigure}
	\begin{subfigure}{0.49\textwidth}
		\centering
		\includegraphics[trim={1 1 1 1 },clip,width=0.9\textwidth]{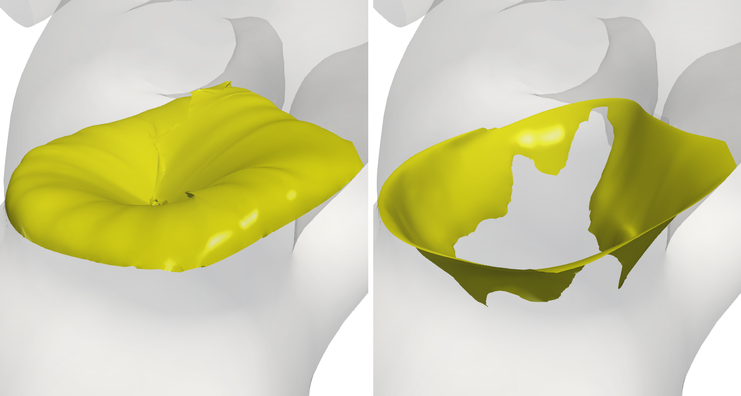}
		\caption{$\Sigma_\mathrm{MV}$}
		\label{mv}
	\end{subfigure}
	\caption{Cardiac valves in their fully closed and fully open configurations.}
	\label{valves_frame}
\end{figure}

As discussed in Section \ref{SEC:NS_ALE_RIIS}, the effect of MV and AV on the fluid is modeled by means of the RIIS method \cite{Fedele_RIIS, Fumagalli2020}. The MV and AV geometries are given by Zygote \cite{zygote} in their closed and open configuration, respectively. To bring the immersed surfaces from their open to closed configuration (and viceversa) we prescribe a displacement field computed via Laplace-Beltrami problems, and algorithms exploiting the closest-point distances \cite{FQ_2021}. We show the cardiac valves in their fully closed and fully open configurations in Figure \ref{valves_frame}.

We do not prescribe a priori the time at which valves open and close, as this is a result of our numerical simulations. If a valve is closed, it opens when the pressure jump across it becomes positive; viceversa, the valve closes when a condition of reversed flow across the orifice area is detected \cite{fernandez_2019}.
The condition on the pressure jump across the valve $\mathrm k$, with $\mathrm{k} = \mathrm{MV}, \, \mathrm{AV}$ is checked by considering two control volumes inside the upstream and downstream chambers, as shown in Figure \ref{control_volumes_valve}. Let $\overline{p}^\text{up}_{\mathrm k}(t)$ and  $\overline{p}^\text{down}_{\mathrm k} (t)$ be the average pressure inside each control volume, the valve $\mathrm k$ opens when the pressure jump is positive, i.e. $\delta \overline{p}_{\mathrm k}(t) = \overline{p}_{\mathrm k}^\text{up}(t)-\overline{p}_{\mathrm k}^\text{down}(t)
>0$.

The closure condition at time $t$ of reversed flow is satisfied when $Q_\mathrm k (t) < 0$, $Q_\mathrm k$ being the flowrate through the valve $\mathrm k$. Hinging upon a mass balance, the sign of $Q_\mathrm k(t)$ is directly related to that of $\dot{V}_\mathrm{LV}(t)$. Therefore, the MV will close when  $\dot{V}_\mathrm{LV}(t)<0$, while the AV when $\dot{V}_\mathrm{LV}(t)>0$.

\section{Coupling the 3D fluid dynamics model with the 0D circulation model}
\label{SEC:coupling_0D_3D}
\begin{figure}[!t]
	\centering
	\includegraphics[trim={1 1 1 1},
	clip,
	width=0.7\textwidth]
	{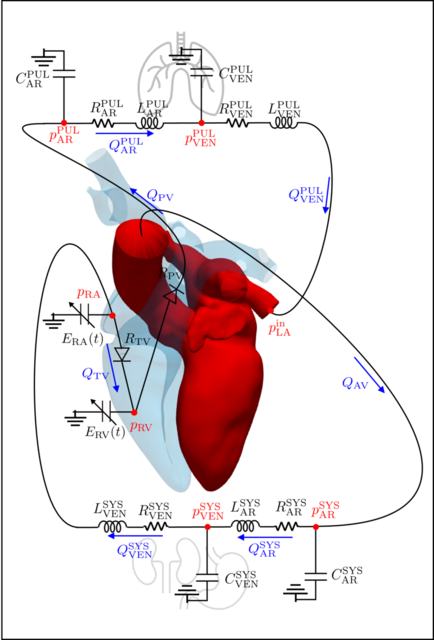}
	\caption{The geometric multiscale model: coupling between the 3D CFD model of the LH and the 0D circulation model of the remaining cardiocirculatory system.}
	\label{circuit_0D_3D}
\end{figure}

In order to couple the 3D CFD model of the LH with the 0D circulation model of the whole cardiovascular system, we first remove from the 0D model Eq. \eqref{0D_VLA}, \eqref{0D_VLV}, \eqref{0D_pressurechambers} for $\mathrm i = \mathrm{LA}, \mathrm{LV}$ and \eqref{0D_flowrates} for $\mathrm j = \mathrm{MV}, \mathrm{AV}$, and we replace them with the 3D model of the LH. The coupling between the 3D and 0D models consists of the enforcement of the  continuity of flowrates and pressures  on the ``artificially chopped'' boundaries, i.e. inlet and outlet sections of the 3D domain: $\Gamma^{\text{PVein}_i}$, $i=1, \dots, 5$ and $\Gamma_t^\mathrm{AA}$.

Referring to Figure \ref{circuit_0D_3D}, the LA pressure $p_\text{LA}(t)$ of the 0D model, appearing in Eq. \eqref{0D_QVENPUL}, represents the pressure downwind the ``PUL VEN'' RLC system, i.e. the one to be prescribed on the inlet sections of our 3D domain. Thus, we rename the latter as $p_\text{LA}^\text{in}(t)$. Analogously, the systemic arterial  pressure $p_\text{AR}^\text{SYS}(t)$, appearing in Eq. \eqref{0D_pARSYS} and \eqref{OD_QARSYS}, acts as the outlet pressure on the AA outlet section $\Gamma_t^\text{AA}$. Moreover, the flowrate $Q_\text{VEN}^\text{PUL}(t)$, appearing in Eq. \eqref{0D_pVENPUL} and \eqref{0D_QVENPUL}, represents the inlet flowrate in the 3D model. Similarly, the flowrate $Q_\text{AV}(t)$ in Eq.~\eqref{0D_pARSYS} represents the outlet flowrate in our 3D model. The interface conditions between the 3D and 0D models are then expressed as:
\begin{equation}
	p_\mathrm{IN}(t) = p_\mathrm{LA}^\mathrm{in}(t), \; \;
	p_\mathrm{OUT}(t) = p_\mathrm{AR}^\mathrm{SYS}(t),
	\label{continuity_p}
\end{equation}
for the continuity of pressures, and
\begin{equation}
	Q_\mathrm{IN}(t)=-Q^\mathrm{PUL}_\mathrm{VEN}(t), \; \;
	Q_\mathrm{OUT}(t)=Q_\mathrm{AV}(t),
	\label{continuity_Q}
\end{equation}
for the continuity of flowrates\footnote{We define the sign of the flowrate in accordance with the outward unit normal $\bm n$. Thus, an inlet flowrate (entering velocity) will be, by definition, negative.}. $Q_\mathrm{IN}(t)$ and $Q_\mathrm{OUT}(t)$ are the inlet and outlet flowrates respectively, defined as
\begin{equation}
	Q_\mathrm{IN}(t)  = \sum_{i=1}^5 \int_{\Gamma^{\mathrm{PVein}_i}} u_n(\bm x, t)\, \mathrm d \bm x, \; \;
	Q_\mathrm{OUT}(t)  = \int_{\Gamma^\mathrm{AA}_t} u_n(\bm x, t)\, \mathrm d \bm x,
	\label{definition_Q_in_Q_out}
\end{equation}
with $u_n = \left (\bm u - \bm u^\mathrm{ALE} \right ) \cdot \bm n $ the relative normal velocity.

The resulting 0D model is the following:
for any $t \in (0, T_f)$:
\begin{subequations}
	\allowdisplaybreaks
	\begin{align}
		\frac{\mathrm dV_\text{RA}(t)}{\mathrm dt} &  =  Q_\text{VEN}^\text{SYS}(t) - Q_\text{TV}(t),
		\label{0D_3D_VRA}  \\
		\frac{\mathrm dV_\text{RV}(t)}{\mathrm dt} &  =  Q_\text{TV}(t) - Q_\text{PV}(t),
		\label{0D_3D_VRV}  \\
		\frac{\mathrm dp_\text{AR}^\text{SYS}(t)}{\mathrm dt} & = \frac{1}{C_\text{AR}^\text{SYS}} \left (Q_\text{AV}(t) - Q_\text{AR}^\text{SYS}(t)
		\right ), \label{0D_3D_pARSYS} \\
		\frac{\mathrm dp_\text{VEN}^\text{SYS}(t)}{\mathrm dt} & = \frac{1}{C_\text{VEN}^\text{SYS}} \left ( Q_\text{AR}^\text{SYS}(t) - Q_\text{VEN}^\text{SYS}(t) \right ),
		\label{0D_3D_pVENSYS} \\
		\frac{\mathrm dp_\text{AR}^\text{PUL}(t)}{\mathrm dt} & = \frac{1}{C_\text{AR}^\text{PUL}} \left (Q_\text{PV}(t) - Q_\text{AR}^\text{PUL}(t) \right ),
		\label{0D_3D_pARPUL} \\
		\frac{\mathrm dp_\text{VEN}^\text{PUL}(t)}{\mathrm dt} & = \frac{1}{C_\text{VEN}^\text{PUL}} \left( Q_\text{AR}^\text{PUL}(t) - Q_\text{PUL}^\text{VEN}(t) \right ),
		\label{0D_3D_pVENPUL} \\
		\frac{\mathrm dQ_\text{AR}^\text{SYS}(t)}{\mathrm dt} & = \frac{R_\text{AR}^\text{SYS}}{L_\text{AR}^\text{SYS}} \left (- Q_\text{AR}^\text{SYS}(t) - \frac{p_\text{VEN}^\text{SYS}(t)-p_\text{AR}^\text{SYS}(t)}{R_\text{AR}^\text{SYS}} \right ),
		\label{0D_3D_QARSYS} \\
		\frac{\mathrm dQ_\text{VEN}^\text{SYS}(t)}{\mathrm dt} & = \frac{R_\text{VEN}^\text{SYS}}{L_\text{VEN}^\text{SYS}} \left (- Q_\text{VEN}^\text{SYS}(t) - \frac{p_\text{RA}(t)-p_\text{VEN}^\text{SYS}(t)}{R_\text{VEN}^\text{SYS}} \right ),
		\label{0D_3D_QVENSYS} \\
		\frac{\mathrm dQ_\text{AR}^\text{PUL}(t)}{\mathrm dt} & = 	\frac{R_\text{AR}^\text{PUL}}{L_\text{AR}^\text{PUL}} \left (- Q_\text{AR}^\text{PUL}(t) - \frac{p_\text{VEN}^\text{PUL}(t)-p_\text{AR}^\text{PUL}(t)}{R_\text{AR}^\text{PUL}} \right ),
		\label{0D_3D_QARPUL}
	\end{align}
	\label{0D_3D_ODE}
\end{subequations}
being
\begin{subequations}
	\allowdisplaybreaks
	\begin{align}
		p_\text{LA}^\text{in}(t) & = p_\text{VEN}^\text{PUL}(t) - R_\text{VEN}^\text{PUL}Q_\text{VEN}^\text{PUL}(t) - L_\text{VEN}^\text{PUL}\dfrac{\mathrm d Q_\text{VEN}^\text{PUL, 3D}(t)}{\mathrm dt},
		\label{0D_3D_pVENPUL_down} \\
		p_\text{RA}(t) & = p_\text{EX}(t) + E_\text{RA}(t) \left ( V_\text{RA}(t) - V_{0, \text{RA}}\right ),
		\label{0D_3D_pRA} \\
		p_\text{RV}(t) & = p_\text{EX}(t) + E_\text{RV}(t) \left ( V_\text{RV}(t) - V_{0, \text{RV}}\right ),
		\label{0D_3D_pRV} \\
		Q_\text{TV}(t) & = \frac{p_\text{RA}(t)-p_\text{RV}(t)}{R_\text{TV}(p_\text{RA}(t), p_\text{RV}(t))},
		\label{0D_3D_QTV} \\
		Q_\text{PV}(t) & = \frac{p_\text{RV}(t)-p_\text{AR}^\text{PUL}(t)}{R_\text{PV}(p_\text{RV}(t), p_\text{AR}^\text{PUL}(t))}.
		\label{0D_3D_QPV}
	\end{align}
	\label{0D_3D_p_Q}
\end{subequations}
Note that $p_\text{LA}(t)$ has been replaced by $p_\text{LA}^\text{in}(t)$ and, as reported in Eq.~\eqref{0D_3D_pVENPUL_down}, it is obtained by solving Eq.~\eqref{0D_QVENPUL} for  $p_\text{LA}^\text{in}(t)$. We gather the unknown variables of Eq. \eqref{0D_3D_ODE} in a vector  $\bm c(t) =  (
V_\text{RA}(t), \,
V_\text{RV}(t), \,
p_\text{AR}^\text{SYS}(t), \,
p_\text{VEN}^\text{SYS}(t), \,
p_\text{AR}^\text{PUL}(t), \,
p_\text{VEN}^\text{PUL}(t), \, $ $
Q_\text{AR}^\text{SYS}(t), \,
Q_\text{VEN}^\text{SYS}(t), \,
Q_\text{AR}^\text{PUL}(t) )^T$ and the left-hand side of Eq. \eqref{0D_3D_p_Q} in a vector $\tilde {\bm c}(t) = (p_\text{LA}^\text{in}(t), \, $ $
p_\text{RA}(t), \,
p_\text{RV}(t), \,
Q_\text{TV}(t), \,
Q_\text{PV}(t) )^T$. We then collect the right-hand sides of Eq. \eqref{0D_3D_ODE} and \eqref{0D_3D_p_Q}  in the vectors $\bm r(t, \bm c(t), \tilde {\bm c}(t) )$ and  $\tilde{\bm r}(t, \bm c(t) )$, respectively.
Hence, the reduced 0D model, enriched with suitable initial conditions, is expressed in a compact form as
\begin{subequations}
	\begin{empheq}[left=\empheqlbrace]{align}
		&\dfrac{\mathrm d \bm c (t)}{\mathrm d t} = \bm r (t, \bm c (t),  \tilde{\bm c}(t)) & t \in (0, T_f), \\
		&\tilde{\bm c}(t) = \tilde{\bm r}(t, {\bm c}(t)) & t \in (0, T_f),
		\\
		&\bm c(0) = \bm c_0. &
	\end{empheq}
	\label{0D_3D_compact}
\end{subequations}
A graphical representation of the multiscale problem is given in Figure \ref{circuit_0D_3D}. The overall 3D-0D continuous problem is expressed as: find $\bm u, \, p, \, \bm c, \, \tilde{\bm c}$ such that:
\begin{subequations}
	\allowdisplaybreaks
	\begin{empheq}[left=\empheqlbrace]{align}
		& \rho \left ( \frac{\widehat \partial \bm u}{\partial t} +  \left( \left( \bm{u} - \bm{u}^{\text{ALE}} \right) \cdot \nabla \right) \bm{u} \right ) & \notag
		\\
		& - \nabla \cdot \bm{\sigma} (\bm u, p) & \notag
		\\
		&  + \sum_{\mathrm k }\dfrac{R_{\text{k}}}{\varepsilon_{\text{k}}}\delta_{\Sigma_{\text{k}}, \varepsilon_{\text{k}}} (\bm u -\bm u^\mathrm{ALE}) =  \bm 0 &  \text{ in } \Omega_t \times (0, T_f),
		\label{coupled_3D_mom}
		\\
		& \nabla \cdot \bm u  =   0  &   \text{ in } \Omega_t \times (0, T_f),
		\label{coupled_3D_div}
		\\
		& \bm \sigma (\bm u, p) \bm n = - p_\text{LA}^\text{in} \bm n & \text{ on } \Gamma^{\text{PVein}_i}\times (0, T_f), \, i=1, \dots, 5,
		\label{coupled_3D_NeumannBC_IN}
		\\
		& \bm \sigma ( \bm u, p ) \bm{n} =  - p_\text{AR}^\text{SYS}\bm{n} & \text{ on } \Gamma_t^\text{AA}\times (0, T_f),
		\label{coupled_3D_NeumannBC_OUT}
		\\
		& \bm u = \bm u^\mathrm{ALE} & \text{ on } \Gamma^{\text{w}}_t\times (0, T_f),
		\label{coupled_3D_wall}
		\\
		& \bm u = \bm 0 & \text{ in } \Omega_0 \times \{0\},
		\label{coupled_3D_IC}
		\\
		& - \Delta \widehat{\bm{d}} =  \,  \bm 0 & \text{ in } \widehat{\Omega}, \label{coupled_lifting}\\
		& \widehat{\bm d} =  \, \widehat{\bm{d}}_{\Gamma} & \text{ on } \widehat{\Gamma},  \label{coupled_lifting_BC}\\
		& \bm u^\text{ALE} = \left ( \frac{\partial \widehat{\bm d}}{\partial t}\right) \circ \mathcal A_t^{-1}, \label{coupled_uALE_dd_dt}
		\\
		& \dfrac{\mathrm d \bm c (t)}{\mathrm dt} = \bm r (t, \bm c (t),  \tilde{\bm c}(t)) & \text{ for } t \in (0, T_f),
		\label{coupled_ODE}\\
		& \tilde{\bm c}(t) = \tilde{\bm r}(t, {\bm c}(t)) & \text{ for }t \in (0, T_f),
		\label{coupled_ODE_stationary}\\
		& \bm c= \bm c_0 & \text{ for } t=0, \label{coupled_ODE_IC}
	\end{empheq}
	\label{fluid_geo_circ_strong}
\end{subequations}
with $\mathrm k = \text{MV, AV}$.

\section{Numerical methods}
\label{SEC:numerical_methods}
In this section, we present the numerical methods to solve the multiscale problem \eqref{fluid_geo_circ_strong}. Specifically, we present the numerical schemes for the solution of the NS-ALE-RIIS equations, the geometric problem, and the 0D circulation model, respectively.
Finally, we introduce the segregated numerical scheme to solve the the coupled problem.

\subsection{Space and time discretization of the NS-ALE-RIIS equations and VMS-LES method}
\label{SEC:numerical_NS_ALE_RIIS}
For the space discretization of Eq. \eqref{NS_ALE_RIIS_strong}, we introduce the following infinite dimensional function spaces:
\begin{equation}
	\bm{\mathcal V}_{\bm g} := \, \left\{ \bm v \in [H^1(\Omega_t)]^3: \bm v = \bm g \text{ on } \Gamma_t^\mathrm{D}\right\}, \quad
	\mathcal Q := \,  L^2(\Omega_t).
	\label{inf_function_spaces}
\end{equation}
We employ the FE method for the spatial discretization and we denote with the superscript $h$ quantities associated to the FE discretization. We use the VMS-LES method \cite{bazilevs_2007, DF_2015, DMQ_2019} to get an \textit{inf-sup} stable FE approximation of the NS equations. Moreover, this formulation is suitable both to control instabilities arising from advection dominated problems, and to account for the transitional and nearly turbulent regime that typically occurs in cardiac blood flow \cite{bazilevs_2007, ZDMQ_2020}.

For the time discretization, we partition the time domain  in $N_t$ time steps of equal size $\Delta t = \frac{T_f}{N_t}$, and we denote with the subscript $n$ quantities evaluated at time step $n$: $t_n$, with $n = 0, \dots, N_t$. We use the Backward Differentiation Formula (BDF) scheme of order $\sigma_t = 1, 2, 3$ to discretize the problem in time \cite{DF_2015}, and we use a semi-implicit treatment of the non linearities by extrapolating the velocity field  by means of the Newton-Gregory backward polynomials of order $\sigma_t$, yielding the extrapolated velocity $\bm u^h_{n+1, \text{EXT}}$. From our numerical results, and consistently with the findings of \cite{DF_2015}, we found that the semi-implicit VMS-LES formulation guarantees a stable solution also for relatively ``large'' time-step sizes, differently from explicit numerical schemes where the time step restriction is generally more severe \cite{quarteroni2009numerical}. By defining the extrapolated relative velocity as $\bm u^h_* =\bm u^h_{n+1, \text{EXT}} - \bm u^{\text{ALE},h}_{n+1}$, the fully discretized linear semi-implicit VMS-LES formulation of the NS-ALE-RIIS equations with BDF as time integration method reads:
given $\bm u^h_{n}, \dots, $ $\bm u^h_{n+1-\sigma_t}$, for any $n = 0, \dots, N_t -1$, find $(\bm u^h_{n+1}, p^h_{n+1})  \in \mathcal {V}_{\bm g}^h \times \mathcal Q^{h}$  such that:

\begin{equation}
	\begin{aligned}
		\allowdisplaybreaks
		& \left ( \bm v^h, \rho \frac{\alpha_\text{BDF} \bm u^h_{n+1}}{\Delta t}\right )_{\Omega_{n+1}}
		+
		\left ( \bm v^h, \rho \left ( \bm u^h_* \cdot \nabla \right ) \bm u^h_{n+1}\right)_{\Omega_{n+1}}
		\\
		+
		&
		\left (\nabla \bm v^h, \mu \nabla \bm u^h_{n+1} \right)_{\Omega_{n+1}}  		
		-\left ( \nabla \cdot \bm v^h, p^h_{n+1} \right)_{\Omega_{n+1}}
		+
		\left ( q^h, \nabla \cdot \bm u^h_{n+1} \right )_{\Omega_{n+1}}
		\\
		+
		&
		\left( \bm v^h,  \sum_{k=1}^{m}\frac{R_k}{\varepsilon_k}\delta_{\Sigma_{k, n+1}, \varepsilon_k}\left(\bm u^h_{n+1} - \bm u^{\mathrm{ALE},h}_{n+1} \right ) \right )_{\Omega_{n+1}}
		\\
		+
		&
		\, \mathcal{S}(\bm v^h, q^h, \bm u^h_{n+1}, \bm u^{\text{ALE},h}_{n+1},  \bm u^h_{n+1,\text{EXT}}, p^h_{n+1}, p^h_{n+1,\text{EXT}})  =
		\left (\bm v^h, \bm h_{n+1}\right )_{\Gamma^{\mathrm N}_{n+1}}
		\\
		+
		&
		\left (\bm v^h, \, \rho (\bm u^h_* \cdot \bm n)_{-}\bm u^h_{n+1}\right )_{\Gamma^{\mathrm N}_{n+1}}
		+ \left ( \bm v^h, \rho \frac{\bm u^h_{n, \text{BDF}}}{\Delta t}\right )_{\Omega_{n}}, \;
		\\
		&
		\text{for all } (\bm v^h, q^h) \in \bm{\mathcal V}_{\bm 0}^h \times \mathcal Q^h,
		\; \text{for all } n \geq \sigma_t - 1.
	\end{aligned}
	\label{NS_ALE_RIIS_fullydiscretevmsles_semiimplicit_quasistatic}
\end{equation}
In Eq. \eqref{NS_ALE_RIIS_fullydiscretevmsles_semiimplicit_quasistatic}, $\Omega_{n+1}$ is the domain at time step $n+1$, as we detail in Section~\ref{SEC:geometric_numerical}. Moreover, the form $\mathcal S$ includes the stabilization and turbulence terms provided by the VMS-LES method:
\begin{equation}
	\allowdisplaybreaks
	\begin{aligned}
		\allowdisplaybreaks
		& \mathcal{S}(\bm v^h, q^h, \bm u^h_{n+1}, \bm u^{\text{ALE},h}_{n+1},  \bm u^h_{n+1,\text{EXT}}, p^h_{n+1}, p^h_{n+1,\text{EXT}})  =
		\\
		& {\left ( \rho\bm u^h_* \cdot \nabla \bm v^h + \nabla q^h , \, \tau_{\text M}(\bm u^h_*) \bm r_{\text M} (\bm u^h_{n+1}, p^h_{n+1})\right )_{\Omega_{n+1}}}
		\\
		+
		&
		{\left (\nabla \cdot \bm v^h, \tau_{\text C}(\bm u^h_*)r_{\text C}(\bm u^h_{n+1})\right)_{\Omega_{n+1}} }
		\\
		+ &{ \left(\rho \bm u^h_* \cdot (\nabla \bm v^h)^T ,  \tau_{\text M}(\bm u^h_*) \bm r_{\text M}(\bm u^h_{n+1}, p^h_{n+1}) \right )_{\Omega_{n+1}} }
		\\
		- &  {\left( \rho \nabla \bm v^h,  \tau_{\text M}(\bm u^h_*)\bm r_{\text {M}}(\bm u^h_*, p^h_{n+1,\text{EXT}}) \otimes  \tau_{\text M}(\bm u^h_{*})\bm r_{\text{M}}(\bm u^h_{n+1}, p^h_{n+1}) \right )_{\Omega_{n+1}}},
	\end{aligned}
\end{equation}
being
$\bm r_{\text M} (\bm u^h, p^h)$ and $r_{\text C}(\bm u^h)$ the strong residuals of (\ref{eq_ns}) and (\ref{eq_div}), defined respectively as:
\begin{subequations}
	\allowdisplaybreaks
	\begin{align}
		\allowdisplaybreaks
		\bm r_{\text M} (\bm u^h, p^h) = & \, \rho \frac{\widehat \partial \bm u^h}{\partial t} + \rho \left( \left( \bm{u}^h - \bm{u}^{\text{ALE}} \right) \cdot \nabla \right) \bm{u}^h
		+ \nabla p^h -\mu \Delta \bm u^h   \notag
		\\
		& + \sum_{k=1}^{m}\frac{R_k}{\varepsilon_k}\delta_{\Sigma_k, \varepsilon_k}(\varphi_k)\left(\bm u^h - \bm u^\mathrm{ALE} \right ), \label{eq:rm} \\
		r_{\text C}(\bm u^h)  = & \nabla \cdot \bm{u}^h, \label{eq:rc}
	\end{align}
\end{subequations}
where the temporal derivative in Eq. \eqref{eq:rm} is approximated with BDF method of the same order $\sigma_t$ adopted for the temporal discretization in Eq. \eqref{NS_ALE_RIIS_fullydiscretevmsles_semiimplicit_quasistatic}, and the convective term is treated semi-implicitly.  The stabilization parameters are chosen as in \cite{bazilevs_2007, DF_2015,Fedele_RIIS, Fumagalli2020}:
\begin{align*}
	\tau_{\text M}(\bm{u}^h_*) & =  \left( \frac{\sigma_t^2 \rho^2 }{\Delta t^2} +
	\rho^2 \, \bm{u}^h_* \cdot \mathfrak G\bm{u}^h_*  +
	C_r \mu^2 {\mathfrak G} :\mathfrak G + \sum_{k=1}^m\frac{R_k^2}{\varepsilon_k^2}\delta_{\Sigma_{k}, \varepsilon_k}^2(\varphi_k) \right)^{-\frac{1}{2}} \,,
	\\
	\tau_{\text C}(\bm{u}^h_*) & = \left ( \tau_{\mathrm M}(\bm{u}^h_*) \mathfrak g \cdot \mathfrak g \right )^{-1} \,,
	\label{tau}
\end{align*}
being $C_r=15\cdot{2^{r}}$ a constant obtained by an inverse inequality depending on the local velocity polynomial degree $r$ \cite{bazilevs_2007, DF_2015}, while $\mathfrak G$ and $\mathfrak g$ are  the metric tensor and metric vector, respectively (see \cite{bazilevs_2007}).

In addition, to avoid possible instabilities on the boundaries where Neumann BCs are set, we add the  backflow stabilization term $\left (\bm v^h, \, \rho (\bm u^h_* \cdot \bm n)_{-}\bm u^h_{n+1}\right )_{\Gamma^{\mathrm N}_{n+1}}$ in Eq.~\eqref{NS_ALE_RIIS_fullydiscretevmsles_semiimplicit_quasistatic}, being $(w)_- = \min \left (w, \, 0 \right )$. The advection velocity is treated semi-implicitly, coherently with the way we consider the non-linearities in the main formulation \cite{bertoglio2014tangential}.

\subsection{Space discretization of the geometric problem}
\label{SEC:geometric_numerical}

For the ALE lifting problem in Eq. \eqref{geometric_strong}, we introduce the infinite dimensional function space $\bm{\mathcal V}_{\widehat{\bm d}_\Gamma} : = \left\{ \widehat {\bm v} \in \left[ H^1(\widehat{\Omega})\right]^d: \widehat{\bm v} = \widehat{\bm d}_{{\Gamma}} \text{ on } \widehat \Gamma \right\}$
and the Galerkin formulation is expressed as: for any $n=0, \dots, N_t$,
\begin{equation}
	\text{ find } \widehat{\bm d}^h_n \in \bm{\mathcal V}_{\widehat{\bm d}_\Gamma}^h \text{ such that } \left ( \nabla \widehat{\bm d}^h_n, \, \nabla \widehat{\bm w}^h \right) = 0, \quad \text{ for all } \widehat{\bm w}^h \in \bm{\mathcal V}_{\bm 0}^h.
	\label{ALE_laplacian_Galerkin}
\end{equation}
The ALE velocity in the reference configuration is recovered by discretizing Eq. \eqref{uALE_dd_dt} with the Backward Euler Method, as
\begin{equation}
	\widehat{\bm u}^{\text{ALE}, h}_{n+1} = \frac{\widehat{\bm d}^h_{n+1} - \widehat{\bm d}^h_{n}}{\Delta t}.
	\label{uALE_dd_dt_discrete}
\end{equation}
The domain $\Omega_{n+1}$ and the mesh $\mathcal T^h_{n+1}$ are defined according to the ALE map $\mathcal{A}_{n+1}(\bm x) = \widehat{\bm x} + \widehat{\bm d}^h_{n+1}$ defined in Eq. \eqref{ALEmap}, namely
$\Omega_{n+1} = \mathcal{A}_{n+1}(\widehat{\Omega})$ and $ \mathcal T^h_{n+1}=\mathcal{A}_{n+1}( \widehat{\mathcal {T}^h}).$

\subsection{Time discretization of the 0D circulation model}
\label{SEC:0D_ERK}
We solve the system of ODEs in Eq. \eqref{0D_3D_compact} with a 4$^\text{th}$ order explicit Runge-Kutta method \cite{quarteroni2010numerical}. 
The time-step size $\Delta t$ employed for its numerical discretization is the same used for the BDF advancing scheme in the 3D problem.

\subsection{A segregated scheme for the 3D-0D coupling}
\label{SEC:numerics_0D_3D}

\begin{figure}[t!]
	\centering
	\includegraphics[trim={5 5 5 5},
	clip,
	width=\textwidth]
	{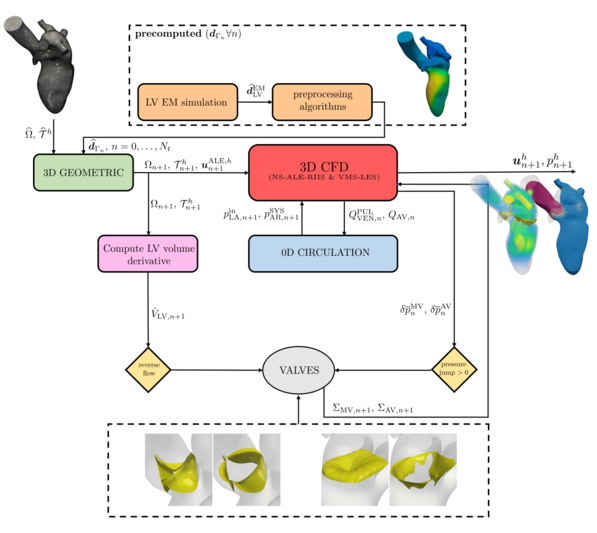}
	\caption{Sketch of the algorithm for the CFD simulation of the LH. The coupling between the CFD and the circulation problem is solved via a segregated numerical scheme.}
	\label{algo_scheme}
\end{figure}

\begin{algorithm}[t!]
	\caption{Segregated scheme for the 3D-0D coupled problem}
	\label{algo_3D_0D}
	\begin{algorithmic}
		\STATE{Initialization, $n= 0$}.
		\STATE{Compute $\widehat{\bm d}_{\Gamma_n}$ for $n = 0, \dots, N_\mathrm{HB}$}.
		\WHILE{$n \leq N_t$}
		\STATE{\vspace{-0.35cm}\begin{align*}
				& \text{Solve geometric problem.}
				\\
				& \text{Compute LV volume derivative and pressure jumps.}
				\\
				& \text{Update valves leaflet position.}
				\\
				& \text{Solve circulation}: \left ( \bm c_{n+1}, \, \tilde{\bm c}_{n+1}\right )^T= \textsc{Circulation}(Q_{\text{VEN},n}^{\text{PUL}}, \, Q_{\text{AV},n}).
				\\
				& \text{Compute interface data (0D $\to$ 3D)}: \left ( \bm c_{n+1}, \, \tilde{\bm c}_{n+1}\right ) \to p^\mathrm{in}_{\mathrm{LA}, n+1}, \, p^\mathrm{SYS}_{\mathrm{AR}, n+1.}
				\\
				& \text{Solve fluid dynamics:} \left ( \bm u^h_{n+1}, \, p^h_{n+1}\right )^T = \textsc{FluidDynamics}(p^\mathrm{in}_{\mathrm{LA}, n+1}, \, p^\mathrm{SYS}_{\mathrm{AR}, n+1}).
				\\
				&\text{Compute interface data  (3D $\to$ 0D)}: \bm u^h_{n+1} \to Q_{\text{VEN},n+1}^{\text{PUL}}, \, Q_{\text{AV},n+1}.
				\\
				& n \leftarrow n+1.
		\end{align*}}
		\ENDWHILE
	\end{algorithmic}
\end{algorithm}

We show the overall numerical scheme for the coupled problem for the CFD simulation of the LH in Algorithm~\ref{algo_3D_0D}. A graphical representation of the whole algorithm is given in Figure~\ref{algo_scheme}.

After initialization, we compute the displacement field for a representative heartbeat ($N_\mathrm{HB}$ time steps), as explained in Section~\ref{SEC:displacement_modeling}. We interpolate the displacement field in time on the temporal grid required by the CFD problem via splines. Furthermore, we assume the displacement field to be periodic for all the heartbeats simulated.  For each time step, we solve the geometric problem as explained in Section~\ref{SEC:geometric_numerical}; we compute the LV volume derivative and the pressure jumps across each valve, which represent the indicators employed to determine the valves status, as we discussed in Section~\ref{SEC:valves_modeling}. Given the inlet and outlet flowrates from the previous time step ($Q_{\text{VEN},n}^{\text{PUL}},\, Q_{\text{AV},n}$), we solve the discretized circulation problem (see Section~\ref{SEC:0D_ERK}). The inlet atrial and systemic arterial pressures, computed in the 0D model, represent the interface conditions from the 0D to the 3D model. Specifically, the pressures $ p^\mathrm{in}_{\mathrm{LA}, n+1}, \, p^\mathrm{SYS}_{\mathrm{AR}, n+1}$ are used as Neumann data for the inlet and outlet BCs, respectively (as we explain in Section~\ref{SEC:coupling_0D_3D}). Once the fluid dynamics problem is solved as explained in Section~\ref{SEC:numerical_NS_ALE_RIIS}, we integrate the normal relative velocity to compute $Q_{\text{VEN},n+1}^{\text{PUL}}, \, Q_{\text{AV},n+1}$ as reported in Eq.~\eqref{continuity_Q}. The flowrates stand as interface data from the 3D to the 0D problem for the following time step. Thus, to advance the coupled problem in time, we use a segregated numerical scheme in which we solve, with the same time step size, the circulation and the fluid dynamics problems. The restriction on the time-step size is given by the CFD problem, being also the more expensive part to be solved. From our numerical tests, we found that the usage of the same $\Delta t$ for the two problems guarantees that the solution is always stable, motivating hence the usage of a segregated numerical scheme.

\section{Numerical simulations}
\label{SEC:numerical_results}
\begin{figure}[t]
	\centering
	\includegraphics[trim={0 0 0 0},clip, width=\textwidth]{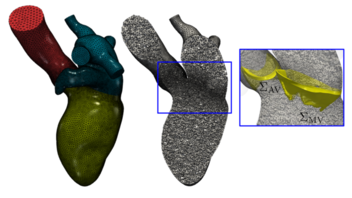}
	\caption{The LH tetrahedral mesh made of three conforming meshes for the LA, LV and AA subdomains; a clip of the mesh showing the local mesh refinement near the MV and AV.}
	\label{LH_mesh}
\end{figure}

\begin{table}[t!]
	\centering
	\begin{tabular}{c|c|c|c|c|c|c|c|c}
		$\rho$ & $\mu$ &  $\bm u_0$ & \multicolumn{2}{c|}{ $R_{\mathrm k}$}  & \multicolumn{2}{c|}{$\varepsilon_\mathrm k$}  & $T_f$& $T_\mathrm{HB}$  \\
		\footnotesize{$[\mathrm{kg/m}^3]$} & 	\footnotesize{$[\mathrm{kg/(m}\cdot\mathrm{s})]$} &
		\footnotesize{$[\mathrm m / \mathrm s]$} & 	\multicolumn{2}{c|}{ \footnotesize{[kg/(m$^2\cdot$s)]}}  & \multicolumn{2}{c|}{\footnotesize{[mm]}} & \footnotesize{[s]} & \footnotesize{[s]}
		\\
		& & & \footnotesize{MV} & \footnotesize{AV} & \footnotesize{MV} & \footnotesize{AV} & &
		\\
		\hline
		\hline
		$1.06 \cdot 10^3$ & $3.5 \cdot 10^{-3}$ & $\bm 0$ & $10^4$  & $10^4$ & 0.6 & 0.6 & 2.0 & 1.0  \\
	\end{tabular}
	\\ \vspace{0.5cm}
	\begin{tabular}{c|c|c|c|c|c|c|c|c}
		\multicolumn{3}{c|}{$h$} & {cells} & \multicolumn{3}{|c|}{{DOFs ($ \mathbb{P}_1 - \mathbb{P}_1 $)}} & BDF & $\Delta t$\\
		\multicolumn{3}{c|}{\footnotesize{[mm]}} & [-] &  \multicolumn{3}{|c|}{[-]} & \footnotesize{[-]}& \footnotesize[s]\\
		\footnotesize{min} & \footnotesize{avg} & \footnotesize{max}  && \footnotesize{$\bm u$} & \footnotesize{$p$} & \footnotesize{total} & &\\
		\hline
		\hline
		0.4 & 1.2 & 4.1 & 1'627'795 & 806'295 & 268'765 & 1'075'060 & 1 &$2.5 \cdot 10^{-4}$\\
	\end{tabular}
	\caption{Parameters for the setup of the numerical simulations.}
	\label{table_setup}
\end{table}

We perform numerical simulations on a LH mesh built from the Zygote Solid 3D heart model \cite{zygote}. To generate the mesh, we use $\mathsf{vmtk}$ \cite{vmtk} and the pre-processing tools for cardiac geometries proposed in \cite{FQ_2021}. We mesh each portion of the domain (LA, LV and AA) separately, with a non-uniform mesh size, refining the mesh near the valves by an algorithm based on the closest-point distances \cite{FQ_2021}. The three meshes are then connected in a conforming fashion to get the tetrahedral mesh in Figure \ref{LH_mesh}.  Specifically, to connect the three meshes, we use the \textit{mesh-connection} algorithm introduced in \cite{FQ_2021}. We use linear FE spaces for velocity and pressures ($ \mathbb{P}_1 - \mathbb{P}_1 $); as time integration scheme, we employ the BDF method of order 1, and a time-step size $\Delta t = 2.5 \cdot 10^{-4} \, \mathrm{s}$.  We summarize the parameters we use for the setup of our numerical simulations  in Table~\ref{table_setup}.

For the RIIS method, following arguments of \cite{Fedele_RIIS, Fumagalli2020}, we set $R_\text{k}=10^4$ kg/(m$^2 \cdot$s) and $\varepsilon_{\text{k}}= 0.6 $ mm, with $\mathrm k = \mathrm{MV, \, AV}$. Our choice of $\varepsilon_{\text{k}}$ allows to have a physiologic representation of the valves' leaflets in a healthy subject and, with our setting of $R_\text{k}$, to avoid flow penetration across the immersed surfaces \cite{Fedele_RIIS}.  In order to accurately represent the valves by means of the RIIS method, the mesh size $h$ must be chosen small enough in the immersed surfaces region. Specifically, the minimum value of $h$ should be set such that $\varepsilon_k$ is at least 1.5 times $h$ \cite{Fedele_RIIS}. Furthermore, we found that the condition number of the linear system associated to the FE discretization of the NS-ALE-RIIS equations becomes larger as $R_\mathrm{k}/\varepsilon_\mathrm{k}$ increases. Thus, for $R_{\mathrm k}$ we choose the minimum value that guarantees impervious valves.
As displayed in Figure \ref{LH_mesh}, we have refined the mesh in the valves region to accurately represent the leaflets as described by the RIIS method with our choices of $\varepsilon_{\text{k}}$. Since we neglect the isovolumetric phases, we cannot set physiological opening and closing times of the valves, thus we open and close them instantaneously (i.e. in one time step). The physical parameters for  blood are density $\rho = 1.06 \cdot 10^3$ kg/m$^3$ and dynamic viscosity $\mu=3.5\cdot 10^{-3}$ kg/(m$\cdot$s); we start our simulation from a zero velocity initial condition $\bm u_0 = \bm 0$.

The mathematical models and numerical methods described in Sections \ref{SEC:0D_CIRCULATION}, \ref{SEC:3D_CFD},  \ref{SEC:coupling_0D_3D} and \ref{SEC:numerical_methods} have been implemented in $\mathsf{life^x}$ \cite{lifex, lifex_arxiv}, a high-performance \texttt{C++} library developed within the iHEART project\footnote{\textit{iHEART - An Integrated Heart model for the simulation of the cardiac function}, European Research Council (ERC) grant agreement No 740132, P.I. Prof A. Quarteroni, 2017-2022.}, mainly focused on cardiac simulations, and based on the $\mathsf{deal.II}$ FE core \cite{dealii}. Numerical simulations are run in a parallel framework\footnote{Numerical simulations were run either on the cluster iHEART (Lenovo SR950 8 x 24-Core Intel Xeon Platinum 8160, 2100 MHz and 1.7TB RAM) available at MOX, Dipartimento di Matematica, Politecnico di Milano and on the cluster GALILEO supercomputer (IBM NeXtScale cluster, 1022 nodes (Intel Broadwell), 2 x 18-Cores Intel Xeon E5-2697 v4 at 2.30 GHz, 36 cores/node, 26.572 cores in total with 128 GB/node) by CINECA.
}. The linear system arising from Eq. \eqref{NS_ALE_RIIS_fullydiscretevmsles_semiimplicit_quasistatic} is preconditioned with the aSIMPLE preconditioner \cite{Deparis_2014}, and each of its blocks are preconditioned with an algebraic multigrid preconditioner based on Trilinos \cite{trilinos}. The linear system is then solved at each time step with the GMRES method.

In the following, we  present our numerical results for the case of a healthy LH.

\begin{figure}[t]
	\centering
	\includegraphics[trim={1cm 1 1cm 1},clip, width=\textwidth]{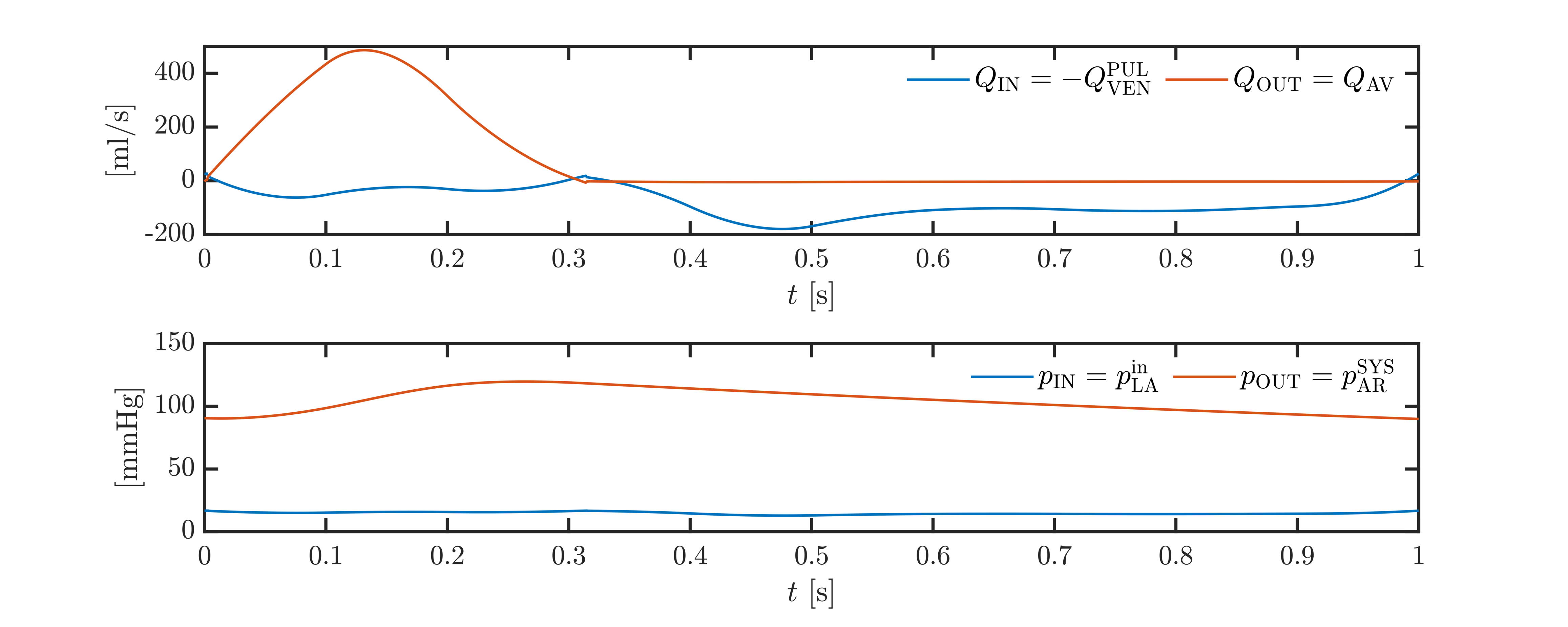}
	\caption{Flowrates (top) and pressures (bottom) at the interfaces of the 3D-0D model.}
	\label{flowrates_pressure_bcs}
\end{figure}


\begin{figure}[!t]
	\centering
	\includegraphics[trim={4 4 4 4},clip, width=\textwidth]{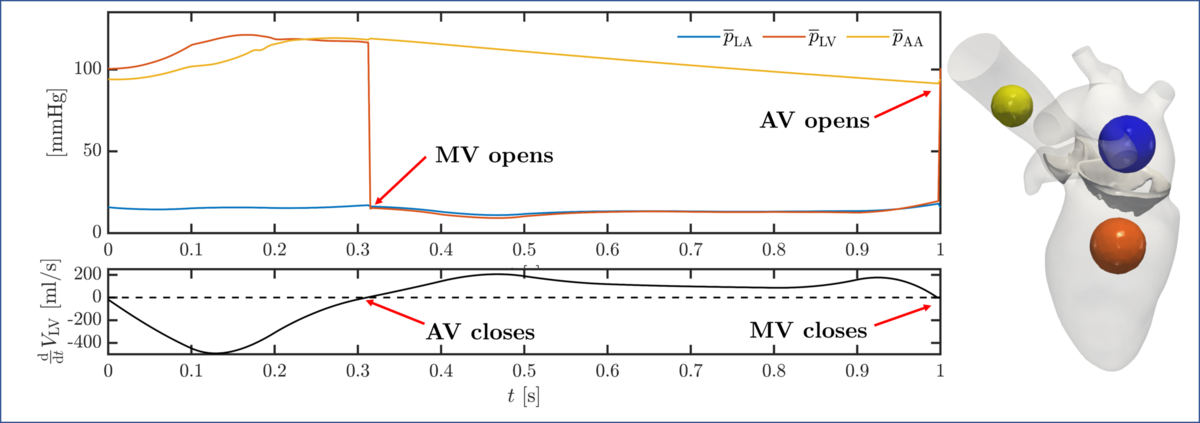}
	\caption{Flow properties to determine opening and closure of valves. Top: average pressure in LA, LV and AA; bottom: time derivative of LV volume.  Average pressures in the chambers are computed in the control volumes in the LH displayed on the right.}
	\label{pressure_chambers_volume_derivative}
\end{figure}

\begin{figure}[t!]
	\begin{subfigure}{0.19\textwidth}
		\centering
		\includegraphics[trim={6 6 6 6},clip, width=1.1\textwidth]{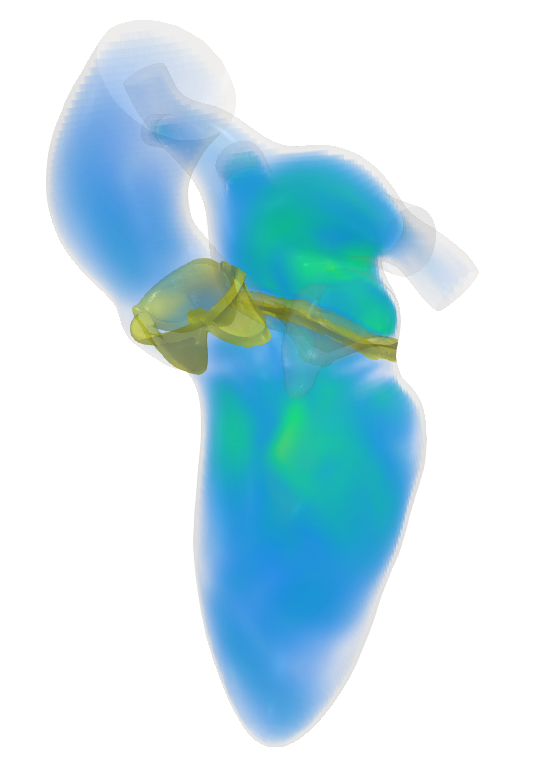}
		\caption{$t = 0.0$ s}
		\label{velocity.0000}
	\end{subfigure}
	\begin{subfigure}{0.19\textwidth}
		\centering
		\includegraphics[trim={6 6 6 6},clip, width=1.1\textwidth]{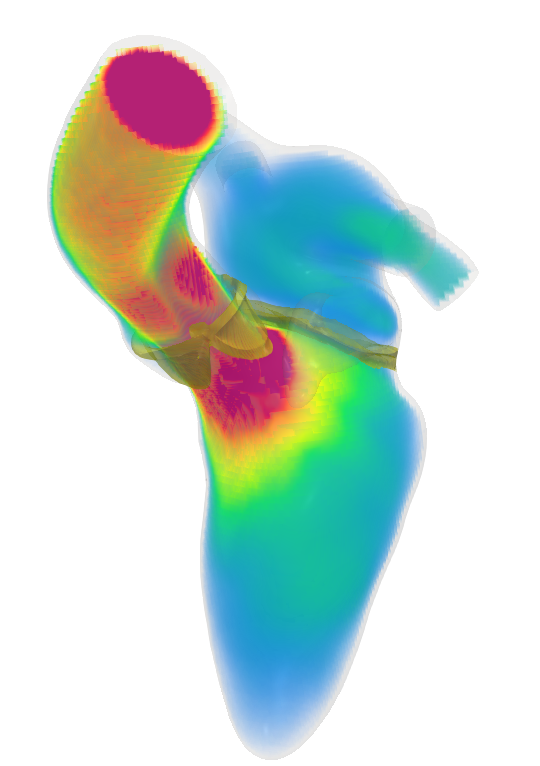}
		\caption{$t = 0.1$ s}
		\label{velocity.0001}
	\end{subfigure}
	\begin{subfigure}{0.19\textwidth}
		\centering
		\includegraphics[trim={6 6 6 6},clip, width=1.1\textwidth]{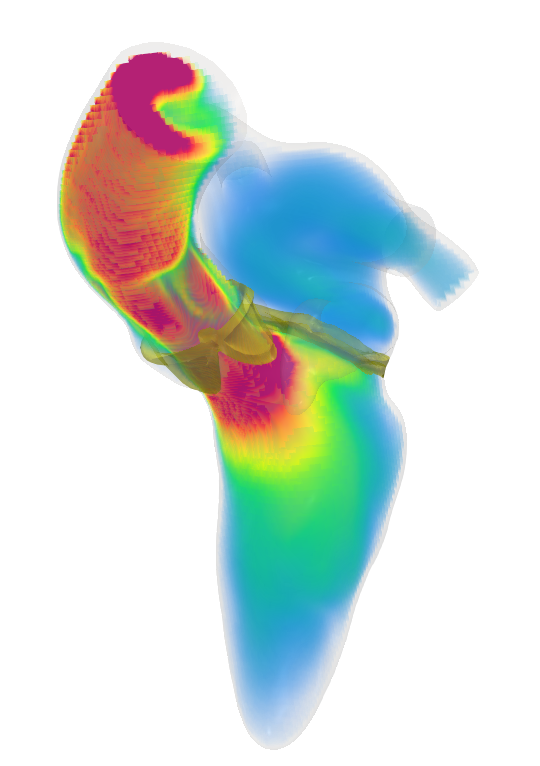}
		\caption{$t = 0.2$ s}
		\label{velocity.0002}
	\end{subfigure}
	\begin{subfigure}{0.19\textwidth}
		\centering
		\includegraphics[trim={6 6 6 6},clip, width=1.1\textwidth]{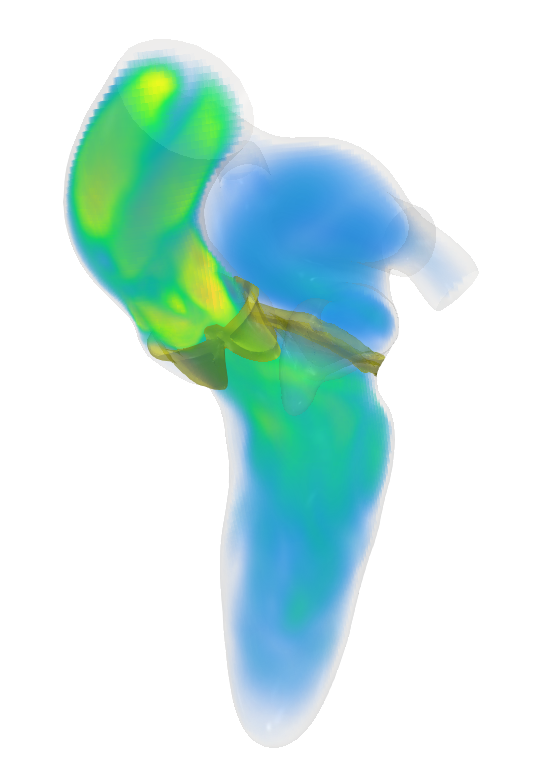}
		\caption{$t = 0.3$ s}
		\label{velocity.0003}
	\end{subfigure}
	\begin{subfigure}{0.19\textwidth}
		\centering
		\includegraphics[trim={6 6 6 6},clip, width=1.1\textwidth]{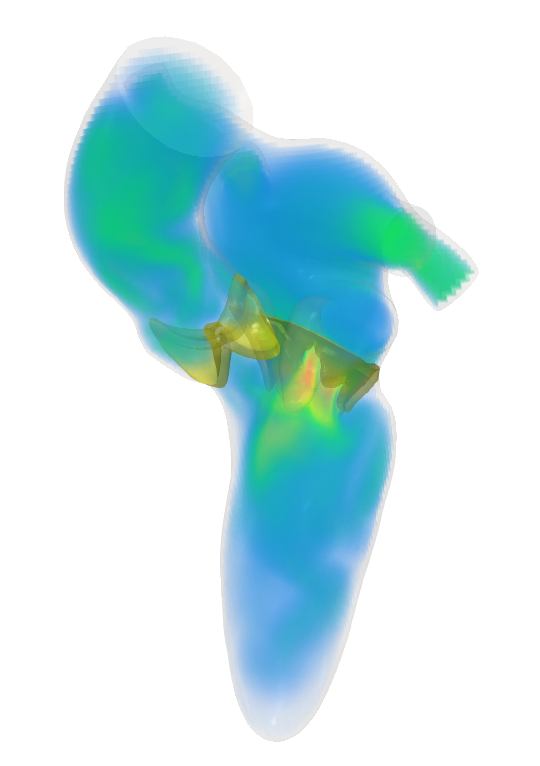}
		\caption{$t = 0.4$ s}
		\label{velocity.0004}
	\end{subfigure}
	\\
	\begin{subfigure}{0.19\textwidth}
		\centering
		\includegraphics[trim={6 6 6 6},clip, width=1.1\textwidth]{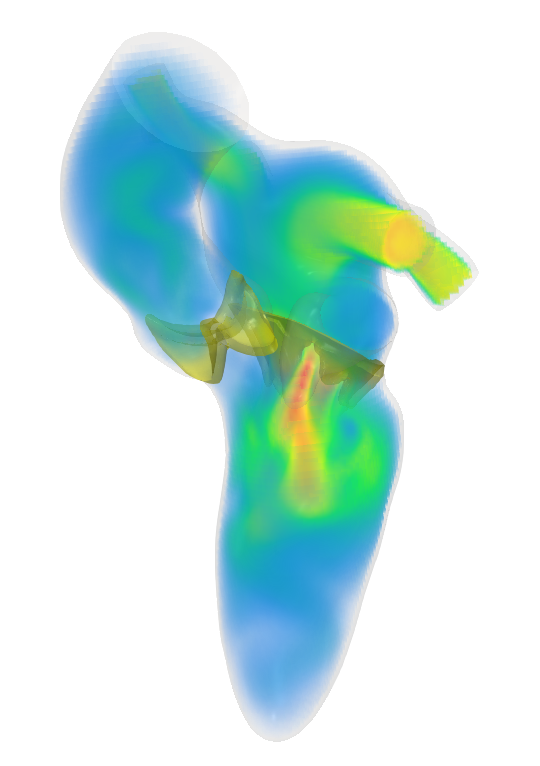}
		\caption{$t = 0.5$ s}
		\label{velocity.0005}
	\end{subfigure}
	\begin{subfigure}{0.19\textwidth}
		\centering
		\includegraphics[trim={6 6 6 6},clip, width=1.1\textwidth]{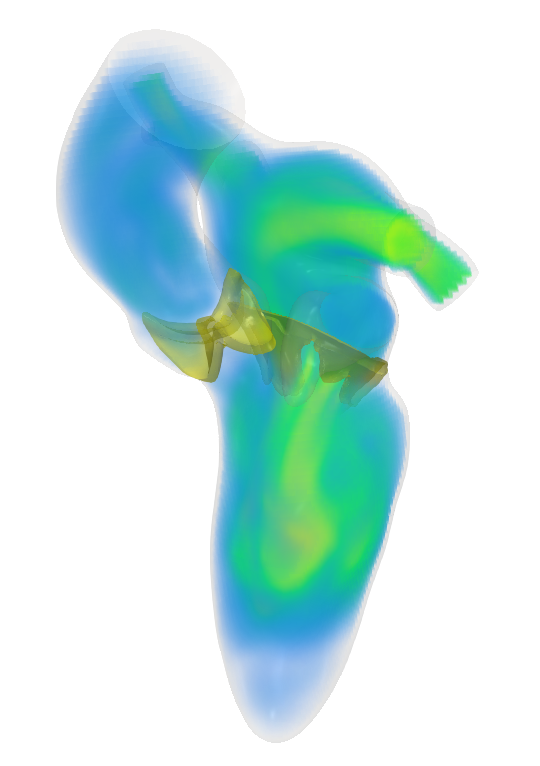}
		\caption{$t = 0.6$ s}
		\label{velocity.0006}
	\end{subfigure}
	\begin{subfigure}{0.19\textwidth}
		\centering
		\includegraphics[trim={6 6 6 6},clip, width=1.1\textwidth]{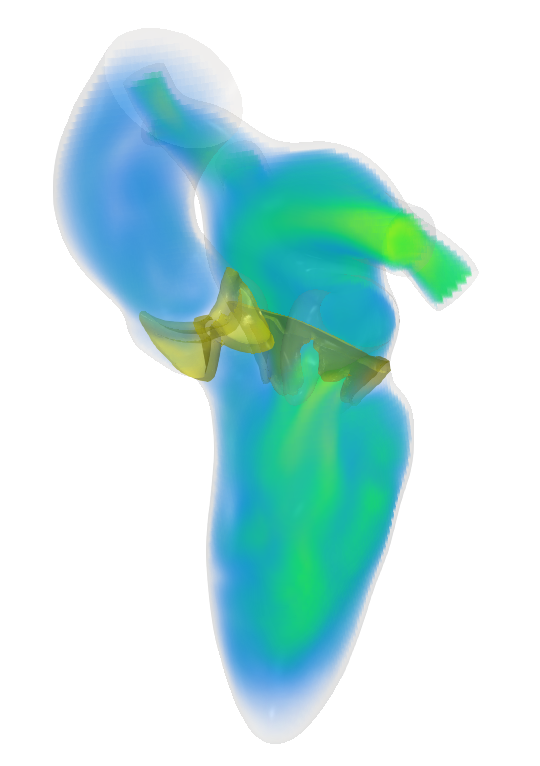}
		\caption{$t = 0.7$ s}
		\label{velocity.0007}
	\end{subfigure}
	\begin{subfigure}{0.19\textwidth}
		\centering
		\includegraphics[trim={6 6 6 6},clip, width=1.1\textwidth]{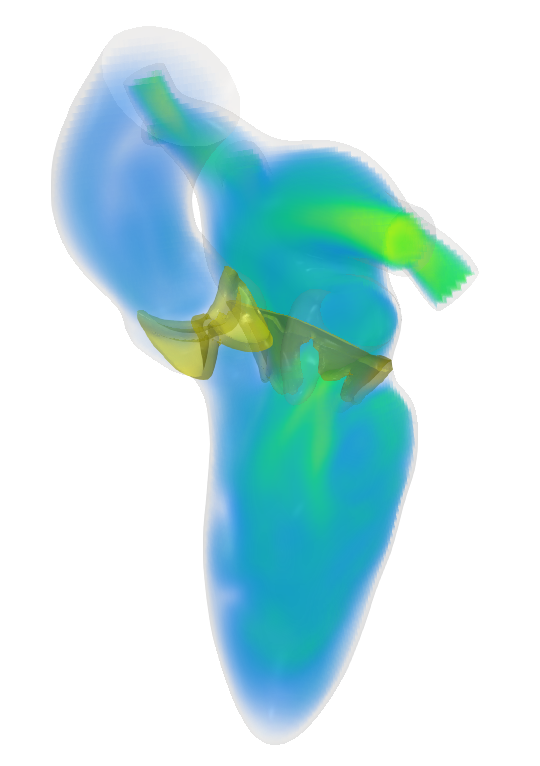}
		\caption{$t = 0.8$ s}
		\label{velocity.0008}
	\end{subfigure}
	\begin{subfigure}{0.19\textwidth}
		\centering
		\includegraphics[trim={6 6 6 6},clip, width=1.1\textwidth]{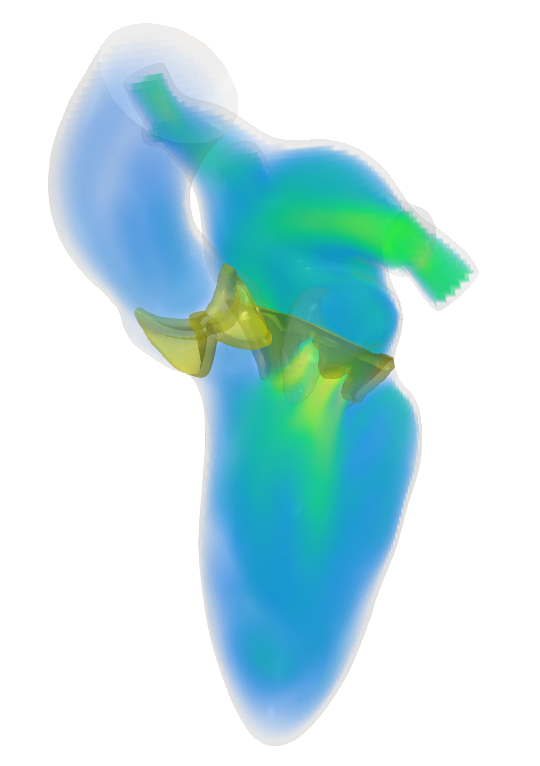}
		\caption{$t = 0.9$ s}
		\label{velocity.0009}
	\end{subfigure}
	\\
	\begin{subfigure}{\textwidth}
		\centering
		\includegraphics[trim={6 6 6 6},clip, width=0.2\textwidth]{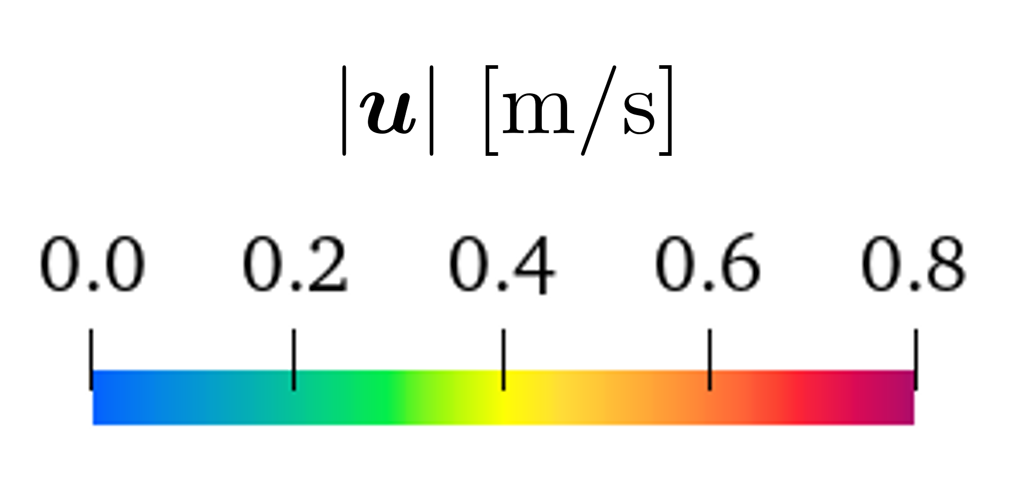}
	\end{subfigure}
	\caption{Volume rendering of velocity magnitude during the whole heartbeat.}
	\label{clip_velocity}
\end{figure}

\begin{figure}[t!]
	\begin{subfigure}{0.19\textwidth}
		\centering
		\includegraphics[trim={6 6 6 6},clip, width=1.1\textwidth]{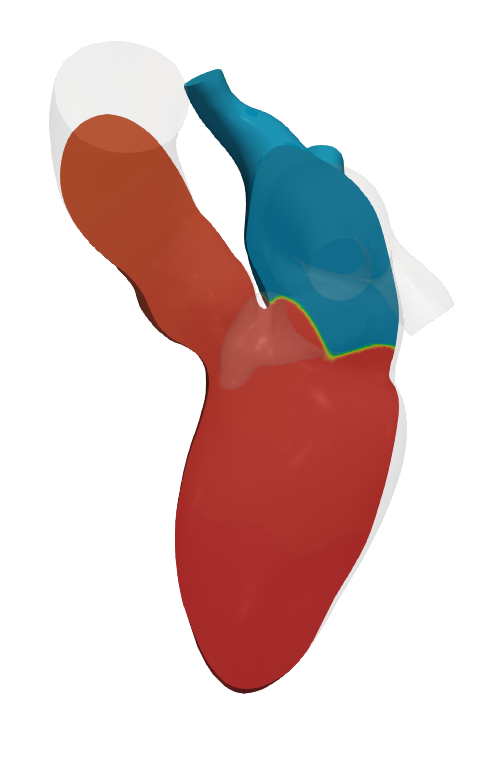}
		\caption{$t = 0.0$ s}
		\label{pressure.0000}
	\end{subfigure}
	\begin{subfigure}{0.19\textwidth}
		\centering
		\includegraphics[trim={6 6 6 6},clip, width=1.1\textwidth]{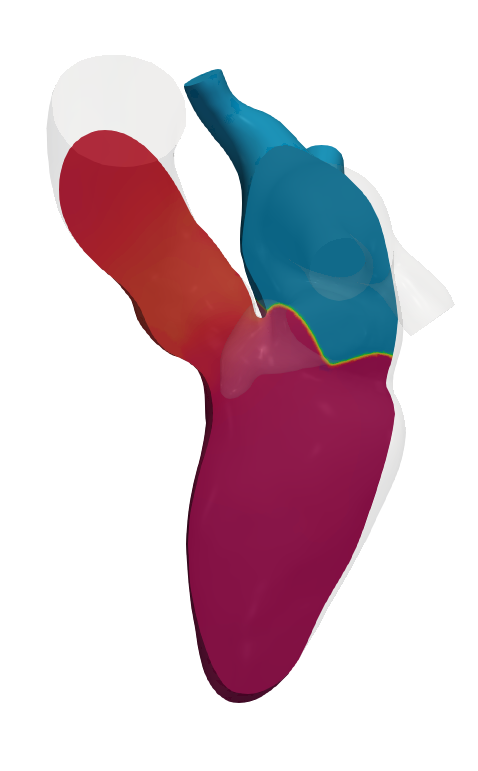}
		\caption{$t = 0.1$ s}
		\label{pressure.0001}
	\end{subfigure}
	\begin{subfigure}{0.19\textwidth}
		\centering
		\includegraphics[trim={6 6 6 6},clip, width=1.1\textwidth]{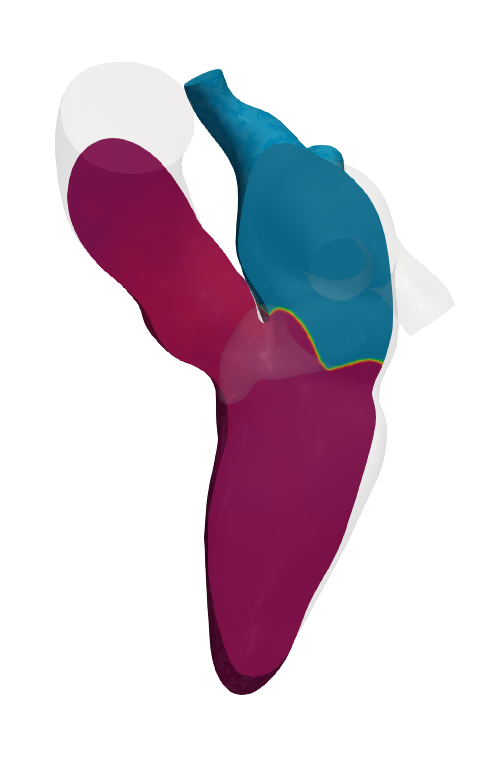}
		\caption{$t = 0.2$ s}
		\label{pressure.0002}
	\end{subfigure}
	\begin{subfigure}{0.19\textwidth}
		\centering
		\includegraphics[trim={6 6 6 6},clip, width=1.1\textwidth]{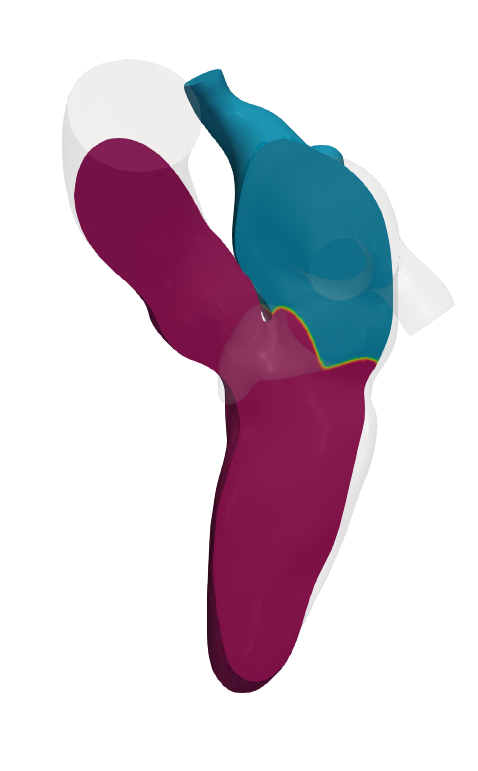}
		\caption{$t = 0.3$ s}
		\label{pressure.0003}
	\end{subfigure}
	\begin{subfigure}{0.19\textwidth}
		\centering
		\includegraphics[trim={6 6 6 6},clip, width=1.1\textwidth]{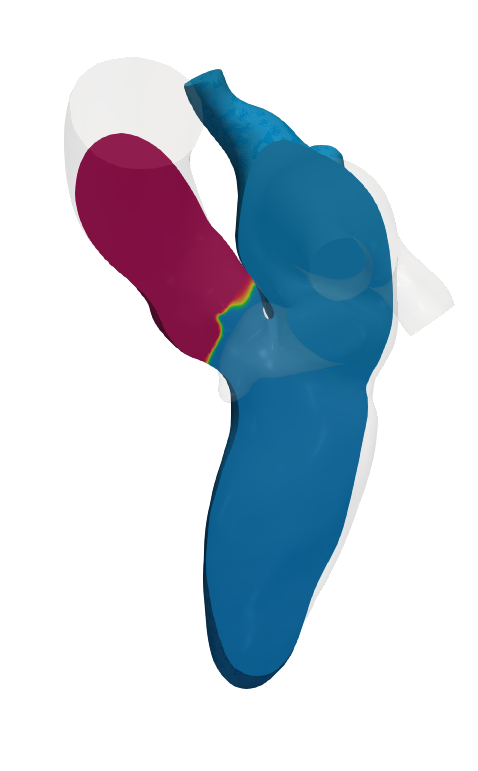}
		\caption{$t = 0.4$ s}
		\label{pressure.0004}
	\end{subfigure}
	\\
	\begin{subfigure}{0.19\textwidth}
		\centering
		\includegraphics[trim={6 6 6 6},clip, width=1.1\textwidth]{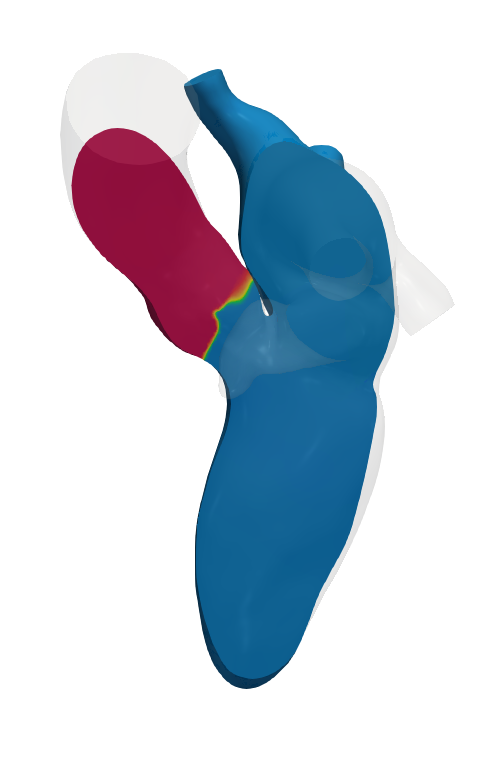}
		\caption{$t = 0.5$ s}
		\label{pressure.0005}
	\end{subfigure}
	\begin{subfigure}{0.19\textwidth}
		\centering
		\includegraphics[trim={6 6 6 6},clip, width=1.1\textwidth]{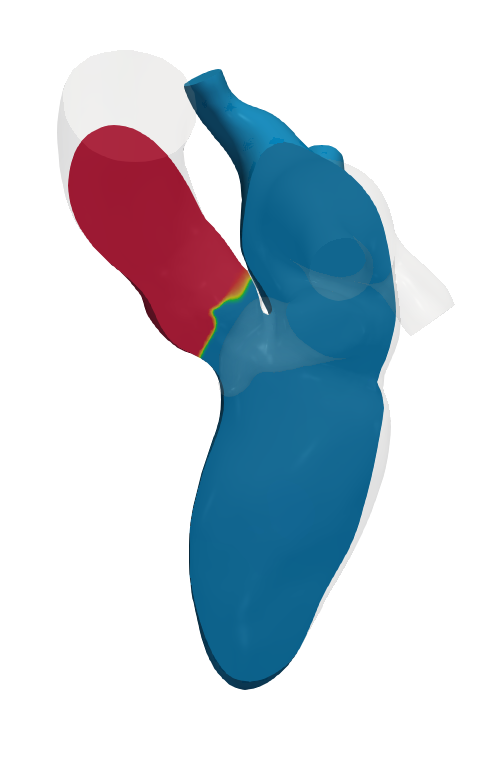}
		\caption{$t = 0.6$ s}
		\label{pressure.0006}
	\end{subfigure}
	\begin{subfigure}{0.19\textwidth}
		\centering
		\includegraphics[trim={6 6 6 6},clip, width=1.1\textwidth]{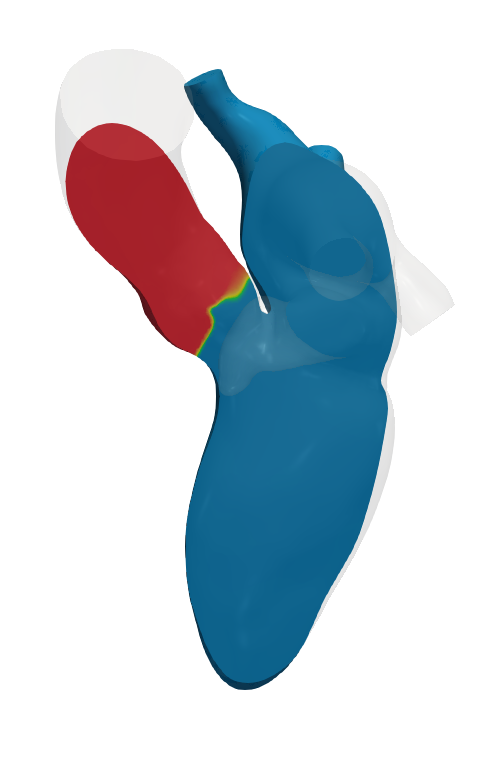}
		\caption{$t = 0.7$ s}
		\label{pressure.0007}
	\end{subfigure}
	\begin{subfigure}{0.19\textwidth}
		\centering
		\includegraphics[trim={6 6 6 6},clip, width=1.1\textwidth]{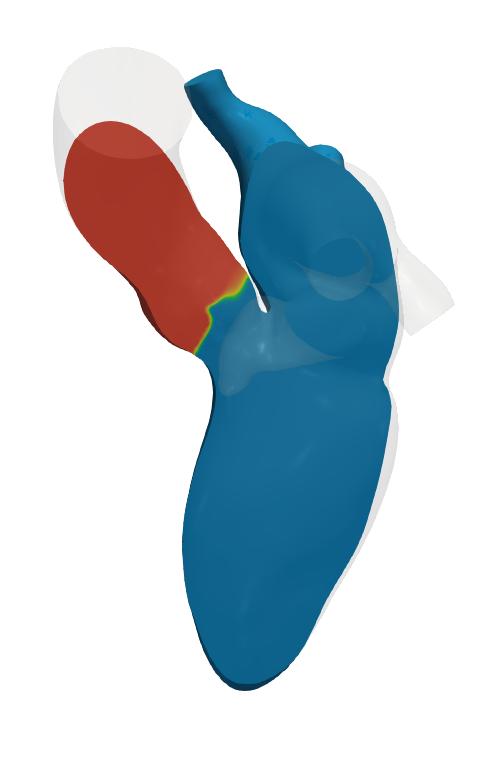}
		\caption{$t = 0.8$ s}
		\label{pressure.0008}
	\end{subfigure}
	\begin{subfigure}{0.19\textwidth}
		\centering
		\includegraphics[trim={6 6 6 6},clip, width=1.1\textwidth]{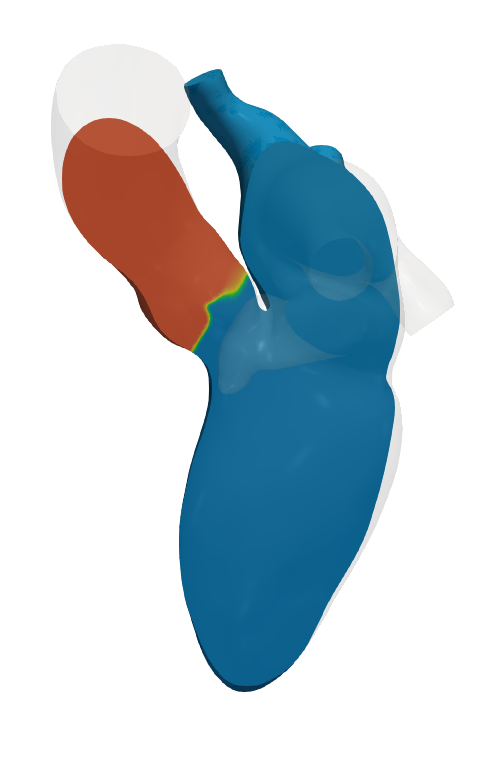}
		\caption{$t = 0.9$ s}
		\label{pressure.0009}
	\end{subfigure}
	\\
	\begin{subfigure}{\textwidth}
		\centering
		\includegraphics[trim={6 6 6 6},clip, width=0.2\textwidth]{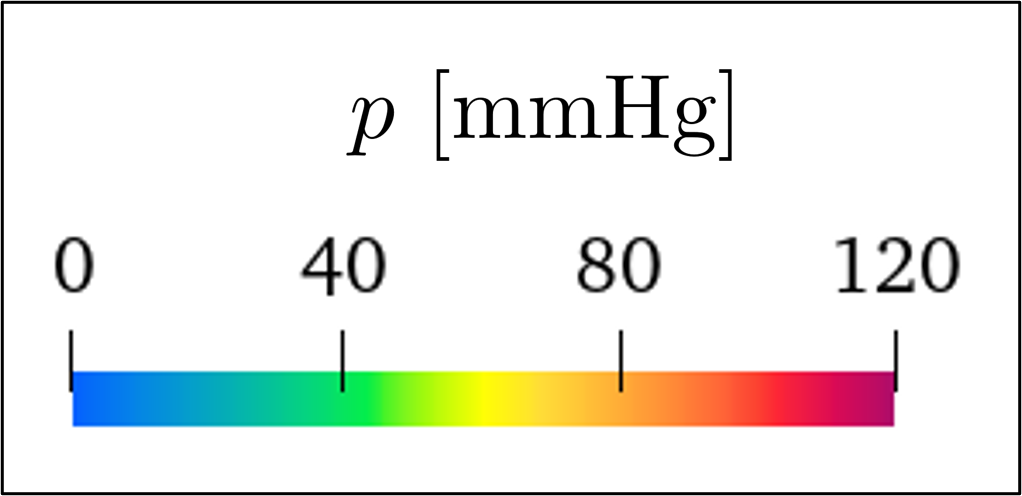}
	\end{subfigure}
	\caption{Pressure on a clip in the LV apico-basal direction during the whole heartbeat.}
	\label{clip_pressure}
\end{figure}

\begin{figure}[t!]
	\begin{subfigure}{0.16\textwidth}
		\centering
		\includegraphics[trim={1cm 1cm 4 5cm},clip, width=\textwidth]{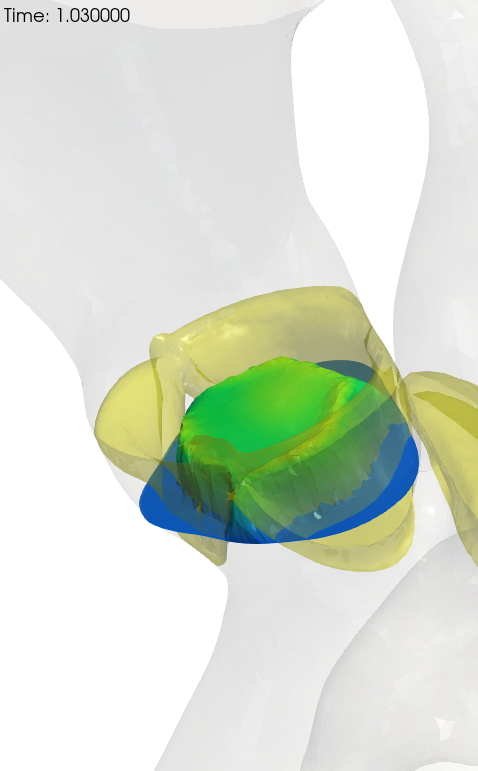}
		\caption{$t = 0.03$ s}
		\label{velocity_av.0001}
	\end{subfigure}
	\begin{subfigure}{0.16\textwidth}
		\centering
		\includegraphics[trim={1cm 1cm 4 5cm},clip, width=\textwidth]{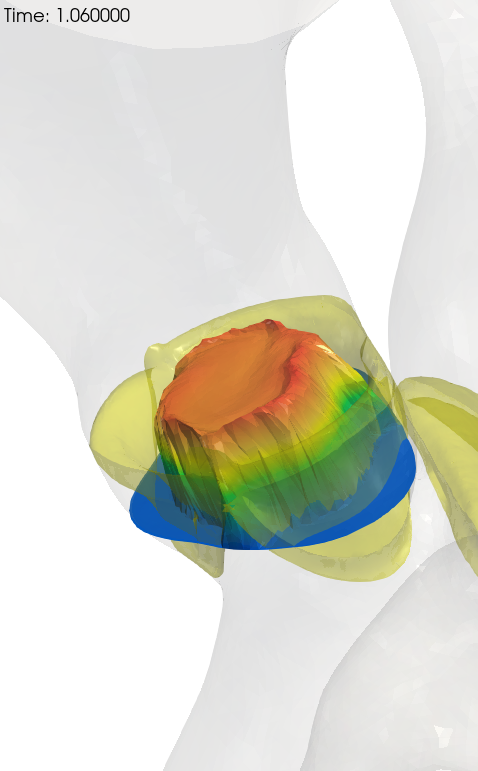}
		\caption{$t = 0.06$ s}
		\label{velocity_av.0002}
	\end{subfigure}
	\begin{subfigure}{0.16\textwidth}
		\centering
		\includegraphics[trim={1cm 1cm 4 5cm},clip, width=\textwidth]{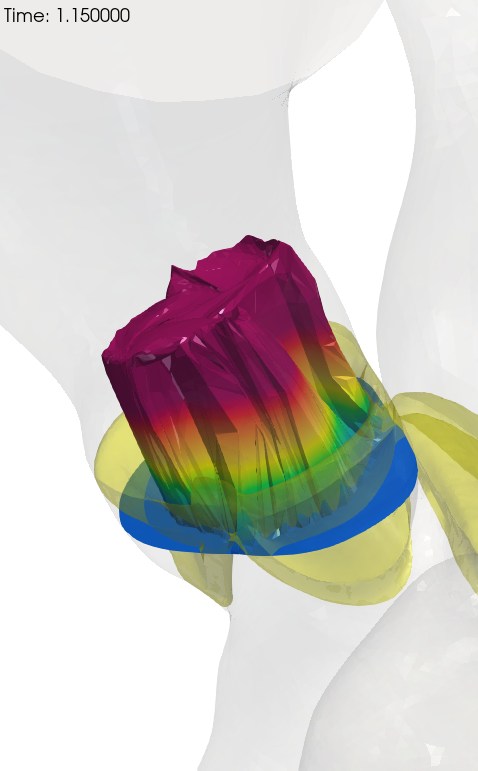}
		\caption{$t = 0.15$ s}
		\label{velocity_av.0005}
	\end{subfigure}
	\begin{subfigure}{0.16\textwidth}
		\centering
		\includegraphics[trim={1cm 1cm 4 5cm},clip, width=\textwidth]{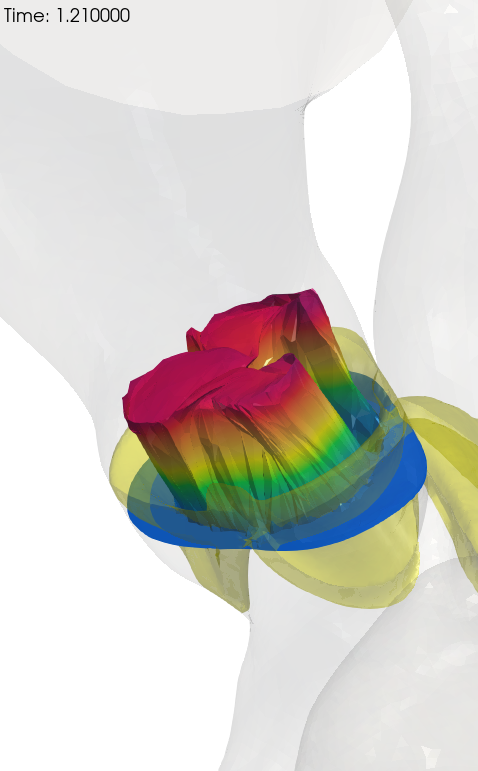}
		\caption{$t = 0.21$ s}
		\label{velocity_av.0007}
	\end{subfigure}
	\begin{subfigure}{0.16\textwidth}
		\centering
		\includegraphics[trim={1cm 1cm 4 5cm},clip, width=\textwidth]{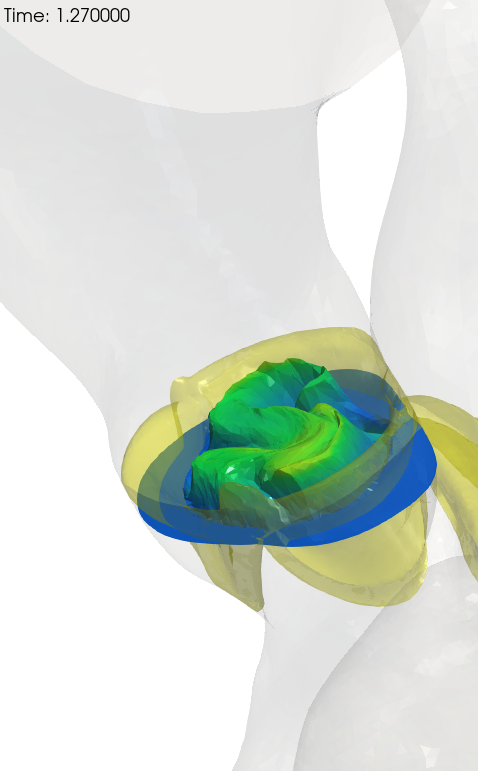}
		\caption{$t = 0.27$ s}
		\label{velocity_av.0009}
	\end{subfigure}
	\\
	\begin{subfigure}{\textwidth}
		\centering
		\includegraphics[trim={4 4 4 4},clip, width=0.2\textwidth]{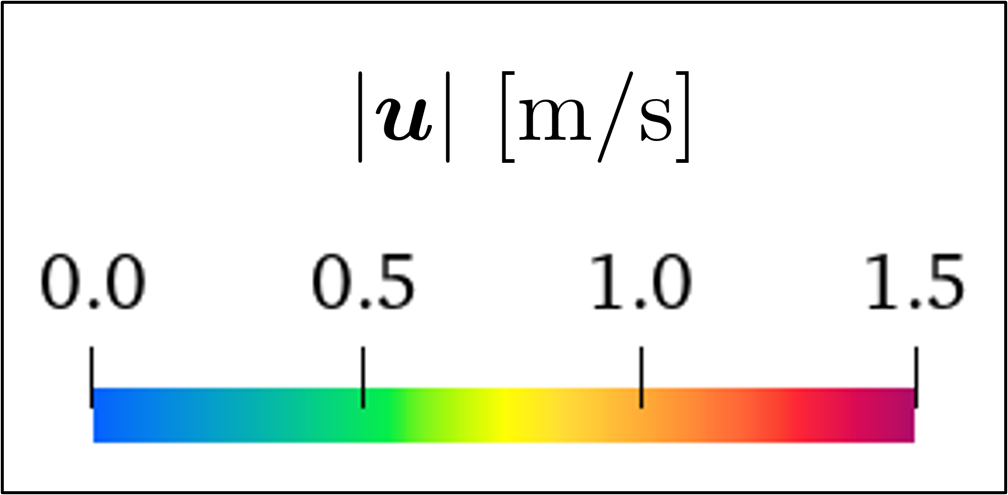}
	\end{subfigure}
	\caption{AV section warped by velocity during the ejection phase. In trasparency: the AV represented by the RIIS method.}
	\label{av_velocity}
\end{figure}

\begin{figure}[!t]
	\centering
	\includegraphics[trim={5 5 5 5},clip, width=\textwidth]{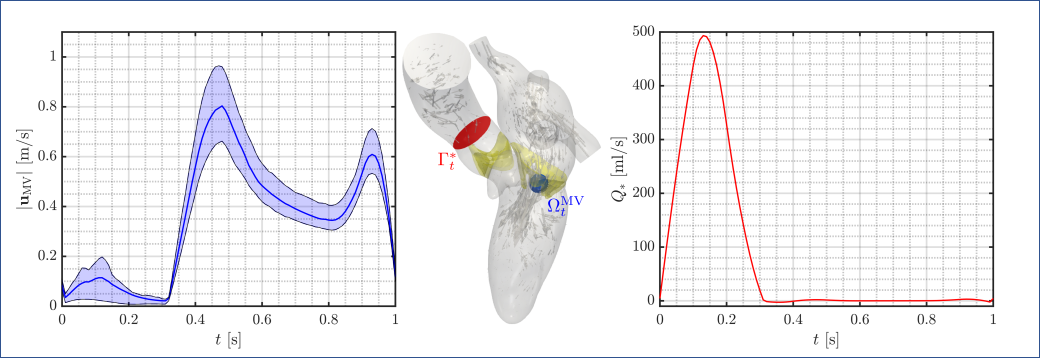}
	\caption{Left: mean, first and third quartile of the space distribution of the velocity magnitude $\bm u_\mathrm{MV}$ in a control volume $\Omega_t^\mathrm{MV}$ downstream the MV section. Right: flowrate computed on a section downwind the AV ($\int_{\Gamma_t^*} \bm u \cdot \bm n$) }
	\label{validation_MV_velocity}
\end{figure}

We simulate two heartbeats ($T_f = 2$ s) of period $T_\text{HB}=1$ s and we discard the first beat in order to disregard the effects of the unphysical initial condition $\bm u_0=0$.  This simulation required about 183 hours on a node with 56 cores. We show numerical results on the last heartbeat and we shift the temporal domain in $(0, T_\text{HB})$.
In Figure \ref{flowrates_pressure_bcs}, we report flowrates and pressures at the inlet and outlet sections of the 3D domain, which are the interfaces of the 3D-0D coupling. In Figure \ref{pressure_chambers_volume_derivative}, we show the average pressure in LA, LV and AA computed by averaging the pressure field in control volumes in each chamber, as explained in Section \ref{SEC:valves_modeling}. On the bottom box, we report the LV volume time-derivative, which is the indicator adopted to check for reversed flow across a valve and hence to determine the valves closing. Considering that we open and close the valves istantaneously, the pressure exhibits an oscillation immediately after this moments. From our numerical simulations, we found that closing the cardiac valves by checking for a reverse flow condition (instead of a negative pressure jump) completely prevents the pressure peaks after valve closing. However, the peak at the opening of the valve is still present, and we believe it might be related to the fact the valve opens in one time step, instead of respecting a dynamic given by the interaction between the leaflets and the fluid. For visualization purposes, the pressure fluctuations at the opening of the valves, lasting for maximum four time steps, are removed from the picture.

The previous simulated heartcycle ends with the closing of the MV  when the total flow through its section changes sign; the pressure in the LV starts rising and the AV opens when the pressure in the LV becomes larger than the one in the AA. This time marks the beginning of the ejection phase in systole. We report the volume rendering of velocity magnitude and the pressure on a clip in the LV apico-basal direction in Figures~\ref{clip_velocity} and \ref{clip_pressure}, respectively.  During the LV contraction, as shown in Figures~\ref{velocity.0001}-\ref{velocity.0003}, the blood flows from the LV to the AA, reaching large values of velocity at the systolic peak $t$ = 0.13 s. In particular, we measured a maximum flowrate in the AV section equal to 493.30 ml/s, a maximum pressure in the AA equal to 119.20 mmHg and a maximum ventricular pressure equal to 121.15 mmHg. The numerical results achieved are consistent with standard physiological data \cite{Gulsin_2017, carey2018prevention, sugimoto2017echocardiographic}: a peak systolic pressure in the range $119 \pm 13$ mmHg \cite{sugimoto2017echocardiographic}; a maximum flowrate of about 489 ml/s \cite{Gulsin_2017}.  A comparison between the biomarkers we compute from our numerical simulations and clinical values found literature is given in Table~\ref{table_comparison}.

\begin{table}[t]
	\centering
	\begin{tabular}{c|c|c|c}
		Biomarker & In-silico result & In-vivo measurements & Reference \\
		\hline
		LV stroke volume [ml] & 82.6 & $95 \pm 14$ & \cite{maceira2006normalized} \\
		LV ejection fraction [\%] & 55.8 & $57.5\pm 7.5$ &  \cite{kumar2014robbins} \\
		Peak AV flowrate  [ml/s] & $493.3$ &   $\approx 489$ & \cite{Gulsin_2017} \\
		LV peak pressure [mmHg] & 121.2 & $119 \pm 13$ & \cite{sugimoto2017echocardiographic} \\
		Peak E-wave velocity [m/s] & $0.96$ & $0.89 \pm 0.15$ & \cite{liza1998peak}\\
		Peak A-wave velocity [m/s] & $0.71$ & $0.78 \pm 0.26$ & \cite{liza1998peak}\\
		EA ratio $[-]$ & $1.35$ & $1.30 \pm 0.57$ & \cite{liza1998peak}
	\end{tabular}
	\caption{Biomarkers: comparison between numerical results and clinical values acquired in healthy individuals.}
	\label{table_comparison}
\end{table}

Specifically on the ejection phase, as shown in Figure \ref{av_velocity}, during the acceleration phase (Figures \ref{velocity_av.0001}, \ref{velocity_av.0002}), the spatial profile of the velocity field is almost flat, suggesting the development of a turbulent flow inside the AA. However, as the systolic peak is reached, the blood continues to flow from the LV to the AA (Figures \ref{velocity_av.0005}-\ref{velocity_av.0009}) but it decelerates. Differently from the acceleration stage, the velocity profile is no longer flat but the flow is partially oriented towards the LV. The AV closes when the flow becomes completely reverse on its section. Once the AV is closed, the pressure in the LV suddenly decreases, until it becomes smaller than the one in the LA: this marks the beginning of diastole with the opening of the MV.
The diastole is characterized by two filling stages: the E-wave and the A-wave.
As the E-wave starts ($t \approx 0.47$ s), a high-speed flux coming from the LA is observed at the MV section: the LV volume increases and the LA volume decreases. During the A-wave (atrial kick, $t \approx 0.92$ s), we observe a rapid contraction of the LA, as also observed in Figure \ref{volumes}, producing a second high-speed flux through the MV section, but milder than the one characterizing the E-wave. Once the atrial kick is over, and a reverse flow condition is detected on the MV section, the MV closes at the beginning of a new heartbeat.

Figure \ref{validation_MV_velocity} (left) shows the velocity magnitude $\bm u_\mathrm{MV}$ computed in our numerical simulations in a control volume immediately below the MV section, denoted as $\Omega_t^\mathrm{MV}$ in the Figure. We compare our result with the velocity profile acquired through trans-mitral valve spectral Doppler in a normal subject from \cite{Tripathi_2014} (normally acquired in a sample volume between the mitral leaflet tips \cite{appleton_chapter_10}). Specifically, we refer to Figure 4 of paper \cite{Tripathi_2014}.  We observe that our model is able to correctly reproduce amplitudes and shapes of the two characteristic waves in diastole, namely E-wave and A-wave: as for our reference data, the first wave shows a stronger peak than the second one. From the numerical results we measured maximum amplitudes of E-wave and A-wave equal to 0.96 m/s and 0.71 m/s, respectively;  an EA ratio equal to 1.35 (ratio among E-wave and A-wave peak velocities). Our results are consistent with physiological values commonly acquired in healthy subjects, as we summarize in Table~\ref{table_comparison}.

In Figure \ref{validation_MV_velocity} (right), we show the flowrate on a section downind the AV. We compare our result with the flowrate in a normal AA acquired through cardiac magnetic resonance through-plane phase-contrast velocity mapping. Specifically, we compare our numerical results with Figure 3 of paper \cite{Gulsin_2017}. By computing the flowrate at a section $\Gamma_t^*$ downwind the AV, we match accurately the systolic peak, with amplitude equal to 493.26 ml/s against a peak approximately equal to 490 ml/s in the reference data.

We wish to point out that our numerical simulations are run on a template -- even though realistic -- geometry \cite{zygote} and fed with data coming from EM simulations and 0D circulation model tuned for a generic healthy subject. Thus, from our analysis we can conclude that a qualitative good and satisfactory agreement can be found with the in-vivo results available in literature, making hence the whole computational model significant and reliable from an hemodynamic view point.

\begin{figure}[t!]
	\begin{subfigure}{0.190\textwidth}
		\centering
		\includegraphics[trim={1 1 1 1},clip, width=1.1\textwidth]{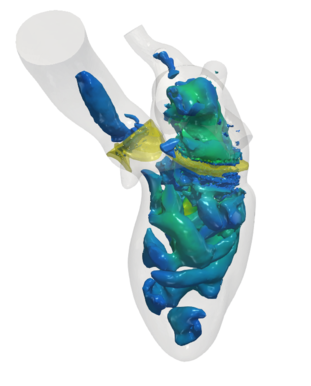}
		\caption{$t = 0.0$ s}
		\label{qcriterion_0000}
	\end{subfigure}
	\begin{subfigure}{0.190\textwidth}
		\centering
		\includegraphics[trim={1 1 1 1},clip, width=1.1\textwidth]{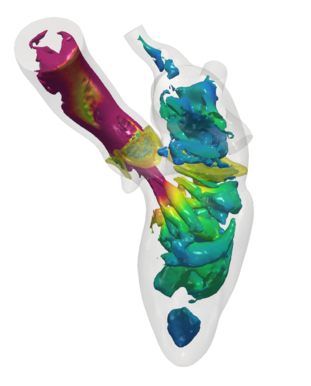}
		\caption{$t = 0.1$ s}
		\label{qcriterion_0001}
	\end{subfigure}
	\begin{subfigure}{0.190\textwidth}
		\centering
		\includegraphics[trim={1 1 1 1},clip, width=1.1\textwidth]{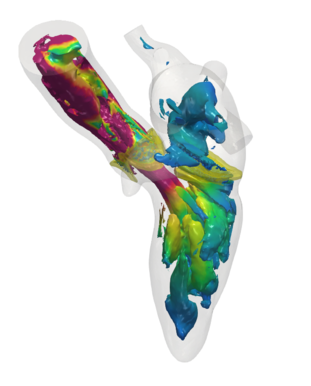}
		\caption{$t = 0.2$ s}
		\label{qcriterion_0002}
	\end{subfigure}
	\begin{subfigure}{0.190\textwidth}
		\centering
		\includegraphics[trim={1 1 1 1},clip, width=1.1\textwidth]{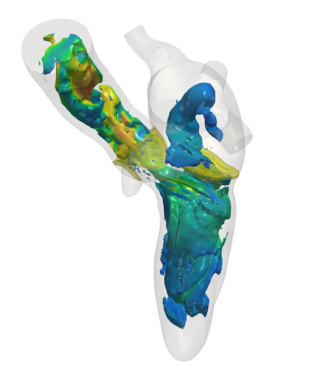}
		\caption{$t = 0.3$ s}
		\label{qcriterion_0003}
	\end{subfigure}
	\begin{subfigure}{0.190\textwidth}
		\centering
		\includegraphics[trim={1 1 1 1},clip, width=1.1\textwidth]{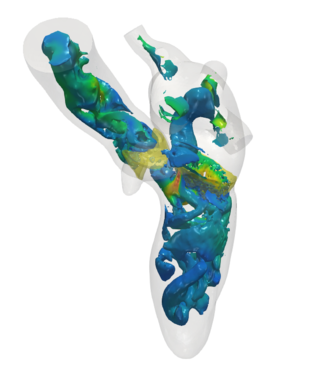}
		\caption{$t = 0.4$ s}
		\label{qcriterion_0004}
	\end{subfigure}
	\\
	\begin{subfigure}{0.190\textwidth}
		\centering
		\includegraphics[trim={1 1 1 1},clip, width=1.1\textwidth]{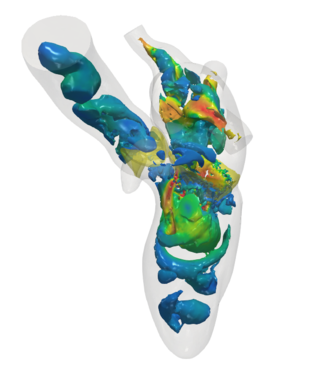}
		\caption{$t = 0.5$ s}
		\label{qcriterion_0005}
	\end{subfigure}
	\begin{subfigure}{0.190\textwidth}
		\centering
		\includegraphics[trim={1 1 1 1},clip, width=1.1\textwidth]{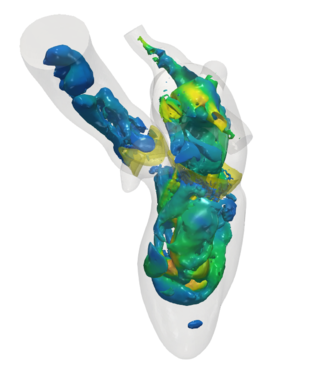}
		\caption{$t = 0.6$ s}
		\label{qcriterion_0006}
	\end{subfigure}
	\begin{subfigure}{0.190\textwidth}
		\centering
		\includegraphics[trim={1 1 1 1},clip, width=1.1\textwidth]{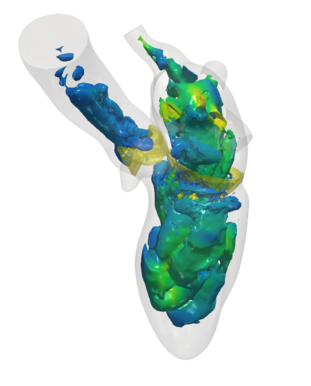}
		\caption{$t = 0.7$ s}
		\label{qcriterion_0007}
	\end{subfigure}
	\begin{subfigure}{0.190\textwidth}
		\centering
		\includegraphics[trim={1 1 1 1},clip, width=1.1\textwidth]{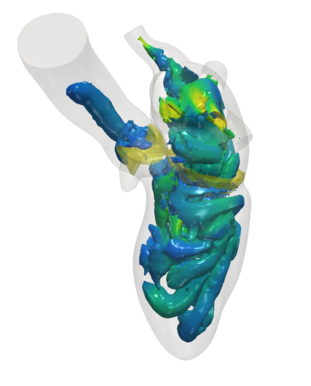}
		\caption{$t = 0.8$ s}
		\label{qcriterion_0008}
	\end{subfigure}
	\begin{subfigure}{0.1905\textwidth}
		\centering
		\includegraphics[trim={1 1 1 1},clip, width=1.1\textwidth]{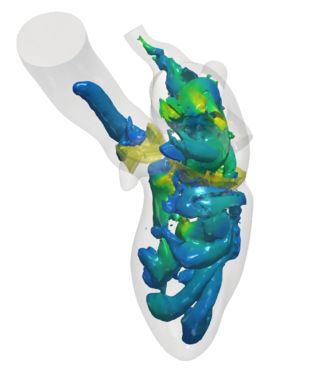}
		\caption{$t = 0.9$ s}
		\label{qcriterion_0009}
	\end{subfigure}
	\\
	\begin{subfigure}{\textwidth}
		\centering
		\includegraphics[trim={6 6 6 6},clip, width=0.2\textwidth]{images/velocity_scale.png}
	\end{subfigure}
	\caption{Iso-contours of Q-criterion ($\mathbb Q(\bm u) = 40 \, \mathrm{s}^{-2}$) colored according to velocity magnitude during a whole heartbeat.}
	\label{Q_criterion}
\end{figure}

\begin{figure}[!t]
	\centering
	\includegraphics[trim={4 4 4 4},clip, width=\textwidth]{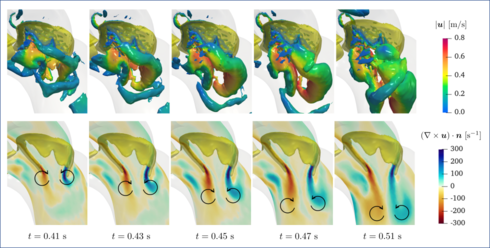}
	\caption{Formation of ring shaped vortex during early diastole. Top: iso-contours of Q-criterion (with $\mathbb Q(\bm u) =1000\,  \text{s}^{-2}$) coloured according to velocity magnitude. Bottom: projection of the vorticity on the normal direction (pointing towards the reader) of a slice in the LV apico-basal direction.}
	\label{O_vortex}
\end{figure}

To identify coherent structures, we introduce the scalar function \cite{q_criterion}:
\begin{equation*}
	\mathbb Q (\bm u) = \frac{1}{2} \left( |\bm \epsilon (\bm u)|_F^2 - |\bm \omega (\bm u)|_F^2 \right ),
\end{equation*}
where $\bm \epsilon (\bm u)$ is the strain rate tensor introduced in Section \ref{SEC:NS_ALE_RIIS}, $\bm \omega (\bm u) = \frac{1}{2}\left ( \nabla \bm u - \nabla ^T \bm u \right )$ the rotation tensor and $|\cdot|_F$ the Frobenius norm of a tensor. The Q-criterion consists of analysing the iso-contours of the positive part of $\mathbb Q(\bm u)$: when $\mathbb Q(\bm u)>0$, the rotation of a fluid is predominant with respect to its stretching. In Figure \ref{Q_criterion}, we show the iso-contours of Q-criterion with $\mathbb Q(\bm u) = 40 \, \mathrm{s}^{-2}$ at different instants during the whole heartbeat. At the beginning of the heartbeat (Figure \ref{qcriterion_0000}), we observe the residual vortical structures in the LH from late diastole. As the AV opens and the ejection phase starts, coherent structures are flushed out in the AA (Figures \ref{qcriterion_0001}-\ref{qcriterion_0003}). During systole, the LA is filled and coherent structures coming from the pulmonary veins impact in the middle of the LA, with still some visible structures at the end of systole. As the MV opens, we observe the formation of four vortex rings coming from the pulmonary veins in the LA. As in the systolic filling, they impact in the LA producing smaller coherent structures, during the E-wave and A-wave, as described in \cite{ZDMQ_2020}. At the same time, during E-wave, a big vortex ring rolls through the MV leaflets (Figure~\ref{qcriterion_0005}). The observed organized vortical pattern rapidly evolves into a chaotic complex flow that fills the whole LV reaching its apex (Figures~\ref{qcriterion_0006}-\ref{qcriterion_0008}). During the A-wave, a new vortex downstream the MV section is formed, but weaker than the one seen during the E-wave (Figure~\ref{qcriterion_0009}). The MV closes, the vortex under its section is suddenly broken and, with the opening of the AV, the new cycle begins and all the coherent structures are flushed out again.

\begin{figure}[t!]
	\centering
	\begin{subfigure}{0.49\textwidth}
		\centering
		\includegraphics[trim={4 4 4 4},clip, width=\textwidth]{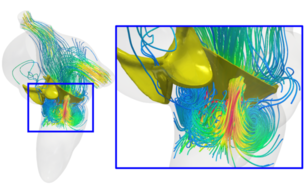}
		\caption{Streamlines}
		\label{streamline}
	\end{subfigure}
	\begin{subfigure}{0.49\textwidth}
		\centering
		\includegraphics[trim={4 4 4 4},clip, width=\textwidth]{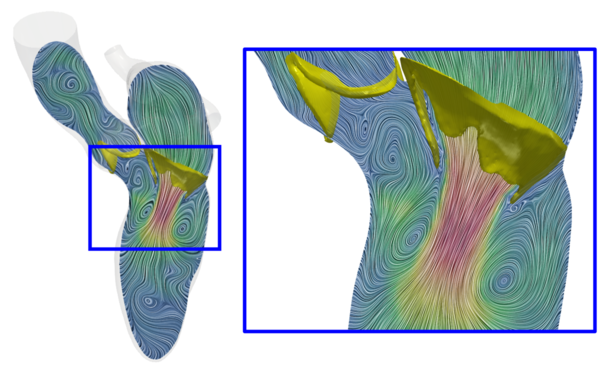}
		\caption{Surface LIC}
		\label{surface-lic}
	\end{subfigure}
	\\
	\begin{subfigure}{\textwidth}
		\centering
		\includegraphics[trim={6 6 6 6},clip, width=0.2\textwidth]{images/velocity_scale.png}
	\end{subfigure}
	\caption{Vortex formation under the MV section at early diastolic peak: (A) Streamlines colored according to velocity magnitude ($t=0.45$ s); (B) Surface LIC visualization -- with velocity as integrator -- on a slice colored according to velocity magnitude ($t=0.48$ s).}
	\label{LH_streamline_4DflowMRI}
\end{figure}

In each subfigure of Figure \ref{O_vortex},  we show the iso-contours of Q-criterion with $\mathbb Q(\bm u)=1000 \, \mathrm{ s}^{-2}$ on the top; the projection of the vorticity on the normal direction of a slice in the LV apico-basal direction on the bottom. Specifically, we focus on the region immediately under the MV section to better asses the formation of vortical structures during diastole. We observe the formation of shear layers on the leaflets of the MV producing different velocities on the two sides of the leaflets. Shear layers roll on MV leaflets and enter in the main cavity producing low-pressure circulation regions. The flow rotates in clockwise direction under the anterior leaflet, counterclockwise under the posterior, as shown in the projection of the vorticity in the slice's normal direction. This process characterizes the formation of the vortex ring that rolls through the MV leaflets (also shown in Figure \ref{Q_criterion}), with a larger velocity in the inner part with respect to the exterior part of the ring. Once the vortex is formed, it evolves towards the LV apex and it breaks into smaller coherent structures. Figure \ref{LH_streamline_4DflowMRI} shows the vortex formation under the MV section via streamlines colored according to the velocity magnitude and a surface Line Integral Convolution (LIC) visualization. We compare our numerical results with 4D flow MRI visualizations available in Figure 3C of \cite{dyverfeldt_2015}, and Figure 4 of \cite{van_der_geest_2016}, respectively. The numerical simulation well predicts the high speed jet across the MV and the formation of vortical flow under its section described by the 4D flow MRI data during the E-wave peak.
Thus, our model captures the formation and the dynamics of the vortex ring (also referred to as ``O vortex'' \cite{DMQ_2019, Chnafa_2014}) during diastole, a well studied cardiac hemodynamic feature whose characteristics and interactions with the LV wall provide information about the diastolic function \cite{Ardvisson_2016, Kim_1995, Chnafa_2014, Mittal_review_2016}.

\section{Conclusions and future developments}
\label{SEC:conclusions}
We proposed a computational model for the  assessment of the left heart hemodynamics. Our model accounted for the displacement of the domain boundary, the motion of the aortic and mitral valve, the dynamics of the circulatory system, and transition to turbulence effects. We carried out an electromechanical simulation of the left ventricle and, for the description of the motion over the whole domain boundary, we introduced an extension procedure combined with a volume-based reduced model for the atrium. The 3D CFD model of the left heart consisted of the Navier-Stokes equations in an Arbitrary Lagrangian Eulerian framework, with valves accounted for by means of the Resistive Immersed Implicit Surface method. We coupled the 3D CFD model to a closed-loop 0D circulation model of the whole cardiovascular system.
We simulated the blood flow in the left heart in physiological conditions by means of the proposed methodology. The results were analyzed in terms of flow and pressure distribution, velocity profiles through the valves, and turbulent coherent structures.
Biomarkers such as the stroke volume, the ejection fraction, the EA ratio and the pressure peaks were compared with those available in the literature and accurately reproduced; a comparison with Doppler echocardiography data and 4D flow magnetic resonance imaging allowed to validate the computational model.

The novelty of the paper is threefold.
\begin{enumerate}
	\item We introduced an original preprocessing procedure that combines i) an extension of a given left ventricular displacement on the whole left heart by means of Laplace-Beltrami equations with physiological kinematic constraints; ii) a reduced model for the motion of the atrium based on the volume variation dictated by a lumped-parameter circulation model. This yields an integrated system in which fluid dynamics is one-way coupled to electromechanics in the left ventricle. The extension procedure can also be employed to merge displacement fields coming from different sources, such as a reconstruction from diagnostic images, and may, thus, be applied also to patient-specific studies. As the imaging data routinely acquired in diagnostic exams is often mainly focused on the ventricle, the extension procedure that we proposed can be used to complete the missing data. \item We devised a coupled 3D-0D model made of the 3D CFD model of the left heart and a 0D circulation model of the whole cardiovascular system. We solved the coupled model with a segregated scheme and we developed computational strategies to solve the integrated system made of fluid dynamics, displacement, valves and circulation models.
	\item We found that our numerical simulations yielded a qualitative and quantitative good agreement with clinical data from different sources, making the whole integrated multiscale model significant and realiable from an hemodynamic view point.
\end{enumerate}

A number of further directions of investigation will be considered. Introducing an electromechanical model of the complete left heart would complete the description of the cardiac contractility and hence allow the relaxation of the hypotheses we made on the atrium motion. In addition, in the present work, we prescribed the electromechanics-based displacement on the boundary of the fluid dynamics system, hence enforcing a kinematic coupling condition between the two physics. Since we are neglecting the dynamic balance (that is the equilibration of forces), the pressure would not be well defined when both valves are closed. Thus, we neglect the isovolumetric phases of the heart cycle. As done in \cite{fernandez_2019}, one could circumvent this issue by introducing an additional penalty term in the Navier-Stokes equations to control the intraventricular pressure during the isovolumetric phases. 

To conclude, the present work can be viewed as a significant step forward to simulate the blood flow dynamics of the whole human heart, allowing to have a complete description of the hemodynamics of the whole organ. This represents a challenging task from a modeling point of view, that we aim to reach with intermediate stages.

\section*{Acknowledgments}
This work has been supported by the ERC Advanced Grant iHEART, ``An Integrated Heart Model for the simulation of the cardiac function'', 2017–2022,  P.I. A. Quarteroni (ERC–2016– ADG, project ID: 740132).
We gratefully acknowledge the CINECA award under the ISCRA B initiative, for the availability of high performance computing resources and support under the project HP10BD303V, P.I. A. Quarteroni, 2020-2021

\clearpage


\appendix

\section{Setup of the circulation model}
\label{append_circulation}
We detail the setup of the 0D circulation model we use for the multiscale 3D-0D CFD simulation. Specifically, we report the parameters employed in Table~\ref{parameters_circulation} and the initial state variables in Table~\ref{initial_state_circulation}. The external pressure $p_\mathrm{EX}(t)$ is set to 0.
\begin{table}[h!]
	\centering
	\begin{tabular}{c|c|c|c|c}
	Compartment & Parameter  & Description & Unit of measure  & Value
	\\
	\hline
	& $E_\mathrm{A}$	&  Active elastance & [mmHg/ml] & 0.06 \\
	& $E_\mathrm{B}$ &  Passive elastance	& [mmHg/ml] & 0.07 \\
	Right & $d_\mathrm{c}$ & Duration of contract. & relative w.r.t. $T_\mathrm{HB}$ & 0.335 \\
	atrium & $d_\mathrm{r}$ & Duration of relax. & relative w.r.t. $T_\mathrm{HB}$ & $1.45 \cdot 10^{-2}$ \\
	& $t_\mathrm{c}$ & Initial time of contract. & relative w.r.t. $T_\mathrm{HB}$ &  0.80 \\
	& $V_{0, \mathrm{RA}}$ & Resting volume & [ml] & 4.00\\
	\hline
	& $E_\mathrm{A}$	&  Active elastance & [mmHg/ml] & 0.65 \\
	& $E_\mathrm{B}$ &  Passive elastance	& [mmHg/ml] & 0.05 \\
	Right & $d_\mathrm{c}$ & Duration of contract. & relative w.r.t. $T_\mathrm{HB}$ & 0.335 \\
	ventricle & $d_\mathrm{r}$ & Duration of relax. & relative w.r.t. $T_\mathrm{HB}$ & $1.45 \cdot 10^{-2}$ \\
	& $t_\mathrm{c}$ & Initial time of contract. & relative w.r.t. $T_\mathrm{HB}$ &  0.00 \\
	& $V_{0, \mathrm{RV}}$ & Resting volume & [ml] & 10.00\\
	\hline
	Pulmonary & $R^\mathrm{PUL}_\mathrm{AR}$ &  Resistance & $[\mathrm{mmHg}\cdot\mathrm{s}/\mathrm{ml}]$ & 0.25 \\
	arterial & $C^\mathrm{PUL}_\mathrm{AR}$ &  Capacitance & $[\mathrm{ml}/\mathrm{mmHg}]$ & 5.00 \\
	system & $L^\mathrm{PUL}_\mathrm{AR}$ & Inductance & $[\mathrm{mmHg}\cdot \mathrm{s}^2/\mathrm{ml}]$ & $5.00 \cdot 10^{-4}$ \\
	\hline
	Pulmonary & $R^\mathrm{PUL}_\mathrm{VEN}$ &  Resistance & $[\mathrm{mmHg}\cdot\mathrm{s}/\mathrm{ml}]$ & 0.02\\
	venous & $C^\mathrm{PUL}_\mathrm{VEN}$ &  Capacitance & $[\mathrm{ml}/\mathrm{mmHg}]$ & 100.00 \\
	system & $L^\mathrm{PUL}_\mathrm{VEN}$ & Inductance & $[\mathrm{mmHg}\cdot \mathrm{s}^2/\mathrm{ml}]$ & $5.00 \cdot 10^{-5}$ \\
	\hline
	Systemic & $R^\mathrm{SYS}_\mathrm{AR}$ &  Resistance & $[\mathrm{mmHg}\cdot\mathrm{s}/\mathrm{ml}]$ & 1.00 \\
	arterial & $C^\mathrm{SYS}_\mathrm{AR}$ &  Capacitance & $[\mathrm{ml}/\mathrm{mmHg}]$ & 2.00 \\
	system & $L^\mathrm{SYS}_\mathrm{AR}$ & Inductance & $[\mathrm{mmHg}\cdot \mathrm{s}^2/\mathrm{ml}]$ & $5.00 \cdot 10^{-3}$ \\
	\hline
	Systemic & $R^\mathrm{SYS}_\mathrm{VEN}$ &  Resistance & $[\mathrm{mmHg}\cdot\mathrm{s}/\mathrm{ml}]$ & 0.24 \\
	venous & $C^\mathrm{SYS}_\mathrm{VEN}$ &  Capacitance & $[\mathrm{ml}/\mathrm{mmHg}]$ & 60.00 \\
	system & $L^\mathrm{SYS}_\mathrm{VEN}$ & Inductance & $[\mathrm{mmHg}\cdot \mathrm{s}^2/\mathrm{ml}]$ & $5.00 \cdot 10^{-4}$ \\
	\hline
	Tricuspid & $R_\mathrm{min}$ & Minimum resistance & [$\mathrm{mmHg}\cdot \mathrm{s} / \mathrm{ml}$] & $7.5 \cdot 10^{-3}$ \\
	valve & $R_\mathrm{max}$ & Maximum resistance & [$\mathrm{mmHg}\cdot \mathrm{s} / \mathrm{ml}$] & $75006.2$ \\
	\hline
	Pulmonary & $R_\mathrm{min}$ & Minimum resistance & [$\mathrm{mmHg}\cdot \mathrm{s} / \mathrm{ml}$] & $7.5 \cdot 10^{-3}$ \\
	valve & $R_\mathrm{max}$ & Maximum resistance & [$\mathrm{mmHg}\cdot \mathrm{s} / \mathrm{ml}$] & $75006.2$ \\
	\end{tabular}
\caption{Parameters used in the circulation model.}
\label{parameters_circulation}
	\end{table}

\begin{table}[t]
	\begin{tabular}{c|c|c|c|c}
	Compartment & Parameter  & Description & Unit of measure  & Value \\
	\hline
	Right atrium & $V_\mathrm{RA}$ & Volume & [ml] & 78.95\\
	\hline
	Right ventricle & $V_\mathrm{RV}$ & Volume & [ml] & 154.00 \\
	\hline
	Pulmonary arterial  & $p^\mathrm{PUL}_\mathrm{AR}$ & Pressure & [mmHg] & 33.50 \\
	system &  $Q^\mathrm{PUL}_\mathrm{AR}$ & Flowrate & [ml/s] & 69.44 \\
		\hline
	Pulmonary venous  & $p^\mathrm{PUL}_\mathrm{VEN}$ & Pressure & [mmHg] & 16.16 \\
	system &  $Q^\mathrm{PUL}_\mathrm{VEN}$ & Flowrate & [ml/s] & 0.00 \\
	\hline
	Systemic arterial  & $p^\mathrm{SYS}_\mathrm{AR}$ & Pressure & [mmHg] & 91.68\\
	system &  $Q^\mathrm{SYS}_\mathrm{AR}$ & Flowrate & [ml/s] & 63.71 \\
	\hline
	Systemic venous  & $p^\mathrm{SYS}_\mathrm{VEN}$ & Pressure & [mmHg] & 23.99\\
	system &  $Q^\mathrm{SYS}_\mathrm{VEN}$ & Flowrate & [ml/s] & 65.40 \\
	\hline
	\end{tabular}
	\caption{Initial state of the circulation model.}
	\label{initial_state_circulation}
\end{table}

\section{The ventricular EM model}
\label{appendix_em}
To carry out EM simulations, we employ the EM model of the left ventricle proposed in \cite{REGAZZONI2022111083}. Electrophysiology is modeled through the Monodomain equation \cite{franzone2014mathematical} coupled with the ten Tusscher-Panfilov ionic model \cite{ten2006alternans}. Subcellular generation of active force is modeled through the activation model proposed in \cite{regazzoni2020machine}. The passive behaviour of the tissue is modeled via the Guccione strain energy density function \cite{guccione1991finite}. We set generalized Robin BCs at the epicardium and we impose energy-consistent BCs \cite{regazzoni2020machine} on the rings of the AV and MV.
The fibers distribution is generated by the rule-based Bayer-Blake-Plank-Trayanova algorithm \cite{bayer2012novel, piersanti2021modeling}. The 3D EM model of the LV is coupled to the 0D closed-loop circulation model of the whole cardiovascular system introduced in Section \ref{SEC:0D_CIRCULATION}. For the numerical approximation of the problem, we use linear finite elements (FE) for the space discretization. Specifically, in the EM model an intergrid transfer operator \cite{SDQ_2020} is adopted to employ a coarser grid for the elastodynamic and a finer one for the electrophysiology, due to the higher resolution required by the latter.  The coupled EM problem is solved by means of the staggered numerical scheme presented in \cite{REGAZZONI2022111083}.

\end{document}